\documentclass[preprint,3p,fleqn,sort&compress,numbers]{elsarticle}
\usepackage{amsthm,amsmath,amsfonts,amssymb}
\usepackage{booktabs}
\usepackage{multirow}
\usepackage{mathtools}
\usepackage{esint}
\usepackage{graphicx}
\usepackage{colortbl}
\usepackage{empheq}
\usepackage{bm}
\usepackage{wasysym}
\usepackage[usenames,dvipsnames]{xcolor}
\usepackage{caption,subcaption}
\usepackage{dsfont}
\usepackage[toc,page]{appendix}
\usepackage[colorlinks]{hyperref}
\usepackage{cleveref}

\newdefinition{remark}{Remark}

\renewcommand{\u}{{\mathbf{u}}}
\newcommand{\n}{{\mathbf{n}}}
\newcommand{\f}{{\mathbf{f}}}
\renewcommand{\b}{{\mathbf{b}}}
\newcommand{\U}{{(\u,p)}}
\newcommand{\G}{{\mathbf{G}}}
\newcommand{\F}{{\mathbf{F}}}
\newcommand{\bxi}{{\boldsymbol{\xi}}}
\newcommand{\grad}{\nabla}
\newcommand{\IBa}{{$\textnormal{IB}_{\sqrt{\cdot}}$}}
\newcommand{\IBb}{{$\textnormal{IB}_{\textnormal{4}}$}}
\definecolor{review_color}{rgb}{0.0, 0.0, 0.0}

\crefname{appsec}{Appendix}{Appendices}

\graphicspath{{./fig/}}

\begin{document}

\definecolor{py_green}{rgb}{0.0, 0.5019607843137255, 0.0}
\definecolor{py_blue}{rgb}{0.0, 0.0, 1.0}
\definecolor{py_purple}{rgb}{0.5019607843137255, 0.0, 0.5019607843137255}
\definecolor{py_orange}{rgb}{1.0, 0.6470588235294118, 0.0}


\title{Immersed Boundary Smooth Extension (IBSE): A high-order method for solving incompressible flows in arbitrary smooth domains}

\author[ucd]{David~B.~Stein\corref{cor1}}
\ead{dbstein@math.ucdavis.edu}

\author[ucd]{Robert~D.~Guy}

\author[ucd]{Becca~Thomases}

\address[ucd]{Department of Mathematics, University of California, Davis, Davis, CA 95616-5270, USA}
\cortext[cor1]{Corresponding author}

\begin{abstract}
	The Immersed Boundary method is a simple, efficient, and robust numerical scheme for solving PDE in general domains, yet for fluid problems it only achieves first-order spatial accuracy near embedded boundaries for the velocity field and fails to converge pointwise for elements of the stress tensor. In a previous work we introduced the Immersed Boundary Smooth Extension (IBSE) method, a variation of the IB method that achieves high-order accuracy for elliptic PDE by smoothly extending the unknown solution of the PDE from a given smooth domain to a larger computational domain, enabling the use of simple Cartesian-grid discretizations. In this work, we extend the IBSE method to allow for the imposition of a divergence constraint, and demonstrate high-order convergence for the Stokes and incompressible Navier-Stokes equations: up to third-order pointwise convergence for the velocity field, and second-order pointwise convergence for all elements of the stress tensor. The method is flexible to the underlying discretization: we demonstrate solutions produced using both a Fourier spectral discretization and a standard second-order finite-difference discretization.
\end{abstract}

\begin{keyword}
	Embedded boundary, Immersed Boundary, Incompressible Navier-Stokes, Fourier spectral method, Complex geometry, High-order
\end{keyword}

\maketitle



\section{Introduction}

The Immersed Boundary (IB) method was originally developed for the study of moving, deformable structures immersed in a fluid, and it has been widely applied to such problems since its introduction \cite{Peskin2002,Mittal2005,olson2014simulating}.  Recently, the method has been adapted to more general fluid-structure problems, including the motion of rigid bodies immersed in a fluid \cite{kallemov2016immersed} and fluid flow through a domain with either stationary boundaries or boundaries with prescribed motion \cite{Taira2007,Teran2009}.  In this broadened context, we use the term \emph{Immersed Boundary method} to refer only to methods in which ({\small\emph{i}}) the boundary is treated as a Lagrangian structure embedded in a geometrically simple domain, ({\small\emph{ii}}) the background PDE (e.g. the Navier-Stokes equations) are solved on a Cartesian grid everywhere in that domain, and ({\small\emph{iii}}) all communication between the Lagrangian structure and the underlying PDE {\color{review_color}is} mediated by convolutions with regularized $\delta$-functions.  These methods have many desirable properties: they make use of robust and efficient Cartesian-grid methods for solving the underlying PDE, are flexible to a wide range of problems, and are simple to implement, requiring minimal geometric information and processing to describe the boundary.

The IB method belongs to the broad category of methods known as \emph{embedded boundary} (EB) methods, including the Immersed Interface \cite{li2006immersed}, Ghost Fluid \cite{Fedkiw1999}, and Volume Penalty methods \cite{Angot1999}.  These methods share a common feature: they enable solutions to PDE on nontrivial domains to be computed using efficient and robust structured-grid discretizations; yet these methods differ largely in how boundary conditions are enforced and whether or not the solution is produced in the entirety of a simple domain.  Methods which compute the solution everywhere in a $d$-dimensional rectangle admit the simplest discretizations and enable the use of high-order discretizations such as Fourier spectral methods.  Unfortunately, this simplicity comes coupled with a fundamental difficulty: the \emph{analytic} solution to these problems is rarely globally smooth on the entire domain.  Consider the one-dimensional Poisson problem $\Delta u=f$ on the periodic interval $\mathbb{T}=[0,2\pi]$ with Dirichlet boundary conditions $u(a)=u(b)=0$ for $a\neq b\in\mathbb{T}$.  Even if $f\in C^\infty(\mathbb{T})$, the solution $u$ will typically display jumps in its derivative at the values $x=a$ and $x=b$.  The lack of regularity in the analytic problem leads to low-order convergence in many numerical schemes, including the Immersed Boundary method: the addition of (regularized) singular forces supported at the boundary causes the solution to be $C^0$, and solutions are accurate only to first-order in the grid spacing $\Delta x$ \cite{kallemov2016immersed,Taira2007}.

The advantages of EB methods are substantial enough that significant effort has been expended on improving their accuracy \cite{Lai2000,Mark2008,Linnick2005,Liu2014,Xu2006,Zhong2007,Yu2007,Zhou2006,Gibou2005,Boyd2005,Bueno-Orovio2006,Lui2009,Sabetghadam2009,Albin2011,Lyon2010a,Lyon2010,Shirokoff2013}.  Two different approaches are generally taken.  The first approach involves locally altering the discretization of the PDE in the vicinity of the boundary to accommodate the lack of smoothness in the solution.  One example of this approach is the Immersed Interface method \cite{li2006immersed}.  Such approaches are particularly useful for interface problems where the solution is required on both sides of the embedded boundary.  When the solution is only needed on one side of the interface, a second approach may be taken in which variables are redefined outside of the domain of interest to obtain higher global regularity.  Improved convergence rates are achieved as a natural consequence of the properties of the discretization scheme when applied to smooth problems.  Variations of this basic idea have been used by the Fourier Continuation (FC) \cite{Lyon2010a,Lyon2010,Albin2011} and the Active Penalty (AP) \cite{Shirokoff2013} methods to provide high order solutions to PDE on general domains.

In \cite{Stein2015}, we introduce another approach: the \emph{Immersed Boundary Smooth Extension} (IBSE) method for the solution of the Poisson and related problems (e.g. the heat equation and Burgers equation). The IBSE method builds off of the basic framework of the direct forcing IB method. Drawing on ideas from the AP and FC methods, the IBSE method remedies the slow convergence of the direct forcing IB method by solving an auxiliary problem to extend the \emph{unknown} solution from the physical domain into the entire computational domain. This extension is used to define the body forcing in the non-physical portion of the domain, forcing the solution to be \emph{globally} $C^k$. High-order accuracy is then achieved naturally, up to $\mathcal{O}(\Delta x^{k+1})$ or the maximum accuracy supported by the underlying discretization for a smooth problem.

In this paper, we extend the work from \cite{Stein2015} to apply to PDE requiring the imposition of a divergence constraint, including the Stokes and Navier-Stokes equations. This is a non-trivial difficulty for a high-order embedded boundary scheme based on smooth extensions. Prior work has either been restricted to the compressible Navier-Stokes equations \cite{Albin2011}, or has dealt with the pressure in an ad-hoc manner, limiting the overall accuracy of the numerical scheme to second-order \cite{Shirokoff2013}. These methods either extend the forcing function or solution from a prior timestep, and are thus unable to provide solutions to the Stokes problem. For time-dependent problems, they require either explicit timestepping or the use of Alternating-Direction-Implicit (ADI) schemes to take implicit timesteps, complicating the ability to achieve high-order accuracy in time.

In contrast, we work directly with the Stokes equations, rather than employing dimensional splitting or decoupling the problem into a set of Poisson equations. The smooth extensions to the solution of the Stokes equation define forcing functions in the non-physical domain in a coupled manner, and the Lagrange multipliers supported on the boundary that define these extensions are fully coupled through the Stokes problem, eliminating any need to impose artificial boundary conditions for a pressure Poisson equation. The IBSE method smoothly extends \emph{the unknown solution}, allowing for efficient direct solutions of the steady incompressible Stokes equation and simple high-order in time implicit-explicit time discretization of the Navier-Stokes equations. Remarkably, the fully coupled problem for the solution to the incompressible Stokes equations (combined with the equations that define the smooth extensions to the velocity and pressure fields) can be reduced to the solution of a relatively small\footnote{With size a small multiple of the number of points used to discretize the boundary} dense system of equations, two Stokes solves on a simple domain ignoring the internal boundary, and a handful of fast Fourier transforms. The dense system depends only on the boundary and the discretization, and so it can be formed and prefactored to allow for efficient time-stepping.

The IBSE method retains the essential robustness and simplicity of the original IB method.  All communication between the Lagrangian boundary and the underlying Cartesian grid is achieved by convolution with regularized $\delta$-functions or normal derivatives of those $\delta$-functions.  This allows an absolute minimum of geometric information to be used.  In the traditional IB method, only the position of the Lagrangian structure must be known; the IBSE method additionally requires normals to that structure and an indicator variable denoting whether Cartesian grid points lie inside or outside of the physical domain where the PDE is defined.  As with the IB method, the IBSE method is flexible to the choice of the underlying PDE discretization, so long as the mesh is uniform in the vicinity of the immersed boundary, facilitating simple coupling to existing code. In this paper, we demonstrate the IBSE solver with an underlying Fourier spectral discretization as well as a second-order staggered-grid finite difference discretization.

This paper is organized as follows. In \Cref{stokessection:methods}, we introduce the methodology, ignoring the details of the numerical implementation. This section includes a review of the relevant methodology introduced in \cite{Stein2015}. In \Cref{stokessection:numerical_implementation}, we discuss the numerical implementation of the method, not including the choice of the underlying discretization of the PDE. In \Cref{stokessection:stokes_test_analytic}, we study a steady Stokes problem with an analytic solution, using a Fourier spectral discretization. We analyze the convergence of the velocity and pressure fields, as well as derivatives of the velocity fields. We demonstrate up to third-order pointwise convergence of the velocity fields, as well as second-order pointwise convergence of all elements of the stress tensor. In \Cref{stokessection:stokes_test:flow_around_cylinder}, we analyze the flow of a viscous incompressible fluid around a cylinder in a confined channel at zero Reynolds number. We solve the problem using both a Fourier spectral discretization and a second-order finite difference discretization, comparing the accuracy and convergence of the solution, as well as the convergence of the scalar drag coefficient to known values from the literature. Finally, in \Cref{stokessection:navier_stokes}, we solve the incompressible Navier-Stokes equations in both steady and unsteady cases, for the flow around a cylinder in a confined channel, demonstrating rapid convergence and quantitative agreement with a number of benchmarks from the literature.



\section{Methods}
\label{stokessection:methods}

We begin by considering the time dependent, incompressible Navier-Stokes equations:
\begin{subequations}
	\begin{align}
		\partial_t\u + \u\cdot\grad\u - \nu\Delta\u + \grad p	&=	\f_\u	&&\text{in }\Omega,	\\
		\grad\cdot\u	&=	0		&&\text{in }\Omega,	\\
		\u						&=	\b	&&\text{on }\Gamma,
	\end{align}
\end{subequations}
where the boundary of the domain $\Omega$ is denoted by $\Gamma=\partial\Omega$. One way to discretize these equations in time is to treat the non-linearity explicitly, and to treat the Stokes operator implicitly. The simplest such time discretization, using forward Euler for the nonlinear terms and backward Euler for the remaining terms, is
\begin{subequations}
	\label{stokeseq:navier_stokes_time_discrete}
	\begin{align}
		(\mathbb{I} - \nu\Delta t\Delta)\u^{t+\Delta t} + \Delta t\grad p^{t+\Delta t}	&=	\u^t + \Delta t\left(\f_\u^{t+\Delta t} - \u^t\cdot\grad\u^t\right)	&&\text{in }\Omega,	\\
		\grad\cdot\u^{t+\Delta t}	&=	0	&&\text{in }\Omega,	\\
		\u^{t+\Delta t}	&=	\b^{t+\Delta t}	&&\text{on }\Gamma.
	\end{align}
\end{subequations}
Although the Navier-Stokes equations are often discretized in time using projection methods \cite{chorin1968numerical}, these discretizations produce large splitting errors near the boundary at low Reynolds number. Implicit-Explicit time-discretizations of the form of \Cref{stokeseq:navier_stokes_time_discrete} eliminate these errors, and are simpler to generalize to higher-order in time schemes (see \Cref{stokessection:navier_stokes}). The unsplit Stokes system may be inverted efficiently using GMRES with a projection preconditioner \cite{cai2014efficient}. We find this inversion strategy to be very efficient (convergence to an algebraic tolerance of $10^{-10}$ achieved in 3-4 iterations) for moderate to high Reynolds numbers ($\textnormal{Re}\geq 1$), and moderately efficient (convergence achieved in approximately 35 iterations) for the Stokes problem where projection methods cannot be used.

Motivated by the form of \Cref{stokeseq:navier_stokes_time_discrete}, in this section we focus primarily on developing the methodology required to solve the general Stokes problem:
\begin{subequations}
	\label{eq:original_stokes}
	\begin{align}
		(\alpha\mathbb{I}-\Delta)\u + \grad p	&=	\f_\u	&&\text{in }\Omega,	\\
		\grad\cdot\u				&=	f_p		&&\text{in }\Omega,	\\
		\u									&=	\b			&&\text{on }\Gamma.
	\end{align}
\end{subequations}
When $\alpha=0$ and $f_p=0$, \Cref{eq:original_stokes} reduces to the steady incompressible Stokes equations. The case where $\alpha>0$ arises in the study of porous media \cite{brinkman1949calculation} and from discretizing the Navier-Stokes equations in time using an implicit-explicit discretization, as in \Cref{stokeseq:navier_stokes_time_discrete}. The non-zero divergence constraint ($f_p\neq 0$) adds no additional complexity and is included for completeness. Before discussing the IBSE method as applied to \Cref{eq:original_stokes}, we first review relevant material from \cite{Stein2015}, in order to introduce the central ideas in the simplified context of the Poisson problem.

\subsection{Review of \cite{Stein2015}: the IBSE method for the Poisson problem}

\begin{figure}
	\centering
	\hspace*{\fill}
	\begin{subfigure}[b]{0.3\textwidth}
		\centering
		\includegraphics[width=\textwidth]{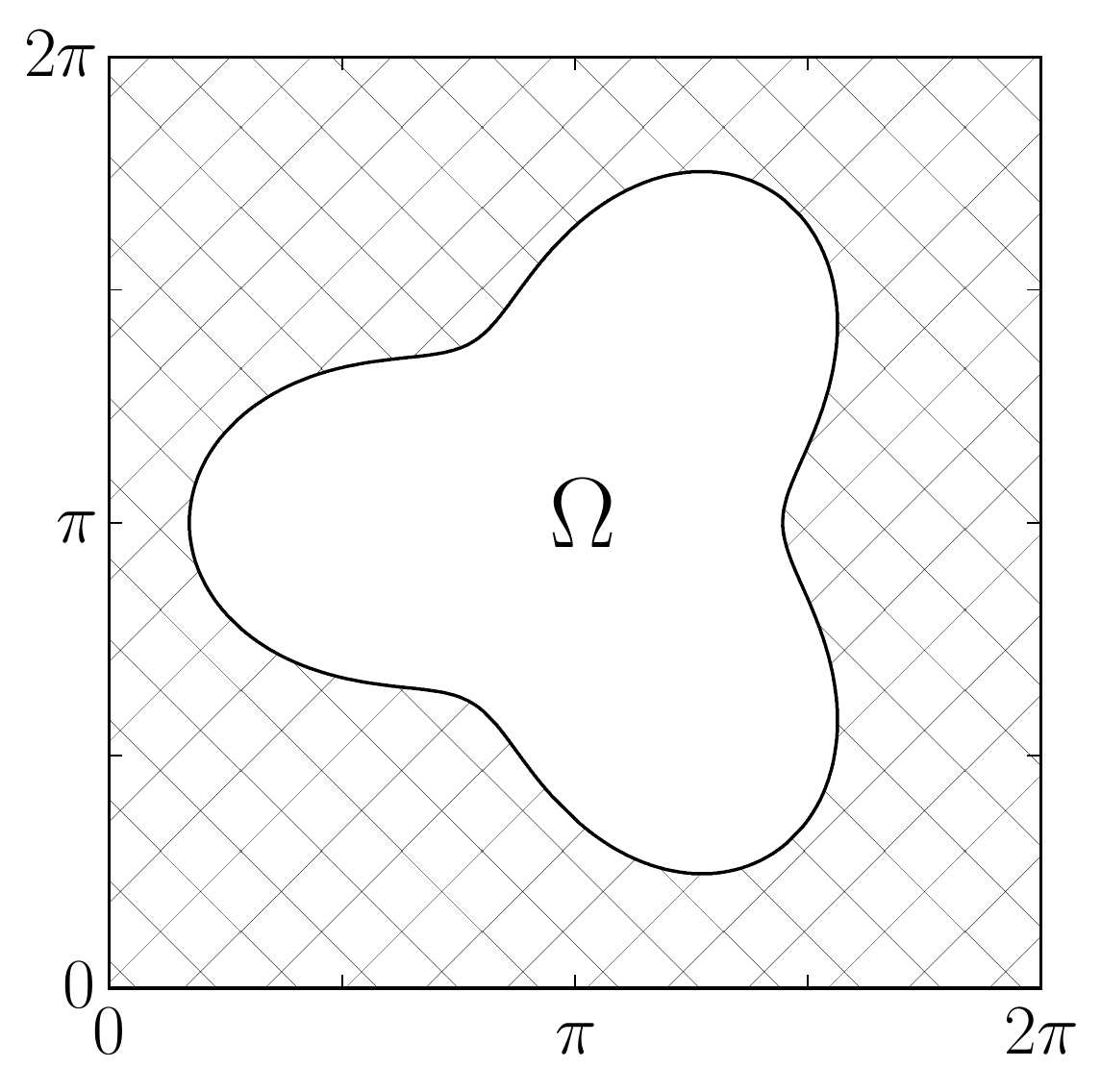}
		\subcaption{Irregular domain}
		\label{stokesfig:domain:1}
	\end{subfigure}
	\hfill
	\begin{subfigure}[b]{0.53\textwidth}
		\centering
		\includegraphics[width=\textwidth]{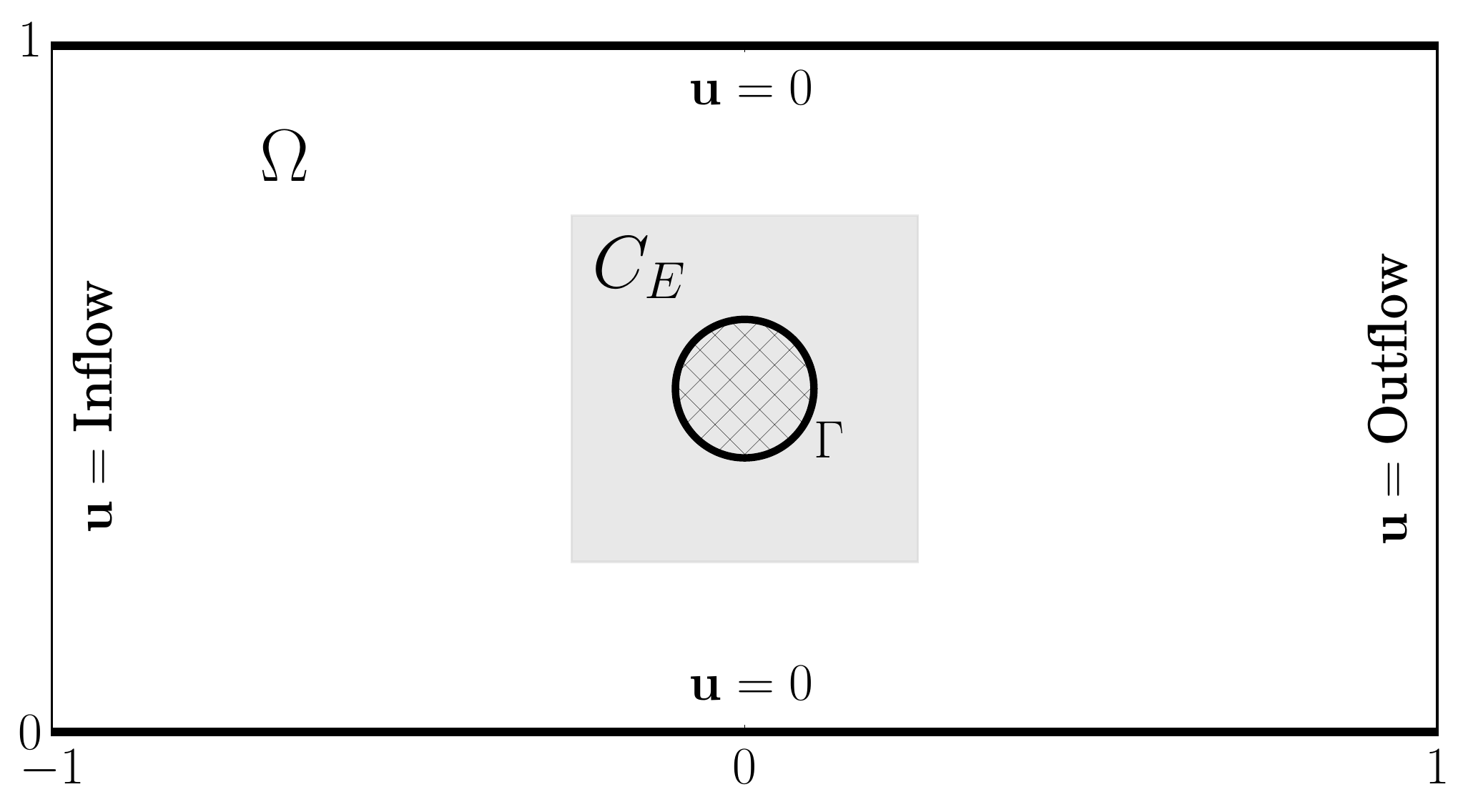}
		\subcaption{Flow around cylinder}
		\label{stokesfig:domain:2}
	\end{subfigure}
	\hspace*{\fill}
	\caption[Two examples of complicated domains for solving partial differential equations.]{Two different and typical domains to solve PDE on.  The physical domain $\Omega$ is shown in white, the extension domain $E$ is denoted by crosshatches. \Cref{stokesfig:domain:1} shows an irregular interior domain. The computational domain $C$ is the 2-torus, $\mathbb{T}^2=[0,2\pi]\times[0,2\pi]$. For this domain it is convenient to take $C_E=C$. \Cref{stokesfig:domain:2} shows a domain that would be used to compute flow around a cylinder. The computational domain $C$ is taken to be $[-1,1]\times[0,1]$, with Dirichlet boundary conditions specified along $\partial C$ (homogeneous for the channel walls at the top and bottom, and specified inflow and outflow conditions at the left and right boundaries of $C$). The extension is computed on the domain $C_E$, shown in gray. For this domain $C_E$ does not coincide with $C$ but does contain $E$. No slip boundary conditions would be imposed on $\Gamma$, the interior boundary that separates $\Omega$ and $E$.}
	\label{stokesfig:domain}
\end{figure}

Suppose that we wish to solve the Poisson problem
\begin{subequations}
	\label{stokeseq:poisson}
	\begin{align}
		\Delta u	&=	f	&	&	\text{in }\Omega,	\label{stokeseq:poisson:momentum}	\\
		u			&=	g	&	&	\text{on }\Gamma, \label{stokeseq:poisson:boundary}
	\end{align}
\end{subequations}
in an arbitrary smooth domain $\Omega$. The IBSE method works by smoothly extending the unknown solution $u$ of a PDE from $\Omega$ to a simplier domain $C$. The domain $\Omega$ will be referred to as the \emph{physical domain} and $C$ will be referred to as the \emph{computational domain}. It is assumed that $\Omega\subset C$. The \emph{interior} boundary $\Gamma=\partial\Omega$ will be assumed to be smooth, not self-intersecting, and must separate $C$ into the two disjoint regions $\Omega$ and $E=C\setminus\overline\Omega$, referred to as the \emph{extension domain}. Neither $\Omega$ nor $E$ need be connected. An additional simple computational domain $C_E$ is required to solve the auxiliary equation that defines the extension of the unknown solution. Although $C_E$ may coincide with $C$, it need only contain $E$. The boundary of the computational domain $C$ will be denoted by $\partial C$. Two typical domains are shown in \Cref{stokesfig:domain}. In \Cref{stokesfig:domain:1}, $C$ is periodic, and no boundary conditions need to be specified on $\partial C$. In \Cref{stokesfig:domain:2}, Dirichlet conditions are imposed on $\partial C$. To prevent confusion with the interior boundary $\Gamma$, the boundary $\partial C$ and its respective boundary condition will only be explicitly included in equations when required to prevent confusion, e.g. in \Cref{stokessection:numerics:summary,stokessection:stokes_test:flow_around_cylinder:finite_difference}.

We will denote the extension of the unknown solution $u$ by $\xi$. This extension of the solution is then used to define a volumetric forcing $\mathcal{F}_e=\Delta\xi$ in the region $E$. With this forcing, an extended problem in all of the simple domain $C$ may be solved:
\begin{subequations}
	\label{stokeseq:poisson:ib_thing}
	\begin{align}
		\Delta u_e - \chi_E\mathcal{F}_e &=	\chi_\Omega f	&&\text{in }C,	\\
		u_e						&=	g	&&\text{on }\Gamma.
	\end{align}
\end{subequations}
The solution $u_e$ gives the desired solution $u$ in $\Omega$ and is equal to $\xi$ in $E$. Because $\xi$ was chosen to be a smooth extension to $u$, the function $u_e$ is globally smooth in $C$.

The extension $\xi$ to the unknown function $u$ is defined as a solution to a high-order PDE which takes for its boundary conditions matching criteria of the form $\frac{\partial^j \xi}{\partial n^j} = \frac{\partial^j u}{\partial n^j}$. This allows the extension to be defined by a small number of unknowns (proportional to the number of points used to discretize the boundary). The extension PDE for $\xi$ is solved efficiently in the simple domain $C_E$ using an Immersed Boundary type method.

In order to succinctly describe the methodology we will require some notation. We define the \emph{spread operator}:
\begin{align}
	\label{stokeseq:derivative_spread_operator}
	(S_{(j)}F)(x) = (-1)^j\int_\Gamma F_j(s)\frac{\partial^j\delta(x-X(s))}{\partial n^j}\,dX(s)
\end{align}
and the \emph{interpolation operator}:
\begin{align}
	\label{stokeseq:derivative_interpolation_operator}
	(S_{(j)}^*\xi)(s) = (-1)^j\int_C \xi(x)\frac{\partial^j\delta(x-X(s))}{\partial n^j}\,dx
\end{align}
for the $j^\text{th}$ normal derivative, 
where $X(s)$ is a parametrization for $\Gamma$ with $s$ in the parameter interval $\mathcal{I}_\Gamma$. The integral on the right hand side of \Cref{stokeseq:derivative_interpolation_operator} is notation for the action of the distribution $\frac{\partial^j\delta}{\partial n^j}$ on the smooth function $\xi$. The traditional Immersed Boundary spread and interpolation operators, $S_{(0)}$ and $S_{(0)}^*$ are denoted by $S$ and $S^*$ respectively. The interpolation operator maps function values (or the normal derivatives of the function) from $C$ to $\Gamma$, while the spread operator maps singular (and hyper-singular) forces supported on $\Gamma$ to $C$. To simplify equations presented in the forthcoming material, we also define the operators $T_k$, $T_k^*$, and $R_k^*$ by:
\begin{subequations}
	\label{eq:define_tks}
	\begin{align}
		T_k &= \sum_{j=0}^{k}S_{(j)},	\\
		T^*_k &= 
		\begin{pmatrix}
			S_{(0)}^*	&	S_{(1)}^*	&	\cdots	&	S_{(k)}^*
		\end{pmatrix}^\intercal,\\
		R^*_k &= 
		\begin{pmatrix}
			S_{(1)}^*	&	\cdots	&	S_{(k)}^*
		\end{pmatrix}^\intercal.
	\end{align}
\end{subequations}
The operator $T_k^*$ provides an interpolation of a function and its first $k$ normal derivatives to the boundary; $R_k^*$ provides an interpolation of the first $k$ normal derivatives to the boundary, but excludes the values; the spread operator $T_k$ represents a set of singular forces ($\delta$-like) and hyper-singular forces  (like the first $k$ normal derivatives of the $\delta$-function) on the boundary.

The central challenge of the IBSE method is to compute the smooth extension to an \emph{unknown} solution. We first discuss how to compute an extension to a given function. Let $v\in C^k(\Omega)$ be given. To compute a $C^k(C)$ extension to $v$, we solve the following high-order PDE in the region $E$:
\begin{subequations}
	\label{stokeseq:extension_part_domain}
	\begin{align}
		\mathcal{H}^k\xi	&=	0	&&\text{in }E,	\\
		\frac{\partial^j \xi}{\partial n^j} &= \frac{\partial^j v}{\partial n^j}	&&\text{on }\Gamma,\ 0\leq j\leq k.\label{stokeseq:extension_part_domain:bc}
	\end{align}
\end{subequations}
Here $\mathcal{H}^k$ is an appropriate differential operator such as the polyharmonic operator $\Delta^{k+1}$; we discuss the choice of this operator in more detail, from a numerical perspective, in \Cref{stokessection:numerics:extension}. This problem may be solved on the simpler domain $C_E$ using methodology directly analogous to the direct forcing Immersed Boundary method. The boundary conditions given in \Cref{stokeseq:extension_part_domain:bc} that force $\xi$ to share its first $k$ normal derivatives with $u$ along $\Gamma$ are enforced by the addition of unknown singular and hyper-singular forces supported on the boundary:
\begin{subequations}
	\label{stokeseq:extension_whole_domain}
	\begin{align}
		\mathcal{H}^k\xi(x) + (T_k F)(x)	&=	0	&	&\text{for }x\in C_E,	\\
		(S_{(j)}^*\xi)(s)	&=	\frac{\partial^j v}{\partial n^j}(s)	&	&\text{for } s\in\mathcal{I}_\Gamma,\ 0\leq j\leq k.
		\label{stokeseq:extension_whole_domain:boundary}
	\end{align}
\end{subequations}
Future equations similar to \Cref{stokeseq:extension_whole_domain:boundary} will be shortened simply to $S_{(j)}^*\xi=\frac{\partial^j v}{\partial n^j}$, and assumed to hold for all $s\in\mathcal{I}_\Gamma$. Notice that $\xi$ is not actually an extension to $v$: that is, $\xi(x)\neq v(x)$ in $\Omega$. We will only be interested in the function $\xi$ in $E$, and so need not form its literal extension (which is $\chi_\Omega v + \chi_E\xi)$.

To solve the Poisson problem given by \Cref{stokeseq:poisson} using the IBSE method, we instead solve the extended problem given in \Cref{stokeseq:poisson:ib_thing}. The forcing function $\mathcal{F}_e$ that is specified in $E$ must be chosen so that it forces the extended solution $u_e$ to be $C^k(C)$. Let $\xi$ smoothly extend $u$, that is, we ask that $\xi$ is globally smooth in $C$ and that it satisfies the constraints
\begin{equation}
	R_k^*\xi = R_k^*u
\end{equation}
at the interface $\Gamma$. These constraints require that the first $k$ normal derivatives of $\xi$ agree with the first $k$ normal derivatives of $u$ on the boundary. Notice that it is not enforced that the \emph{values} of $u$ agree with the values of $\xi$ on $\Gamma$; we will see momentarily that this is unnecessary. The forcing function $\mathcal{F}_e$ is then defined as:
\begin{equation}
	\mathcal{F}_e = \Delta\xi.
\end{equation}
Coupling these equations together, we obtain the IBSE formulation for the Poisson problem given by \Cref{stokeseq:poisson}:
\begin{subequations}
	\label{stokeseq:ibse_poisson}
	\begin{align}
		\Delta u_e - \chi_E\Delta\xi	&=	\chi_\Omega f	&&\text{in }C,	\label{stokeseq:ibse_poisson:a}	\\
		\mathcal{H}^k + T_k F				&=	0							&&\text{in }C_E,	\label{stokeseq:ibse_poisson:b}	\\
		R_k^*\xi	&=	R_k^*u_e,	\label{stokeseq:ibse_poisson:c}	\\
		S^*u_e		&=	0.	\label{stokeseq:ibse_poisson:d}
	\end{align}
\end{subequations}
These equations are \eqref{stokeseq:ibse_poisson:a} the extended Poisson problem, \eqref{stokeseq:ibse_poisson:b} the extension PDE, \eqref{stokeseq:ibse_poisson:c} the interface constraints that force $u$ and $\xi$ to share their first $k$ normal derivatives, and \eqref{stokeseq:ibse_poisson:d} the physical boundary condition to the Poisson problem. Both the interface constraints and the physical boundary condition are enforced by the unknown forces $F$ in the extension PDE. We will subsequently drop the notation $u_e$ and refer to the solution $u_e$ of the IBSE system of equations simply as $u$.

In $\Omega$, $u$ satisfies the physical Poisson equation; while in $E$, $u$ satisfies:
\begin{subequations}
	\begin{align}
		\Delta u	&= \Delta\xi	&	&\text{in }E,	\\
		\frac{\partial u}{\partial n}	&=	\frac{\partial \xi}{\partial n}	&	&\text{on }\Gamma,
	\end{align}
\end{subequations}
and thus $u=\xi$ in $E$ up to a constant. Since $\mathcal{H}^k\xi=-T_k F$ implies that $\xi\in H^{k+1}(C)$, then for $k\geq1$, the forcing $\chi_D f + \chi_E\Delta\xi\in L^2(C)$, and so $u\in H^2(C)$, and is thus at least $C^0(C)$ (for spatial dimensions $d\leq 3$)\footnote{The notation $H^2(\Omega)$ denotes the space of measurable functions with the property that the function itself, as well as its first two weak derivatives, are square integrable. The Sobolev embedding theorem implies that any function in $H^2(\Omega)$ is equal almost everywhere to a continuous function, at least in spatial dimensions $d=2$ and $d=3$ \cite{adams2003sobolev}.}. This continuity, along with the smoothness constraints guaranteed by the constraint that $R^*_k\xi=R^*_k u$ imply that $u$ is globally $C^k(C)$.

In \cite{Stein2015}, we verify that the IBSE formulation of the problem produces $C^k(C)$ solutions that converge at a rate of $\mathcal{O}(\Delta x^{k+1})$ in the grid-spacing $\Delta x$, for the Poisson problem, as well as the heat equation, Burgers equation, and the Fitzhugh-Nagumo equations.

\begin{remark}
	In \cite{Stein2015}, both a volumetric force $\chi_E\Delta\xi$ and a singular force distribution $SG$ are added to the Poisson equation, giving a slightly different formulation for the IBSE method:
	\begin{align}
		\label{stokeseq:the_system_big}
		\Delta u_e - \chi_E\Delta\xi + SG	&=	\chi_\Omega f	&&\text{in }C,	\\
		\mathcal{H}^k + T_k F				&=	0							&&\text{in }C_E,	\\
		T_k^*\xi	&=	T_k^*u_e,	\\
		S^*u_e		&=	0.	
	\end{align}
	The singular force distribution $SG$ acts to enforce that the boundary values of $\xi$ equal those of $u$, this is an unnecessary requirement for the IBSE method. We will proceed with the simpler formulation given in \Cref{stokeseq:ibse_poisson} when developing the methodology for the Stokes equation.
\end{remark}

\subsection{Solution of the general Stokes problem}
\label{stokes_subsection:recap}

We now focus on the modifications to the IBSE method required to solve the general Stokes problem:
\begin{subequations}
	\begin{align}
		\mathcal{L}\u + \grad p	&=	\f_\u	&	&	\text{in }\Omega,	\label{stokeseq:general_brinkmann_problem:momentum}\\
		\grad\cdot\u			&=	f_p	&	&	\text{in }\Omega,	\label{stokeseq:general_brinkmann_problem:brinkmann:conservation}\\
		\mathcal{B}\U			&=	\b	&	&	\text{on }\Gamma.
	\end{align}
	\label{stokeseq:general_brinkmann_problem}
\end{subequations}
We assume that $\f=(\f_u,\ f_p)$ and $\b$ are sufficiently smooth functions defined on $\Omega$ (for $\f$) and $\Gamma$ (for $\b$). Here $\mathcal{L}$ denotes the Helmholtz operator $\mathcal{L}=\alpha\mathbb{I}-\Delta$ and $\mathcal{B}$ denotes a boundary operator. We discuss the boundary operator associated with the Dirichlet problem ($\mathcal{B}\U=\u|_\Gamma$); imposition of other boundary conditions is similar and discussed in \cite{Stein2015}. Note that this equation reduces to the incompressible Stokes problem when $\alpha$ and $f_p$ are zero.

In the direct forcing Immersed Boundary (IB) method, this problem is solved in a simple computational domain $C$ by adding singular forces $\G$ supported on the boundary that act as Lagrange multipliers which enforce the boundary condition. For Dirichlet problems this can be represented as
\begin{subequations}
	\label{stokeseq:immersed_boundary}
	\begin{align}
		\mathcal{L}\u + \grad p + S\G	&=	\f_\u	&	&	\text{in }C,	\\
		\grad\cdot\u			&=	f_p										&	&	\text{in }C,	\\
		S^*\u			&=	\b.
	\end{align}
\end{subequations}
The singular forces supported on $\Gamma$ induce jumps in the normal derivatives of $\u$ at the boundary; the velocities $\u$ produced by the IB method are generically $C^0(C)$ while the pressure field typically has jump discontinuities. For Eulerian discretizations of the underlying PDE that are ignorant of the boundary, this leads to the slow convergence rate $\mathcal{O}(\Delta x)$ for the velocity field in $L^\infty(\Omega)$. The pressure field $p$ fails to converge near to the boundary, but does converge at the rate $\mathcal{O}(\sqrt{\Delta x})$ in $L^2(\Omega)$. 

In order to improve this convergence, the Immersed Boundary Smooth Extension (IBSE) method extends the unknown solutions $\u$ and $p$ to be $C^k(C)$ and $C^{k-1}(C)$, respectively. Suppose that the solution vector $\U$ to \Cref{stokeseq:general_brinkmann_problem} is known in $\Omega$. Define $\bxi_\u$ by
\begin{subequations}
	\label{stokeseq:define_extension:velocity}
	\begin{align}
		\mathcal{H}^{k}\bxi_\u	+ T_k\F_\u &=	0	&&\text{in }C_E,	\\
		R_k^*\bxi_\u	&=	R_k^*\u,
	\end{align}
\end{subequations}
and $\xi_p$ by
\begin{subequations}
	\label{stokeseq:define_extension:pressure}
	\begin{align}
		\mathcal{H}^{k-1}\xi_p	+ T_{k-1}F_p &=	0	&&\text{in }C_E,	\\
		T_{k-1}^*\xi_p	&=	T_{k-1}^* p.
	\end{align}
\end{subequations}
Recall that $\mathcal{H}^k$ is a high-order differential operator such as $\Delta^{k+1}$, which will be defined precisely in \Cref{stokessection:numerics:extension}. We make two remarks concerning the extension of the pressure.
\begin{remark}
	The pressure function $p$ is typically thought of as a Lagrange multiplier: it equals whatever is required to enforce the divergence constraint that $\grad\cdot\u=0$. One may suspect therefore that it is unnecessary to construct an explicit extension of the pressure function $p$ as we do in \Cref{stokeseq:define_extension:pressure}. Unfortunately, this is not true. Consider the problem:
	\begin{subequations}
		\begin{align}
			-\Delta \u + \grad p &=	\f_\u	&&\text{in }\Omega,	\\
			\grad\cdot \u	&=	0	&&\text{in }\Omega, \\
			\u	&= \b	&&\text{on }\Gamma,
		\end{align}
	\end{subequations}
	along with the modified problem where $\f_\u$ is replaced by $\tilde\f_\u=\f_\u + \grad\phi$ for any scalar function $\phi$. The pressure will be changed: $\tilde p = p - \phi$. But the solution $\u$ will be left unchanged, and thus its extension functions $\bxi_\u$ will be unchanged. If the extension of the pressure function depends only locally on $\bxi_\u$, it cannot extend both $p$ and $\tilde p$ smoothly.
\end{remark}
\begin{remark}
	\label{stokesremark:pressure_extension_regularity}
	The choice to extend the pressure function $p$ to be $C^{k-1}(C)$ rather than $C^k(C)$ yields improved numerical stability. When $\u\in C^k(C)$, the structure of the Stokes equations implies that $p\in C^{k-1}(C)$. Applying the operator $T_k^*$ to the function $p$ with that global regularity yields an inconsistent estimate of $\frac{\partial^k p}{\partial n^k}$, i.e. $S^*_{(k)}p = \frac{\partial^k p}{\partial n^k} + \mathcal{O}(1)$, and thus enforcing that $T_k^*\xi_p=T_k^*p$ does not actually enforce that $\frac{\partial^k \xi_p}{\partial n^k}=\frac{\partial^k p}{\partial n^k}$.
\end{remark}

From these extensions, we may define a set of forces $(\boldsymbol{\mathcal{F}}_\u, \mathcal{F}_p)$:
\begin{subequations}
	\begin{align}
		\boldsymbol{\mathcal{F}}_\u	&=	\mathcal{L}\bxi_\u+\grad \xi_p,	\\
		\mathcal{F}_p	&=	\grad\cdot\bxi_\u.
	\end{align}
\end{subequations}
Analogous to \Cref{stokeseq:poisson:ib_thing}, to smoothly extend the Stokes problem given by \Cref{stokeseq:general_brinkmann_problem} from $\Omega$ to $C$, we add $\boldsymbol{\mathcal{F}}_\u$ to \Cref{stokeseq:general_brinkmann_problem:momentum} and $\boldsymbol{\mathcal{F}}_p$ to \Cref{stokeseq:general_brinkmann_problem:brinkmann:conservation} in the extension region $E$ to obtain
\begin{subequations}
	\begin{align}
		\mathcal{L}\u + \grad p	&=	\chi_\Omega\f _\u+ \chi_E\mathcal{L}\bxi_\u + \chi_E\grad\xi_p	&	&	\text{in }C,	\\
		\grad\cdot\u			&=	\chi_\Omega f_p + \chi_E\grad\cdot\bxi_\u								&	&	\text{in }C,	\\
		\mathcal{B}\U			&=	\b	&	&	\text{on }\Gamma.
	\end{align}
\end{subequations}
In the physical domain $\Omega$, we have that
\begin{subequations}
	\begin{align}
		\mathcal{L}\u + \grad p	&=	\f_\u &	&	\text{in }\Omega,	\\
		\grad\cdot\u			&=	f_p &	&	\text{in }\Omega,	\\
		\mathcal{B}\U			&=	\b	&	&	\text{on }\Gamma,
	\end{align}
\end{subequations}
and thus $\U$ solves the original Stokes problem given in \Cref{stokeseq:general_brinkmann_problem}. In the extension domain $E$, we have that
\begin{subequations}
	\begin{align}
		\mathcal{L}(\u-\bxi_\u) + \grad (p-\xi_p)	&=	0 &	&	\text{in }E,	\\
		\grad\cdot(\u-\bxi_\u)			&=	0 &	&	\text{in }E,	\\
		\frac{\partial(\u-\bxi_\u)}{\partial n}			&=	0	&	&	\text{on }\Gamma,
	\end{align}
\end{subequations}
Thus $\u=\bxi_\u$ and $p=\xi_p$, both up to a constant (in $E$). Since $\xi_p$ is enforced to match $p$ at the boundary by \Cref{stokeseq:define_extension:pressure}, they are equal at the boundary, and thus $p=\xi_p$ in $E$, implying that $p$ is globally smooth in $C$, while the global regularity for $\u$ is implied by the same argument given in \Cref{stokes_subsection:recap}. Numerical approximations of these solutions may converge rapidly, even when the PDEs are discretized using an underlying Cartesian grid discretization that is ignorant of the boundary. The difficulty in this formulation comes from the fact that the extensions $\bxi_\u$ and $\xi_p$ \emph{depend on the unknown solution $\U$}.

Combining \Cref{stokeseq:define_extension:velocity,stokeseq:define_extension:pressure,stokeseq:general_brinkmann_problem}, we obtain a coupled system for the solution $\U$ together with the smooth extensions $(\bxi_\u,\xi_p)$ to that solution:\vspace{-0.5em}\\
\begin{minipage}{\textwidth}
	\begin{subequations}
		\label{stokeseq:the_stokes_system}
		\begin{empheq}[left=\text{Stokes eq. }\hspace{0.313in}\empheqlbrace\hspace{0.1in}]{align}
			\vphantom{\mathcal{H}^k}\mathcal{L}\u + \grad p - \chi_E\mathcal{L}\bxi_\u - \chi_E\grad\xi_p	&=	\chi_\Omega\f_\u	\hspace{0.65in}	\text{in }C,	\label{stokeseq:the_stokes_system:a} \\
			\vphantom{\mathcal{H}^k}\grad\cdot\u - \chi_E\grad\cdot\bxi_\u			&=	\chi_\Omega f_p	\hspace{0.65in}	\text{in }C, \label{stokeseq:the_stokes_system:b}
		\end{empheq}
		\begin{empheq}[left=\text{Extension eq. }\hspace{0.1in}\empheqlbrace\hspace{0.753in}]{align}
			\mathcal{H}^k\bxi_\u + T_k\F_\u	&=	0	\hspace{0.885in}\text{in }C_E,	\label{stokeseq:the_stokes_system:c}	\\
			\mathcal{H}^{k-1}\xi_p + T_{k-1} F_p	&=	0	\hspace{0.885in}\text{in }C_E,\label{stokeseq:the_stokes_system:d}
		\end{empheq}
		\begin{empheq}[left=\text{Bdy. matching }\hspace{0.039in}\empheqlbrace\hspace{1.403in}]{align}
			R_k^*\bxi_\u		&=	R_k^*\u,	\label{stokeseq:the_stokes_system:e}\\
			T_{k-1}^*\xi_p	&=	T_{k-1}^*p,\label{stokeseq:the_stokes_system:f}
		\end{empheq}
		\begin{empheq}[left=\text{Phys. Bdy. Cond. }\hspace{1.532in}]{align}
					S^*\u		&=	\b.\label{stokeseq:the_stokes_system:g}
		\end{empheq}
	\end{subequations}
\end{minipage}\\
\vspace{0.5em}

We will refer to this system of equations as the IBSE-$k$ equations. The remainder of this paper is dedicated to demonstrating that these equations provide accurate solutions to the Stokes and Navier-Stokes equations, that they allow a flexible choice of discretization, and that they admit an efficient inversion strategy.
We organize the presentation of information as follows:
\begin{enumerate}
	\item In \Cref{stokessection:numerical_implementation}, we discuss the details of numerical implementation that do not depend upon the underlying choice of discretization.
	\item In \Cref{stokessection:stokes_test_analytic}, we solve a Stokes problem with an analytic solution that allows us to carefully verify the convergence rates of the solution $\U$ (up to third-order for $\u$ and second-order for $p$, in $L^\infty(\Omega)$), as well as the convergence of elements of the fluid stress tensor $\sigma$ (up to second-order in $L^\infty(\Omega)$).
	\item In \Cref{stokessection:stokes_test:flow_around_cylinder}, we describe both a Fourier spectral and a finite-difference discretization to the confined channel flow around a cylinder problem for incompressible Stokes flow. We validate the convergence rates for solutions produced by the IBSE method and compare to known benchmarks.
	\item In \Cref{stokessection:navier_stokes}, we discuss time-stepping in the IBSE framework, solve both steady and unsteady Navier-Stokes problems, and demonstrate rapid convergence of the solutions and agreement with known benchmarks.
\end{enumerate}


\section{Numerical implementation}
\label{stokessection:numerical_implementation}

In this section, we discretize the coupled IBSE-$k$ system given in \Cref{stokeseq:the_stokes_system}. We assume that the functions are discretized using an underlying Cartesian mesh with uniform\footnote{For the discretization described here, the mesh needs to be uniform only in a small neighborhood of the boundary, with the same width as the regularized $\delta$-function used to discretize the operators $S_{(j)}$ and $S^*_{(j)}$. In fact, the mesh need not even be uniform, but this would add some complexity to the discretization \cite{Akiki2016}.} grid-spacing $\Delta x$ but otherwise make no assumptions regarding the underlying discretization of the functions and differential operators. These final details will be treated in \Cref{stokessection:stokes_test_analytic} for a Fourier spectral discretization, and in \Cref{stokessection:stokes_test:flow_around_cylinder:finite_difference} for a standard second-order finite difference scheme. Little varies in the details.

\begin{enumerate}
	\item In \Cref{stokessection:numerics:operators}, we discretize the spread $(S_{(j)})$ and interpolation $(S_{(j)}^*)$ operators for the $j^{th}$ normal derivative that were introduced in \Cref{stokeseq:derivative_spread_operator,stokeseq:derivative_interpolation_operator}. Note that the discretization of these operators automatically induces a discretization for the composite operators $T_k$, $T_k^*$, and $R_k^*$, defined in \Cref{eq:define_tks}.
	\item In \Cref{stokessection:numerics:extension}, we discuss how we choose the extension operator $\mathcal{H}^k$ introduced in \Cref{stokeseq:extension_whole_domain}.
	\item In \Cref{stokessection:numerics:inversion}, we describe an efficient inversion strategy for the IBSE-$k$ system given by \Cref{stokeseq:the_stokes_system}.
	\item In \Cref{stokessection:numerics:summary}, we provide a brief summary of the algorithm.
\end{enumerate}

\subsection{Discretization of singular integrals}
\label{stokessection:numerics:operators}

Let the boundary $\Gamma$ be parametrized by the function $X(s)$.  In all examples in this manuscript, we work with a two-dimensional fluid, and the boundary $\Gamma$ is one-dimensional and closed, with the single parameter $s$ defined on the periodic interval $[0,2\pi]$.  The discrete version of the spread operator for the $j^\text{th}$ normal derivative, $S_{(j)}:\Gamma\to C$ defined in \Cref{stokeseq:derivative_spread_operator} requires a regularized $\delta$-function and a discretization of the integral over $\Gamma$. The regularized $\delta$-function that we use is analytically three times differentiable, satisfies four discrete moment conditions, and has a support width of $16\Delta x$. The $\delta$-function and its properties are introduced in \cite{Stein2015,Stein2016}; this discrete $\delta$ function will be represented by $\tilde\delta$. Multivariate $\delta$-functions are computed as Cartesian products of the univariate $\delta$ and also denoted by $\tilde\delta$. Normal derivatives of $\tilde\delta$ are computed by the formula \cite{john1982partial}
\begin{equation}
	\frac{\partial^j\tilde\delta}{\partial n^j} = n_{i_1}\cdots n_{i_2}\frac{\partial^j\tilde\delta}{\partial x_{i_1}\cdots \partial x_{i_j}},
	\label{stokeseq:normal_derivative_formula}
\end{equation}
where the Einstein summation convention has been used to indicate sums over repeated indices and\\ $\partial^j\tilde\delta/\partial x_{i_1}\cdots \partial x_{i_j}$ is computed as Cartesian products of the appropriate derivatives of the one dimensional $\tilde\delta$.  For example, $\partial\tilde\delta/\partial n$ in two dimensions is computed as
\begin{equation}
	\frac{\partial\tilde\delta}{\partial n} = n_x \tilde\delta' \otimes \tilde\delta + n_y\tilde\delta\otimes\tilde\delta'.
\end{equation}
Discretization of the integral over $\Gamma$ is made by choosing $n_\text{bdy}$ quadrature nodes $\tilde\Gamma=\{X_i\}_{i=1}^{n_\text{bdy}}$, equally spaced in the parameter interval $\mathcal{I}_\Gamma=[0,2\pi]$ so that $X_i=X(s_i)$ and $s_i=(i-1)2\pi/n_\text{bdy}$.  Quadrature weights are computed at each quadrature node to be $\Delta s_i=\left\|\frac{\partial X}{\partial s}(s_i)\right\|_2$; this is a spectrally accurate quadrature rule for the integral of smooth periodic functions on $\Gamma$.  The discrete spread operator $S_{(j)}$ maps functions sampled at points in $\tilde\Gamma$ to $C$ by
\begin{equation}
	(S_{(j)}F)(x) = \sum_{i=1}^{n_\text{bdy}} F(s_i)\frac{\partial^j\tilde\delta(x-X_i)}{\partial n^j}\Delta s_i.
\end{equation}
We do not adopt explicit notation to distinguish between the analytic and discrete operators.  The number of points in the quadrature is chosen so that $\Delta s\approx 2\Delta x$.  This choice of node-spacing is wider than that recommended for the traditional IB method \cite{Peskin2002} but has been observed empirically to be the optimal choice in other studies of \emph{direct-forcing} IB methods \cite{kallemov2016immersed}.
We define the interpolation operator $S_{(j)}^*$ by the adjoint property $\left<u,S_{(j)}F\right>_C=\left<S_{(j)}^*u,F\right>_\Gamma$, but note that the discrete interpolation operator $S_{(j)}^*$ produces a discrete function
\begin{equation}
	(S_{(j)}^*u)(s_k) = \int_C u(x)\frac{\partial^j\tilde\delta(x-X_k)}{\partial n^j}\,dx.
\end{equation}
Discrete integrals over $C$ are straightforward sums computed over the underlying uniform Cartesian mesh, and may be efficiently evaluated due to the finite support of $\tilde\delta$. When acting on vector functions, all operators defined in this section are assumed to operate elementwise.

\subsection{Solution of the extension equations and the choice of extension operator $\mathcal{H}^k$}
\label{stokessection:numerics:extension}

Due to the difficulty in accurately and efficiently inverting the high-order differential operator needed to define the smooth extensions, we solve these equations utilizing a Fourier spectral discretization regardless of the discretization used for the Stokes equations. Let $m$ be the largest wave-number associated with the discrete Fourier transform (DFT) on the discretized domain $C_E$. We choose the extension operator $\mathcal{H}^k$ introduced in \Cref{stokeseq:extension_part_domain} to be:
\begin{equation}
	\label{stokeseq:Helmholtz_definition}
	\mathcal{H}^k = \Delta^{k+1} + (-1)^{k+1}\Theta(k,m).
\end{equation}
Here $\Theta$ is a positive scalar function that depends on the smoothness of the extension ($k$) and the largest wave-number ($m$).  The function $\Theta$ is chosen to mitigate the numerical condition number of the operator $\mathcal{H}^k$:
\begin{equation}
	\kappa = 1 + \frac{m^{2(k+1)}}{\Theta}.
\end{equation}
Minimizing the condition number $\kappa$ (by taking $\Theta$ to be large) must be balanced against the need to resolve the intrinsic length scale $L=\Theta^{-1/2(k+1)}$ introduced to the problem. We set the length scale $L$ as a function of $\Delta x$, in effect allowing the extension to decay to zero over some number of grid cells. We thus choose
\begin{equation}
	\label{stokeseq:Theta_definition}
	\Theta^* = \left(\frac{1}{N\Delta x}\right)^{2(k+1)},
\end{equation}
where $N$ is proportional to the number of gridpoints over which the extension decays. The choice of $N$ can impact the accuracy of the scheme: if $N$ is too large, the condition number of $\mathcal{H}^k$ may be high with respect to the precision of the computer being used, and numerical instability causes a degradation of the quality of the solution. If $N$ is too small, short length scales that are not resolvable by the discretization are introduced to the problem. For all examples presented in this paper, $N$ is chosen to be $200$.
 In \Cref{stokessection:stokes_test_analytic}, we investigate the effect of the choice of $N$ on a specific problem.

\subsection{Inversion of the IBSE-$k$ system, \Cref{stokeseq:the_stokes_system}}
\label{stokessection:numerics:inversion}

In this section, we describe an efficient inversion strategy for \Cref{stokeseq:the_stokes_system}. We will assume that the discrete Stokes operator is invertible\footnote{The null-space in the pressure function $p$ should be fixed by the discrete operator; e.g. by enforcing that $\int_C p=0$.}. The details for when it is not invertible (e.g. when solving on a periodic domain with $\mathcal{L}=-\Delta$) are analogous to the Poisson problem, and the solution is presented in \cite{Stein2015}. In matrix form, \Cref{stokeseq:the_stokes_system} is given by:
\begin{equation}
	\label{stokeseq:the_discrete_system_matrix}
	\left(
	\begin{array}{cccc|cc}
		\mathcal{L}	&	\grad		&	-\chi_E\mathcal{L}	&	-\chi_E\grad	&		&			\\
		\grad\cdot	&				&	-\chi_E\grad\cdot		&							&		&			\\
					&				&	\mathcal{H}^k					&							&	T_k	&			\\
					&				&									&	\mathcal{H}^{k-1}		&		&	T_{k-1}	\\
		\hline
			S^*	&	&	&	\\
			R_k^*	&	&	-R_k^*	&	&	&	\\
					&	T_{k-1}^*	&									&	-T_{k-1}^*				&		&			\\
	\end{array}
	\right)
	\begin{pmatrix}
		\u	\\	p	\\	\bxi_u	\\	\xi_p	\\ \hline	\F_\u	\\	F_p
	\end{pmatrix}
	=
	\begin{pmatrix}
		\chi_\Omega\f_\u	\\ \chi_\Omega f_p	\\	0	\\	0	\\	\hline
			\b	\\
			0	\\
		0
	\end{pmatrix}.
\end{equation}
This system is large, but the number of constraints (the last three equations) is relatively small. By block Gaussian elimination, we can find a Schur-complement for this matrix, which we label $SC$:
\begin{equation}
	\label{stokeseq:define_schur_complement}
	SC = 
	\begin{pmatrix}
			S^*	\\
			R_k^*	&	&	-R_k^*	\\
		&	T_{k-1}^*	&									&	-T_{k-1}^*
	\end{pmatrix}
	\begin{pmatrix}
		\mathcal{L}	&	\grad		&	-\chi_E\mathcal{L}	&	-\chi_E\grad	\\
		\grad\cdot	&				&	-\chi_E\grad\cdot		&							\\
					&				&	\mathcal{H}^k					&							\\
					&				&									&	\mathcal{H}^{k-1}			
	\end{pmatrix}^{-1}
	\begin{pmatrix}
		0	&	0	\\
		0	&	0	\\
		T_k	&	0	\\
		0	&	T_{k-1}	
	\end{pmatrix}
\end{equation}
along with an associated system of equations for the Lagrange multipliers $\F_\u$ and $F_p$,
\begin{equation}
	\bigg(		\quad SC \quad \bigg)
	\begin{pmatrix}
			\F_\u	\\	F_p
	\end{pmatrix}
	=
	\begin{pmatrix}
			S^*	\\	R^*_{k}	\\
			&	T^*_{k-1}
	\end{pmatrix}
	\begin{pmatrix}
		\mathcal{L}	&	\grad	\\
		\grad\cdot
	\end{pmatrix}^{-1}
	\begin{pmatrix}
		\chi_\Omega\f_\u	\\ \chi_\Omega f_p
	\end{pmatrix}
		-
	\begin{pmatrix}
			\b	\\	0	\\	0
	\end{pmatrix}.
	\label{stokeseq:schur_system}
\end{equation}
The size of $SC$ is comparatively small, only $(2(k+1)+k)n_\text{bdy}$ square.  This equation is the key to the efficiency of the algorithm, as it allows the Lagrange multipliers $\F_\u$ and $F_p$ to be computed rapidly without first solving for $\u$, $p$, $\bxi_\u$, or $\xi_p$.  Once $\F_\u$ and $F_p$ are determined, $\bxi_\u$ and $\xi_p$ may be found by solving \Cref{stokeseq:the_stokes_system:c,stokeseq:the_stokes_system:d}. With $\bxi_\u$ and $\xi_p$ known, $\u$ and $p$ may be computed simply by solving the Stokes equation in $C$, ignoring the boundary $\Gamma$ (that is, by solving \Cref{stokeseq:the_stokes_system:a,stokeseq:the_stokes_system:b}).

Rapid inversion of the system of equations for $\F_\u$ and $F_p$ given in \Cref{stokeseq:schur_system} is not a trivial task.  Because we have restricted to problems set on stationary domains and the size of $SC$ is small, it is feasible to proceed using dense linear algebra.  We form $SC$ by repeatedly applying it to basis vectors.  This operation is expensive: for two dimensions it is $\mathcal{O}(N^{3/2}\log N)$ in the total number of unknowns $N=n^2$.  The Schur-complement is then factored by the LU algorithm provided by LAPACK \cite{laug}.  Once this factorization is computed \Cref{stokeseq:schur_system} can be solved rapidly.  This Schur-complement depends only on the domain and the discretization, so the LU-decomposition can be reused to solve multiple problems on the same domain or in each timestep when solving time-dependent problems.

\subsection{Outline of the methodology}
\label{stokessection:numerics:summary}

We provide a brief outline of the steps required for solving the general Stokes problem given by \Cref{stokeseq:general_brinkmann_problem} using the IBSE-$k$ method. The algorithm proceeds as follows:
\begin{enumerate}
	\item Form the right hand side of the Schur-complement equation, as defined in \Cref{stokeseq:schur_system}. This may be done by solving
	\begin{subequations}
		\label{stokeseq:discrete_stokes_system}
		\begin{align}
			\mathcal{L}\u_0 + \grad p_0	&=	\chi_\Omega\f_\u	&&\text{in }C,	\\
			\grad\cdot\u_0			&=	\chi_\Omega f_p		&&\text{in }C,	\\
			\mathcal{B}(\u_0,p_0)		&=	\b_C	&&\text{on }\partial C,
		\end{align}
	\end{subequations}
	for $(\u_0, p_0)$ and then computing $S^*\u_0-\b$, $R^*_k\u_0$ and $T^*_{k-1}p_0$. Note that the \emph{internal} boundary $\Gamma$ is ignored when computing the solution to \Cref{stokeseq:discrete_stokes_system}.
	\item Compute $\F_\u$ and $F_p$ by solving the Schur-complement system given by \Cref{stokeseq:schur_system}, given the right hand side computed above. We assume that the Schur complement $SC$ has been formed and its LU decomposition found as a pre-computation.
	\item Compute $\bxi_\u$ and $\xi_p$ from $\F_\u$ and $F_p$ by solving \Cref{stokeseq:the_stokes_system:c,stokeseq:the_stokes_system:d}. We take $C_E$ to always have periodic boundary conditions, and these equations may be simply inverted using the fast Fourier transform.
	\item With $\bxi_\u$ and $\xi_p$ known, find $\u$ and $p$ by solving \Cref{stokeseq:the_stokes_system:a,stokeseq:the_stokes_system:b}. To be precise, we solve
	\begin{subequations}
		\begin{align}
			\mathcal{L}\u + \grad p	&=	\chi_E\mathcal{L}\bxi_\u + \chi_E\grad\xi_p + \chi_\Omega\f_{\u}	&&\text{in }C,	\\
			\grad\cdot\u			&=	\chi_E\grad\cdot\bxi_\u + \chi_\Omega f_p		&&\text{in }C,	\\
			\mathcal{B}(\u,p)	&=				\b_C	&&\text{on }\partial C.
		\end{align}
	\end{subequations}
	Note again that the \emph{internal} boundary $\Gamma$ is ignored when computing this solution.
\end{enumerate}

Step two of this algorithm requires that the Schur-complement $SC$ be formed and factored as a pre-processing step. The Schur complement operator maps a set of singular (and hyper-singular) forces $(\F_\u, F_p)$ to the constraints given by equations \Cref{stokeseq:the_stokes_system:e,stokeseq:the_stokes_system:f,stokeseq:the_stokes_system:g}. We form the matrix by applying $SC$ to basis vectors. To apply $SC$ to any set of singular forces, first compute $\bxi_\u$ and $\xi_p$ from $\F_u$ and $F_p$ by solving \Cref{stokeseq:the_stokes_system:c,stokeseq:the_stokes_system:d}. Then $\u$ and $p$ may be found by solving
\begin{subequations}
	\begin{align}
		\mathcal{L}\u + \grad p	&=	\chi_E\mathcal{L}\bxi_\u + \chi_E\grad\xi_p	&&\text{in }C,	\\
		\grad\cdot\u			&=	\chi_E\grad\cdot\bxi_\u	&&\text{in }C,	\\
		\mathcal{B}(\u,p)	&=				0	&&\text{on }\partial C,
	\end{align}
\end{subequations}
and finally with $\u$, $p$, $\bxi_\u$ and $\xi_p$ known, then $S^*\u$, $R_k^*\u-R_k^*\bxi_\u$, and $T_{k-1}^*p - T_{k-1}^*\xi_p$ can be evaluated. Note that $\f_\u$, $f_p$ and $\b_C$ do not enter into this computation. The Schur complement $SC$ depends only on the geometry and the discrete operators $\mathcal{L}$, $\grad$, $\grad\cdot$.



\section{Results: Stokes equation with analytic solution}
\label{stokessection:stokes_test_analytic}

\begin{figure}
	\centering
	\hspace*{\fill}
	\begin{subfigure}[b]{0.3\textwidth}
		\centering
		\includegraphics[width=\textwidth]{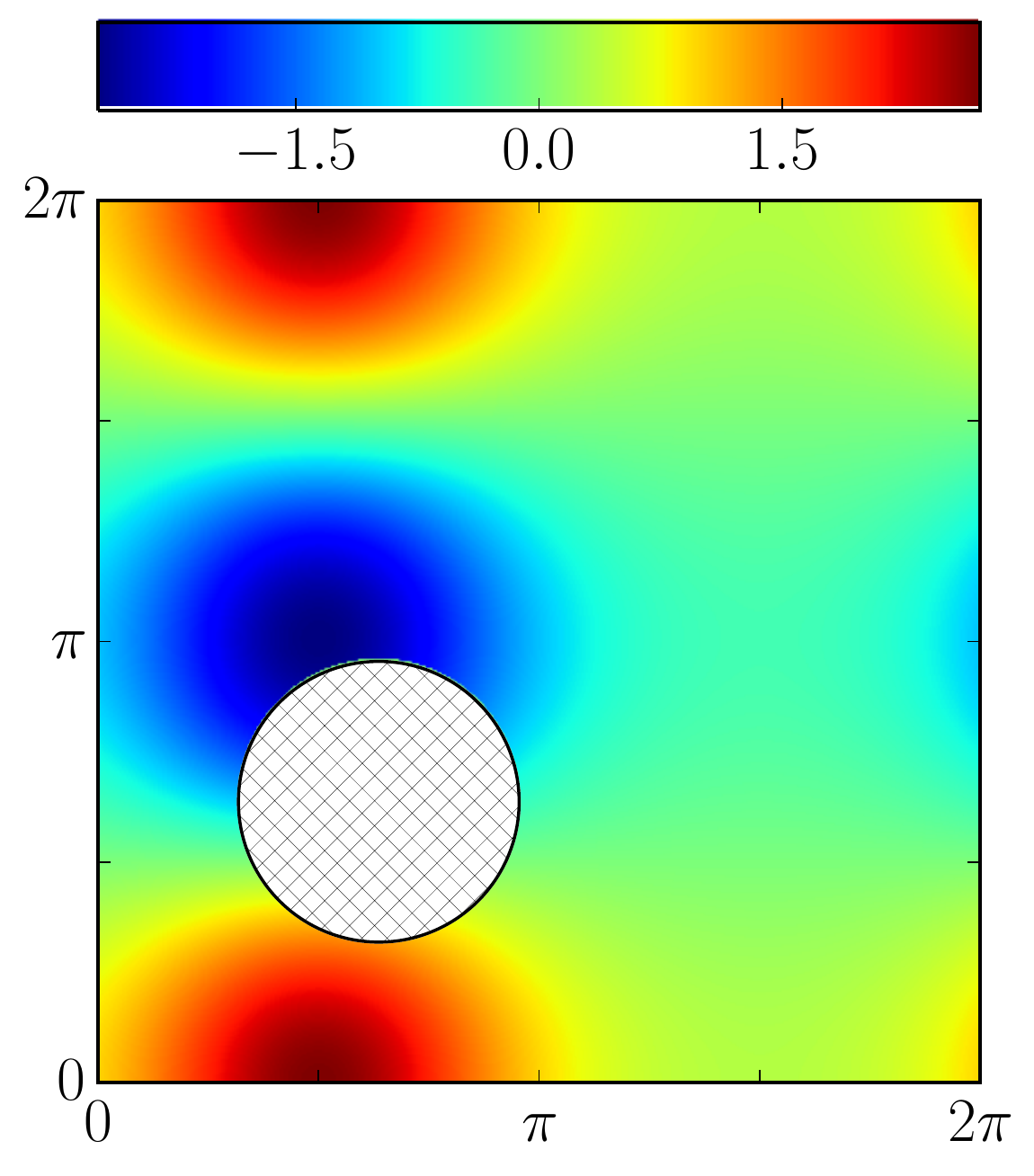}
		\subcaption{$u$}
		\label{stokesfig:manufactured_analytic_solution:u}
	\end{subfigure}
	\hfill
	\begin{subfigure}[b]{0.3\textwidth}
		\centering
		\includegraphics[width=\textwidth]{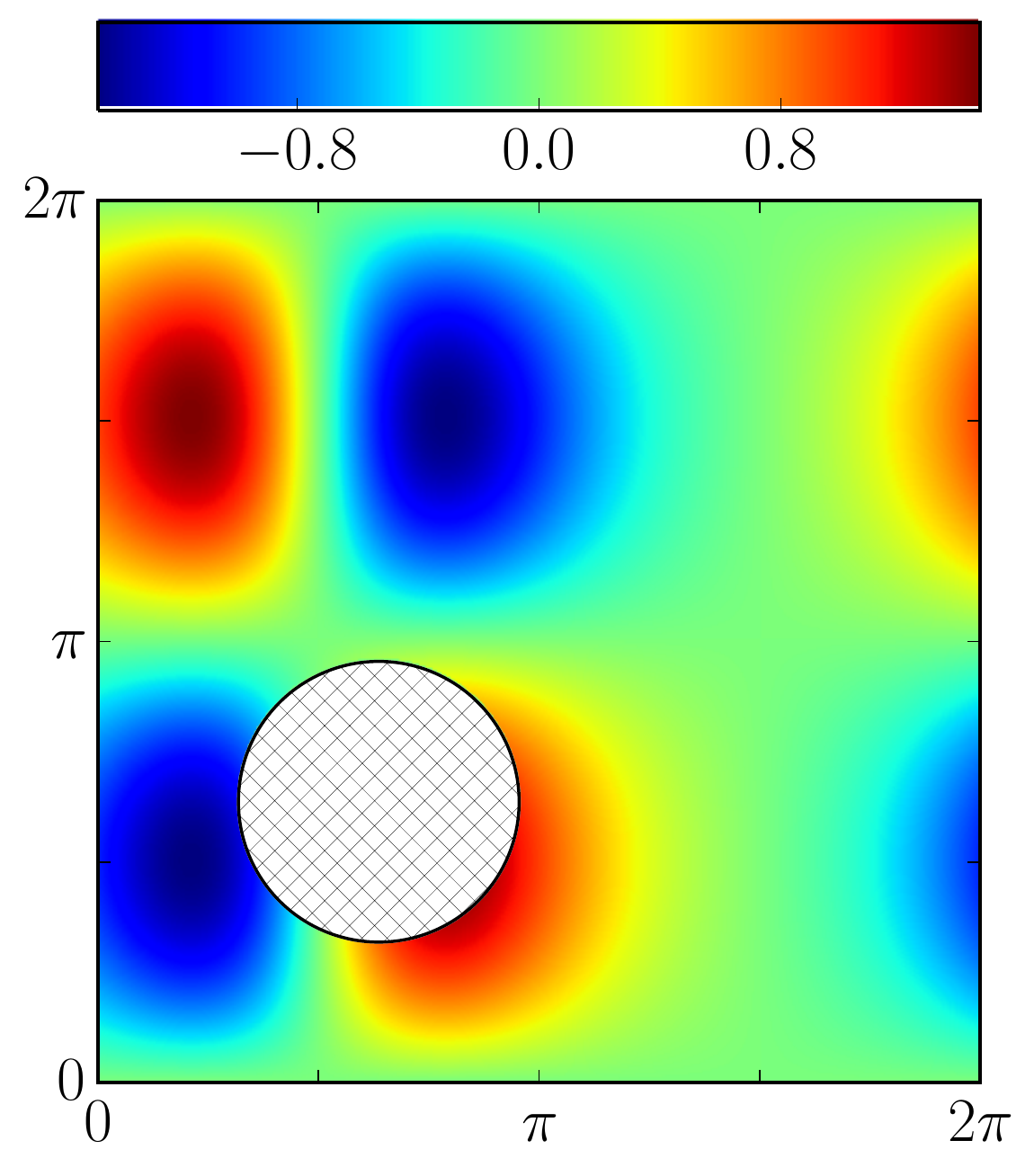}
		\subcaption{$v$}
		\label{stokesfig:manufactured_analytic_solution:v}
	\end{subfigure}
	\hfill
	\begin{subfigure}[b]{0.3\textwidth}
		\centering
		\includegraphics[width=\textwidth]{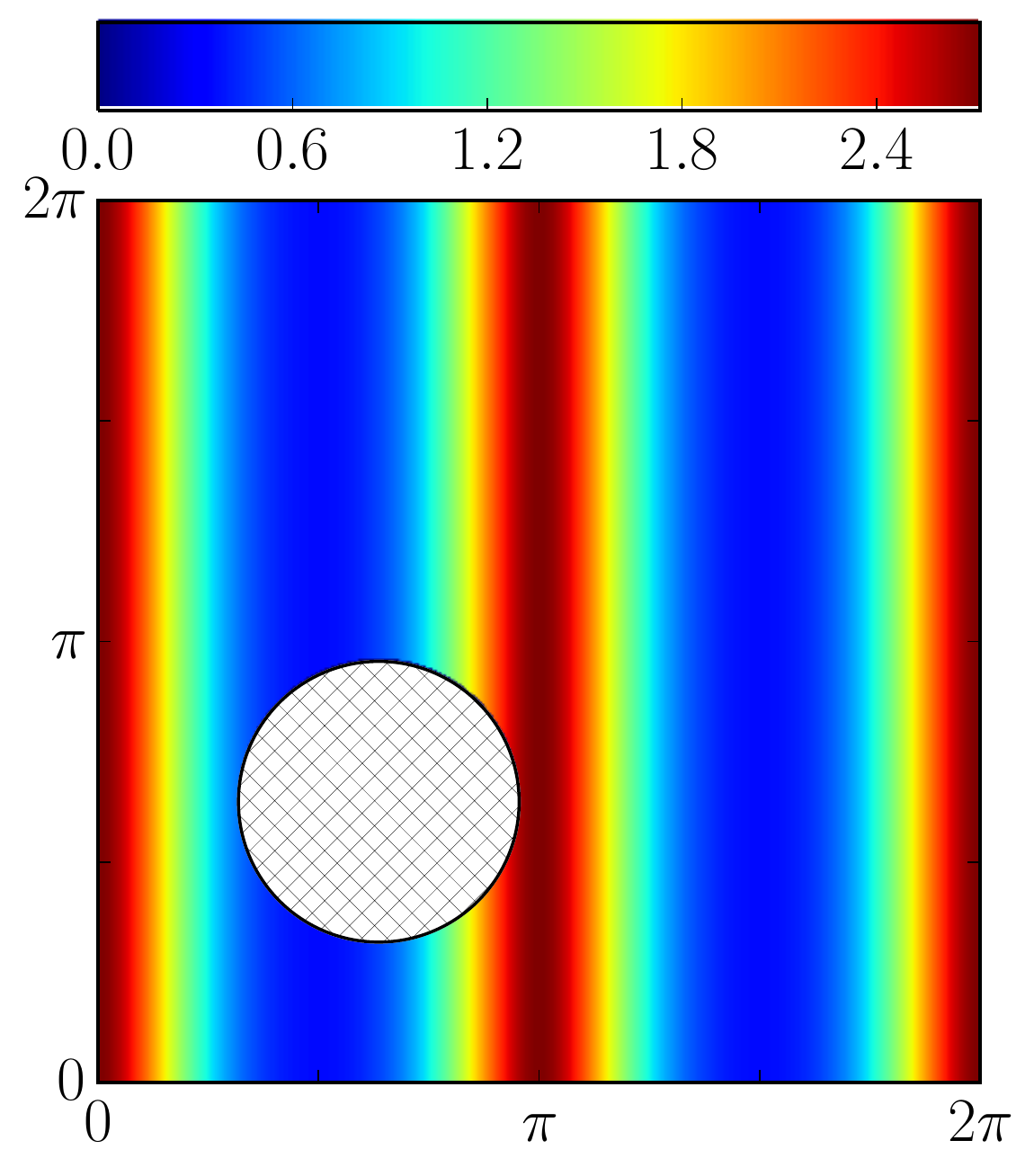}
		\subcaption{$p$}
		\label{stokesfig:manufactured_analytic_solution:p}
	\end{subfigure}
	\hspace*{\fill}
	\caption[The analytic solution and domain for one of the Stokes test problems used to test the accuracy of the Immersed Boundary Smooth Extension (IBSE) method.]{The analytic solution $\U$ and the physical domain $\Omega$ to the problem defined by \Cref{stokeseq:manufactured_problem,stokeseq:manufactured_solutions}. The extension domain $E$ is denoted with crosshatches.}
	\label{stokesfig:manufactured_analytic_solution}
\end{figure}

To demonstrate the improved convergence properties of the IBSE method as compared to the IB method, we construct a simple example with an analytic solution and compare the accuracy of the solutions with those produced by the Immersed Boundary method. For the numerical comparisons with the IB method provided in this section, we show two sets of values: those produced using a regularized $\delta$-function commonly used in IB simulations \cite{Lai2000}:
\begin{equation}
	8\,\delta(x) = 
	\begin{cases}
		3 - 2|x| + \sqrt{1 + 4|x| - 4|x|^2}		&	0\leq x\leq 1,	\\
		5 - 2|x| - \sqrt{-7 + 12|x| - 4|x|^2}	&	1 < x \leq 2,	\\
		0										&	2 < x,
	\end{cases}
	\label{stokeseq:standard_ib_delta}
\end{equation}
as well as those produced using the smoother and more accurate $\delta$ function introduced in \cite{Stein2015} and used in all IBSE simulations in this paper. These two methods will be denoted by \IBa$\ $and \IBb, respectively. We solve the general Stokes problem with $\mathcal{L}=\mathbb{I}-\Delta$ and Dirichlet boundary conditions:
\begin{subequations}
	\label{stokeseq:manufactured_problem}
	\begin{align}
		\u - \Delta \u + \grad p	&=	\f_\u	&&\text{in }\Omega,	\\
		\grad\cdot\u				&=	0	&&\text{in }\Omega,	\\
		\u							&=	\b	&&\text{on }\Gamma,
	\end{align}
\end{subequations}
where the solution $\U$ is given by
\begin{subequations}
	\label{stokeseq:manufactured_solutions}
	\begin{align}
		u(x, y)	&=	e^{\sin x}\cos y,			\\
		v(x, y)	&=	-\cos x\, e^{\sin x}\sin y,	\\
		p(x, y)	&=	e^{\cos 2x},
	\end{align}
\end{subequations}
the forcing function $\f_\u$ is defined as
\begin{equation}
	\f_\u =
	\begin{cases}
		\u-\Delta\u + \grad p	&\text{in }\Omega,\\
		0											&\text{in }E,
	\end{cases}
\end{equation}
and the boundary condition $\b$ is found by evaluating $\u$ on $\Gamma$. The physical domain $\Omega$ is taken to be $\mathbb{T}^2\setminus \overline{B_1(2,2)}$, where $\mathbb{T}^2$ denotes the 2-torus identified with the periodic rectangle $[0,2\pi]\times[0,2\pi]$, and where $B_1(2,2)$ denotes the ball of radius $1$ centered at $(2,2)$. The computational domain is taken to be $C=\mathbb{T}^2$. Plots of the analytic solution $\U$ in the physical domain $\Omega$ are shown in \Cref{stokesfig:manufactured_analytic_solution}.

For this problem, it is convenient to use a Fourier spectral discretization for $C$ and to choose $C_E=C$. We discretize the domain $[0,2\pi]\times[0,2\pi]$ using a regular Cartesian mesh with $n$ points discretizing each dimension: $\Delta x=2\pi/n$ and $x_j=j\Delta x$ for $0\leq j < n$. Differential operators are discretized in the usual way.

\begin{figure}
	\centering
	\begin{subfigure}[b]{0.48\textwidth}
		\centering
		\includegraphics[width=\textwidth]{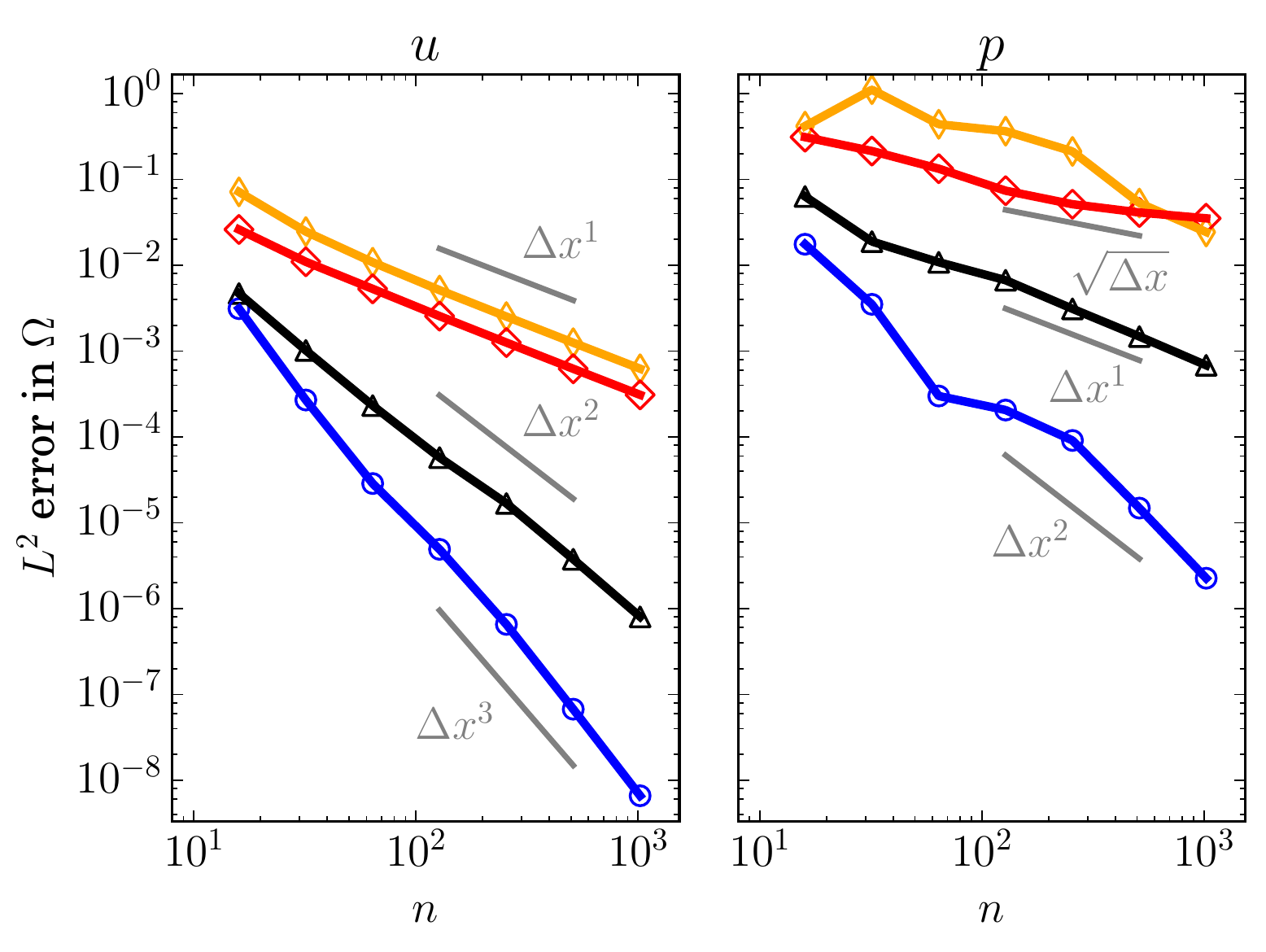}
		\subcaption{$L^2(\Omega)$ Error}
		\label{stokesfig:manufactured_error_solutions:l2}
	\end{subfigure}
	\begin{subfigure}[b]{0.48\textwidth}
		\centering
		\includegraphics[width=\textwidth]{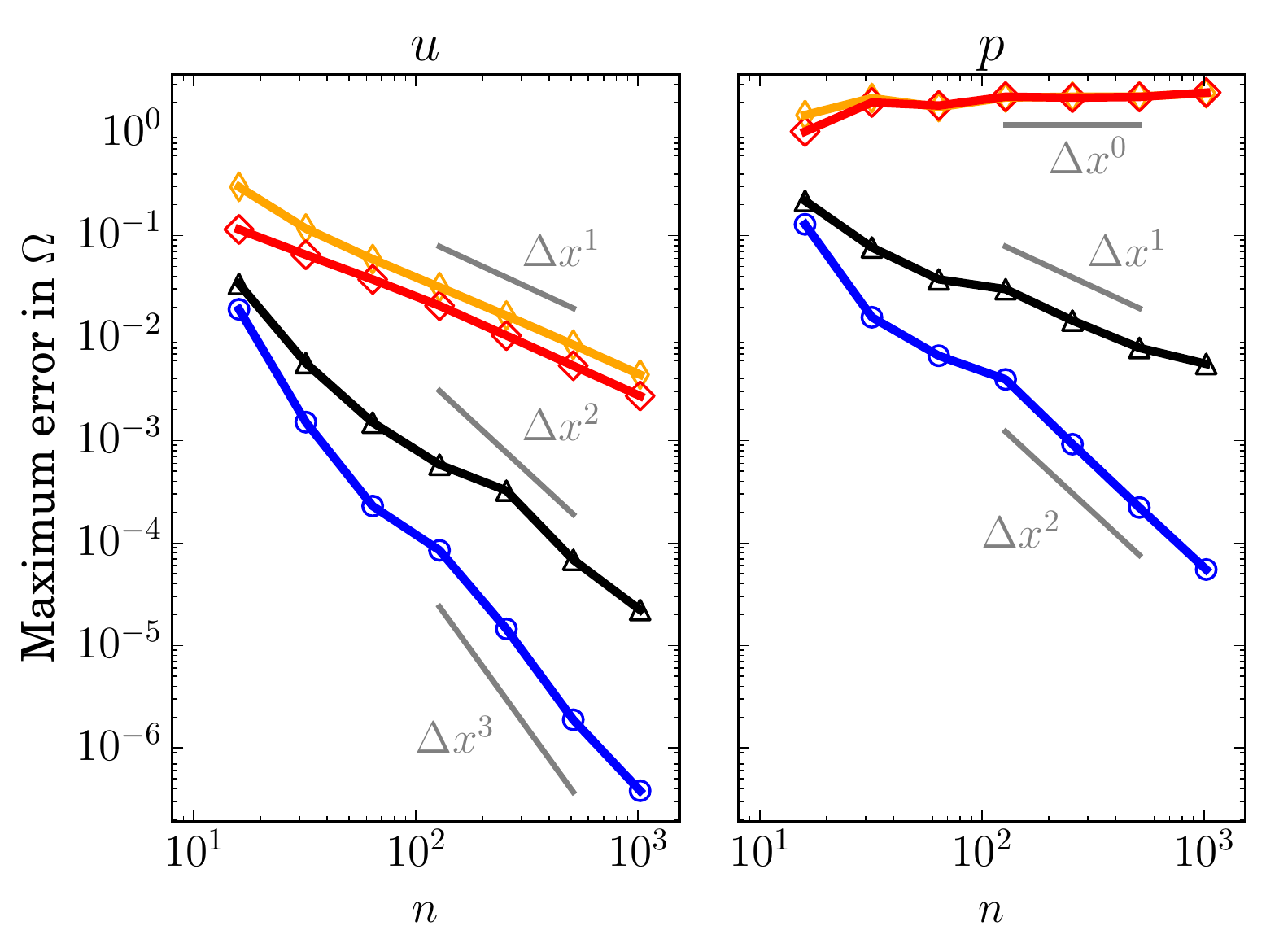}
		\subcaption{$L^\infty(\Omega)$ Error}
		\label{stokesfig:manufactured_error_solutions:max}
	\end{subfigure}
	\caption[A refinement study showing the convergence of the solutions to a Stokes problem, produced by the Immersed Boundary and Immersed Boundary Smooth Extension methods.]{$L^2(\Omega)$ and $L^\infty(\Omega)$ errors for the solution to \Cref{stokeseq:manufactured_problem} for the \IBa ({\color{py_orange} $\lozenge$}), \IBb ({\color{red}\large $\diamond$}), IBSE-$1$ ({{\color{black} $\triangle$}), and IBSE-$2$ ({\color{py_blue} \Circle}}) methods.}
	\label{stokesfig:manufactured_error_solutions}
\end{figure}

A refinement study for the $L^2(\Omega)$ and $L^\infty(\Omega)$ errors in the solutions for $u$ and $p$ is shown in \Cref{stokesfig:manufactured_error_solutions}. The refinement path for $v$ is omitted; it is similar to what is shown for $u$. In $L^2(\Omega)$, the velocity field $u$ converges at first-, second-, and third-order accuracy in $\Delta x$ for the IB, IBSE-$1$ and IBSE-$2$ methods, respectively, while the pressure field $p$ converges at first- and second-order accuracy for IBSE-$1$ and IBSE-$2$, but at the slow rate of $\mathcal{O}(\sqrt{\Delta x})$ for the IB method. In $L^\infty(\Omega)$, the rates of convergence for the velocity remain unchanged in all methods, while for the pressure the difference is starker: the IBSE-$1$ and IBSE-$2$ methods maintain first- and second-order convergence, while the IB method \emph{fails to converge pointwise.}

\begin{remark}
	It may be surprising that the velocity fields produced by the Immersed Boundary method converge at only first-order in $L^2(\Omega)$. Indeed, for a given $\F$ and $\b$, the velocity fields produced by solving
	\begin{subequations}
		\begin{align}
			-\Delta \u + \grad p + S\F	&=	0,	\\
			\grad\cdot\u								&=	0,	\\
			S^*\u												&=	\b,
		\end{align}
	\end{subequations}
	converge at $\mathcal{O}(\Delta x)$, $\mathcal{O}(\Delta x^{3/2})$, and $\mathcal{O}(\Delta x^2)$ in the $L^\infty(\Omega)$, $L^2(\Omega)$, and $L^1(\Omega)$ norms, respectively \cite{Mori2008}. However, $\F$ is not given, it is computed by inverting a Schur-complement equation whose elements include terms of the form $S^*\u$ where $\u\in C^0(C)$, and thus contain $\mathcal{O}(\Delta x)$ errors.
\end{remark}

\begin{figure}
	\centering
	\begin{subfigure}[b]{0.48\textwidth}
		\centering
		\includegraphics[width=\textwidth]{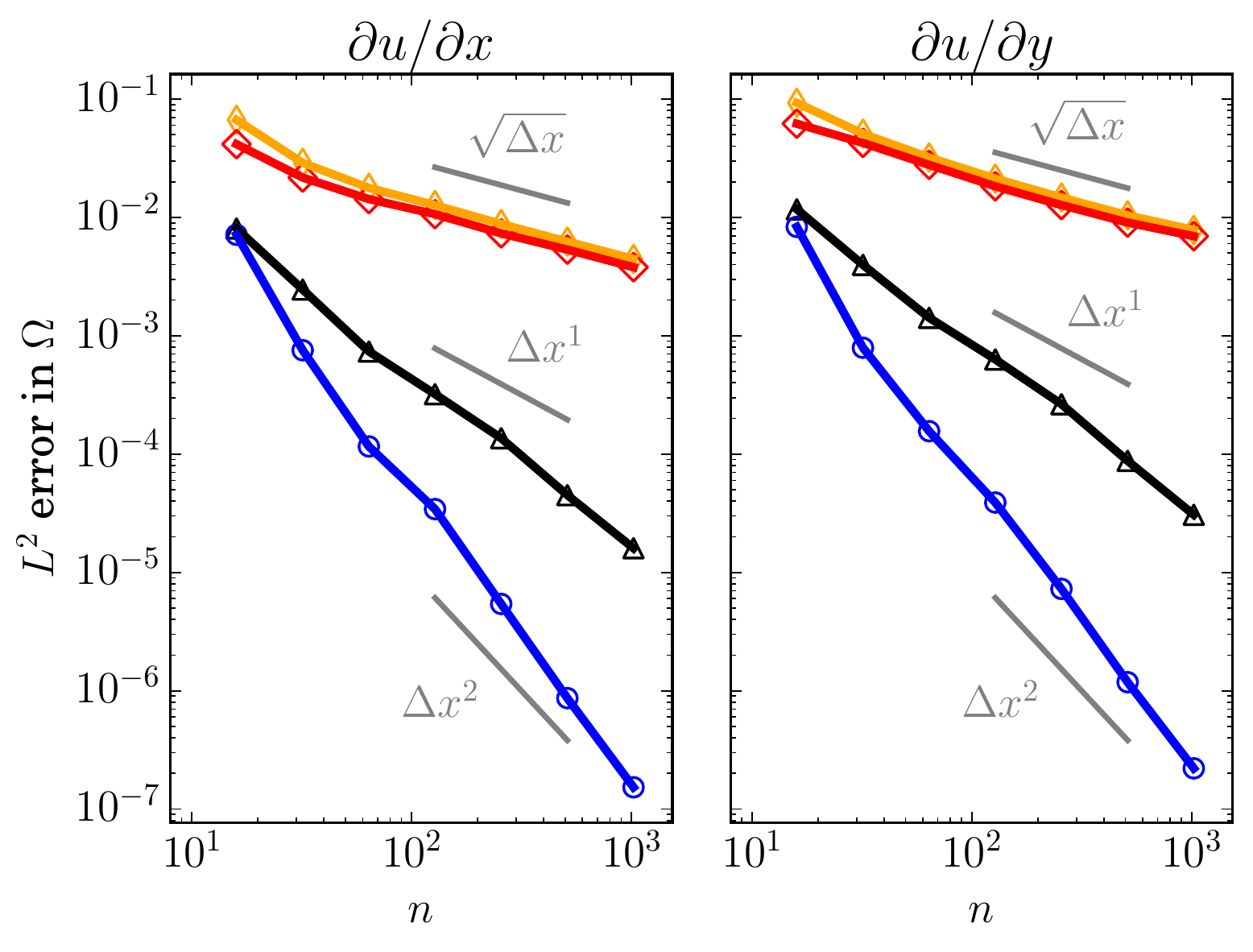}
		\subcaption{$L^2(\Omega)$ Error}
		\label{stokesfig:manufactured_error_derivatives:l2}
	\end{subfigure}
	\begin{subfigure}[b]{0.48\textwidth}
		\centering
		\includegraphics[width=\textwidth]{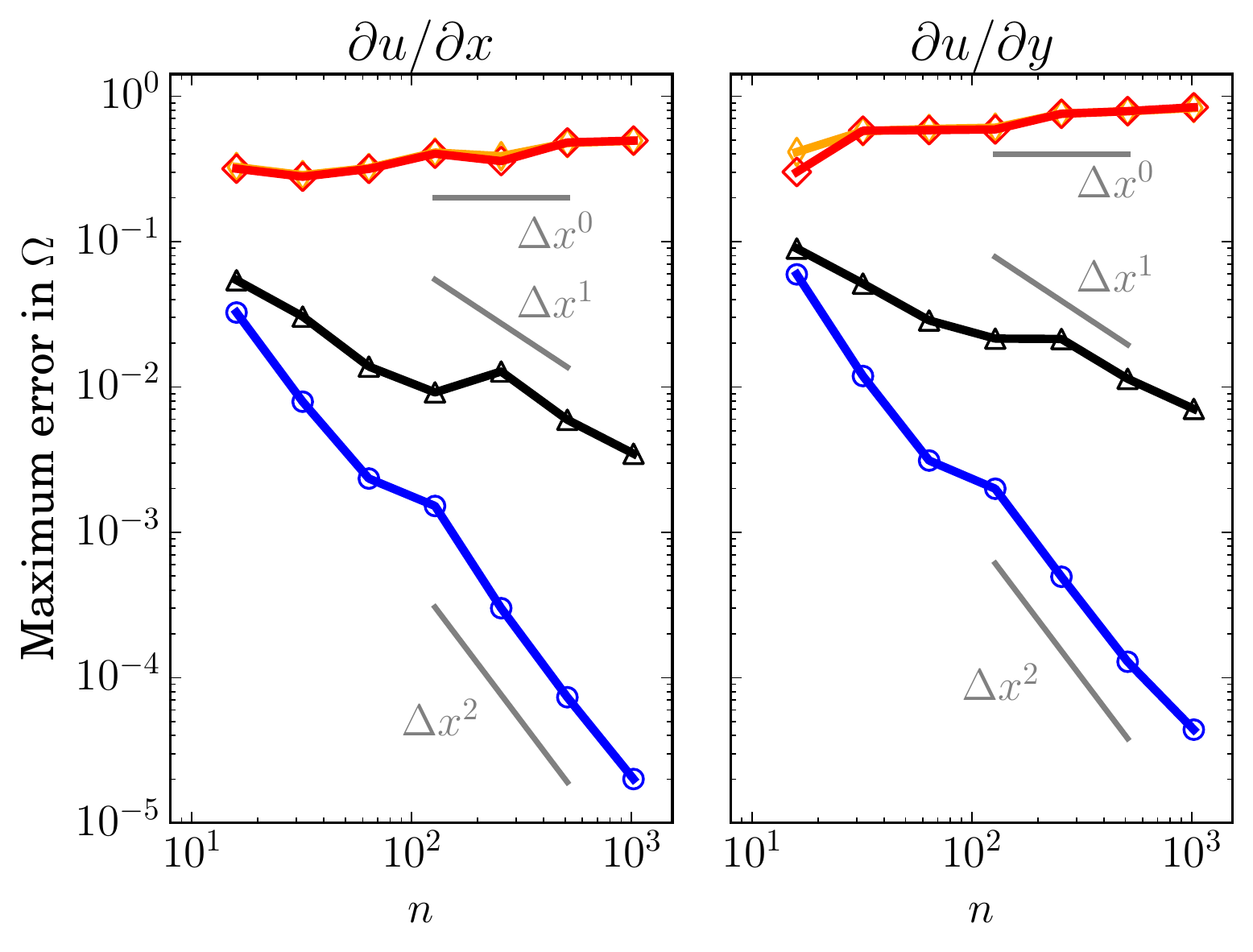}
		\subcaption{$L^\infty(\Omega)$ Error}
		\label{stokesfig:manufactured_error_derivatives:max}
	\end{subfigure}
	\caption[A refinement study showing the convergence of the derivatives of solutions to a Stokes problem, produced by the Immersed Boundary and Immersed Boundary Smooth Extension methods.]{$L^2(\Omega)$ and $L^\infty(\Omega)$ errors for the derivatives of the velocity fields that solve \Cref{stokeseq:manufactured_problem}, as computed by the \IBa ({\color{py_orange} $\lozenge$}), \IBb ({\color{red}\large $\diamond$}), IBSE-$1$ ({{\color{black} $\triangle$}), and IBSE-$2$ ({\color{py_blue} \Circle}}) methods.}
	\label{stokesfig:manufactured_error_derivatives}
\end{figure}

In the refinement study shown in \Cref{stokesfig:manufactured_error_solutions}, the pressure function converges more slowly than the velocity function; in the case of the IB method, it fails to converge pointwise. This will be true in general for elements of the stress tensor, which are one derivative less smooth than the velocity field. In \Cref{stokesfig:manufactured_error_derivatives}, we show a refinement study for the $L^2(\Omega)$ and $L^\infty(\Omega)$ errors in the derivatives $\partial u/\partial x$ and $\partial u/\partial y$ (the convergence paths for $\partial v/\partial x$ and $\partial v/\partial y$ are omitted, but essentially the same). A similar pattern to the convergence of the pressure function is observed: for all of the derivatives, the IB method converges slowly (at $\mathcal{O}(\sqrt{x})$) in $L^2(\Omega)$, but fails to converge pointwise; the IBSE-$1$ and IBSE-$2$ method provide first- and second-order convergence in both $L^\infty(\Omega)$ and $L^2(\Omega)$. We emphasize that with pointwise convergence of the pressure function $p$ as well as pointwise convergence of these derivatives, the IBSE method is able to capture the stress tensor $\sigma=\nu(\grad\u + \grad\u^\intercal) - p\mathbb{I}$ \emph{pointwise up to the boundary $\Gamma$.}

\begin{figure}
	\centering
	\hspace*{\fill}
	\begin{subfigure}[b]{0.25\textwidth}
		\centering
		\includegraphics[width=\textwidth]{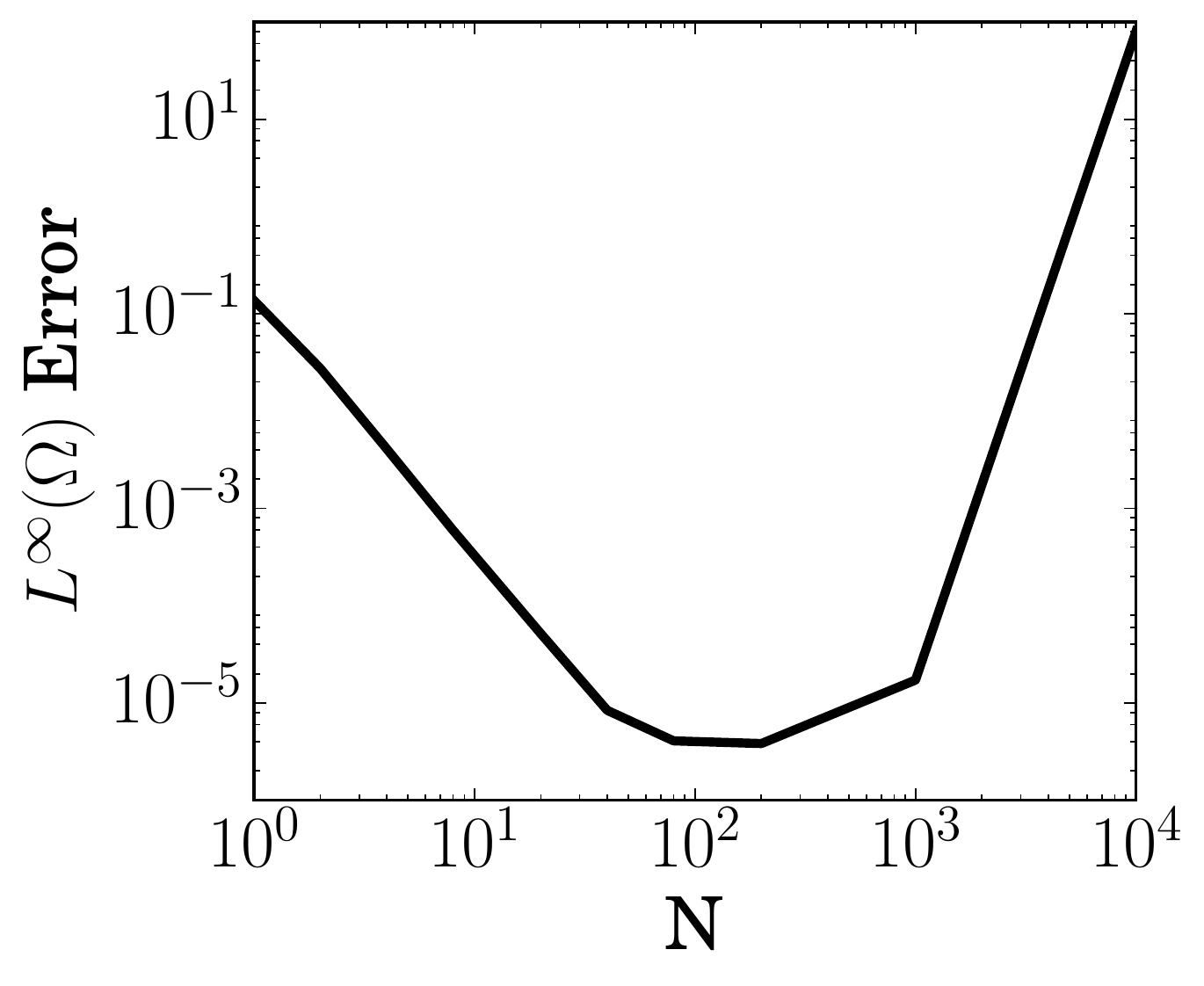}
		\caption{$L^\infty(\Omega)$ error in $u$}
		\label{stokesfig:manufactured_extension:error}
	\end{subfigure}
	\hfill
	\begin{subfigure}[b]{0.25\textwidth}
		\centering
		\includegraphics[width=\textwidth]{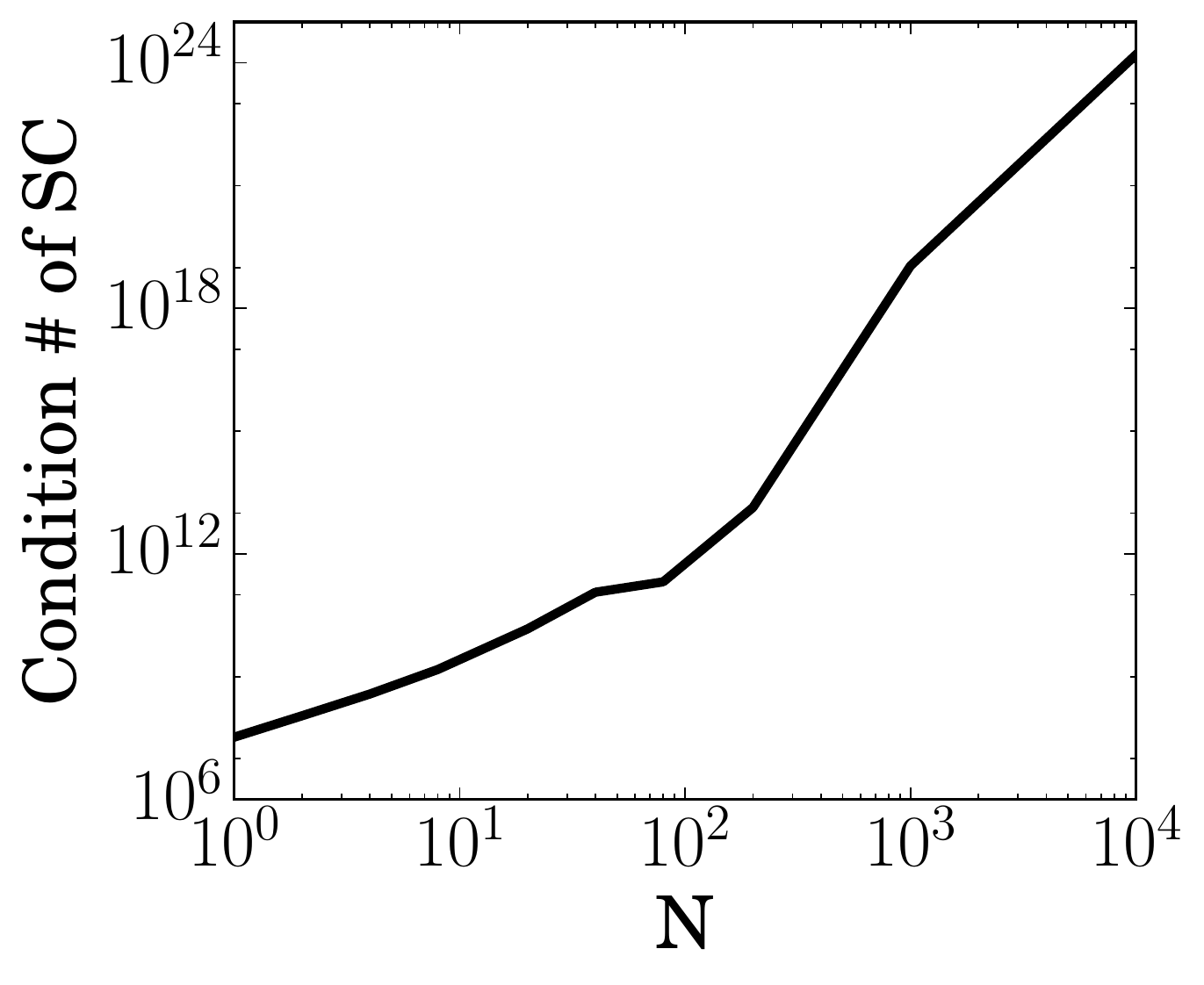}
		\subcaption{Condition \# of SC}
		\label{stokesfig:manufactured_extension:condition}
	\end{subfigure}
	\hfill
	\begin{subfigure}[b]{0.22\textwidth}
		\centering
		\includegraphics[width=\textwidth]{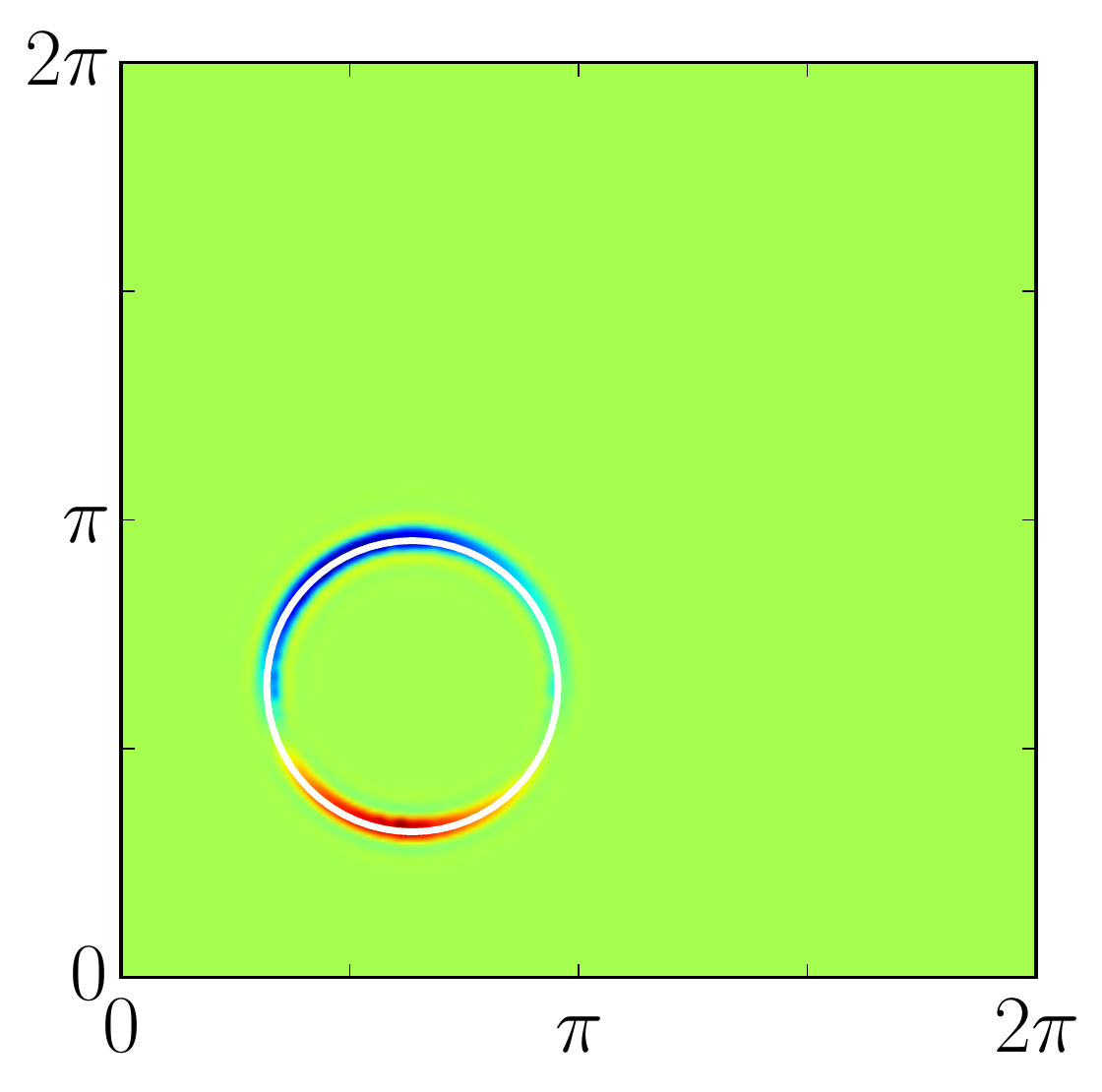}
		\subcaption{$N=1$}
		\label{stokesfig:manufactured_extension:1}
	\end{subfigure}
	\hfill
	\begin{subfigure}[b]{0.22\textwidth}
		\centering
		\includegraphics[width=\textwidth]{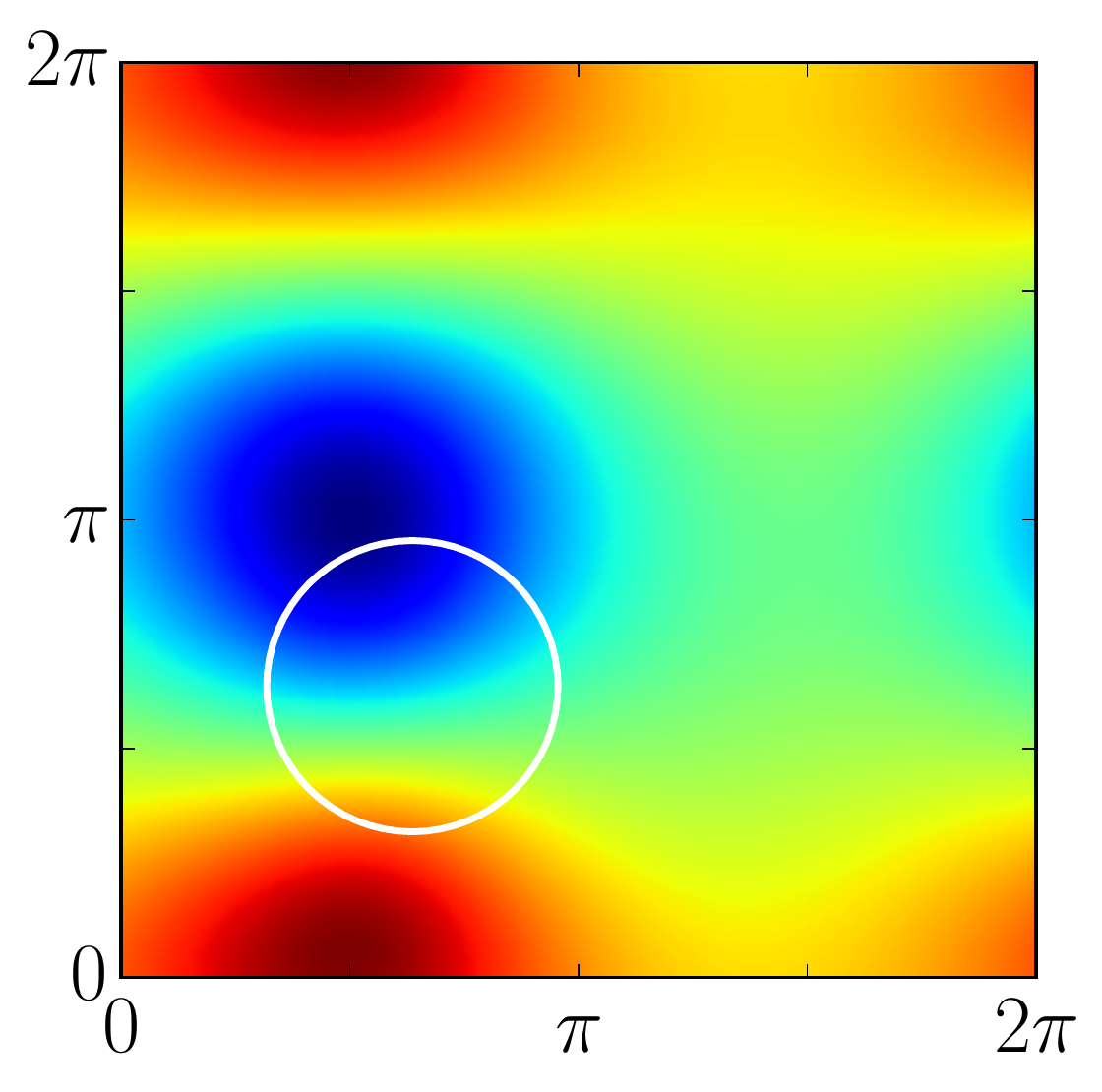}
		\subcaption{$N=200$}
		\label{stokesfig:manufactured_extension:200}
	\end{subfigure}
	\hspace*{\fill}
	\caption[The effect of the regularization parameter for the extension operator $\mathcal{H}^k$ on the error in the solution and on the extension function.]{\Cref{stokesfig:manufactured_extension:error,stokesfig:manufactured_extension:condition} show the $L^\infty(\Omega)$ error in $u$ for the solution to \Cref{stokeseq:manufactured_problem} and the condition number of the Schur Complement defined by \Cref{stokeseq:define_schur_complement} produced using different values of $N$. \Cref{stokesfig:manufactured_extension:1,stokesfig:manufactured_extension:200} show the extension function $\xi_u$ to $u$ for $N=1$ and $N=200$ (the value used in all simulations in this paper, other than this test). When $N$ is too small (e.g. $N=1$) short length scales not fully resolvable by the discretization are introduced to the problem.}
	\label{stokesfig:manufactured_solutions_extensions_by_theta}
\end{figure}

Finally, we use this test problem to numerically investigate the choice of the free parameter $N$ in the definition of the extension operator $\mathcal{H}^k$. This parameter controls the rate at which the extension functions $\bxi_\u$ and $\xi_p$ decay to $0$ away from the boundary. We compute solutions with $n=256$, using extensions across a wide range of choices of $N$ (from $N=1$ to $N=10000$). The $L^\infty(\Omega)$ error for $u$, and the condition number of the Schur complement, as a function of $N$, are shown in \Cref{stokesfig:manufactured_extension:error,stokesfig:manufactured_extension:condition}. For very small values of $N$, the error is dominated by the inability of the discretization to resolve the small length scale introduced to the problem. As $N$ is increased, the error produced by the IBSE method improves, but the condition number of the Schur-complement worsens. As $N$ grows large, the error again begins to degrade, dominated by the error produced when inverting the ill-conditioned Schur-complement. This presents a balance between resolution and conditioning. For this test, the minimum error is somewhere between $N=100$ and $N=200$, a result that has been fairly robust across the range of test problems we have examined. Other than for this test, all numerical results in this paper are produced using with $N=200$ to construct $\mathcal{H}^k$. In \Cref{stokesfig:manufactured_extension:1,stokesfig:manufactured_extension:200}, we show the extension function $\xi_u$ for the solution $u$ to \Cref{stokeseq:manufactured_problem} for $N=1$ and $N=200$.

\begin{remark}
	In this section, we presented results for the Immersed Boundary method using both the \IBa$\ $and \IBb$\ $functions to discretize the $S$ and $S^*$ operators. For all results presented in this section, we observe the same asymptotic behavior, although the \IBb$\ $method typically provides lower errors. For the remainder of this paper, we will present results for the Immersed Boundary method using only the more commonly used \IBa$\ $function.
\end{remark}

\section{Stokes flow around a cylinder}
\label{stokessection:stokes_test:flow_around_cylinder}

In this section we focus on the zero Reynolds number flow around a cylinder in a confined channel. The exact problem we will simulate is on the domain $[-20,20]\times[-2,2]$, with a cylinder with radius $R=1$ centered at $(0,0)$, and a mean stream-wise flow velocity of $u=1$ enforced at the inflow boundary (where $x=-20$). No-slip ($\u=0$) boundary conditions are imposed on the channel and cylinder walls. We choose this setup as this domain has been extensively studied in the context of polymeric flow problems \cite{Dou2006,Alves2001,Claus2013}, and accurate benchmarks are available.

We solve this problem using two separate underlying discretizations. In \Cref{stokessection:stokes_test:flow_around_cylinder:spectral}, we describe a Fourier spectral discretization. The Fourier discretization allows the method to achieve third-order accuracy for the velocity field $\u$ when $C^2$ extensions are used, but when using this discretization, the enforcement of the no-slip $\u=0$ boundary condition along the long channel walls requires the use of the IBSE methodology. This leads to a large number of boundary points $n_\text{bdy}$ and a very large Schur-complement operator that must be formed and inverted.

To help alleviate this problem, in \Cref{stokessection:stokes_test:flow_around_cylinder:finite_difference} we couple the IBSE method to a standard second-order staggered-grid finite difference discretization of the Stokes equations. The no-slip condition on the channel walls is enforced naturally by the discretization, and the IBSE methodology is used only to enforce the no-slip condition along the cylinder wall and to smoothly extend the unknown solution to the interior of the cylinder. Despite the discretization being limited to second-order accuracy, we show that there is still a benefit to using the IBSE-$2$ method over the IBSE-$1$ method, especially when local information regarding the stress is required. If higher-order accuracy is desired, a more accurate finite-difference discretization could be employed.

In \Cref{stokessection:stokes_test:flow_around_cylinder:results}, we compare the accuracy of these two different discretizations, for the \IBa, IBSE-$1$, and IBSE-$2$ methods, and finally in \Cref{stokessection:stokes:inefficiences}, we compare the numerical efficiency of the different discretizations.

\subsection{Fourier spectral discretization}
\label{stokessection:stokes_test:flow_around_cylinder:spectral}

\begin{figure}[htb!]
	\centering
	\includegraphics[width=\textwidth]{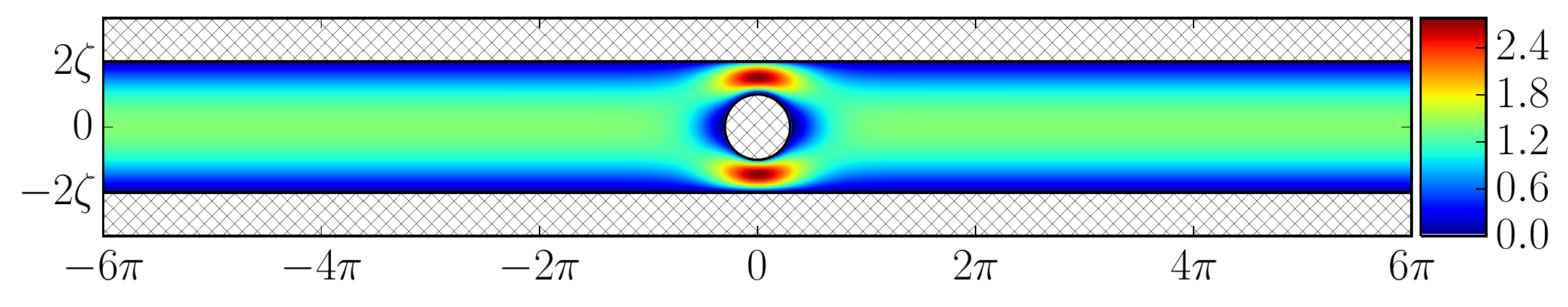}
	\caption[The streamwise velocity and domain for the confined flow around a cylinder test problem used to validate the accuracy of the Immersed Boundary Smooth Extension method for Stokes flow.]{The streamwise velocity $u$ and solution domain for the Fourier spectral discretization to the Stokes channel flow around a cylinder problem defined in \Cref{stokessection:stokes_test:flow_around_cylinder:spectral}. The extension domain $E$ is denoted with crosshatches. The solution that is shown was produced using the IBSE-$2$ method with $n_y=512$. Note that space is rescaled and $\zeta$ is chosen so that this domain is equivalent to the $[-20,20]\times[-2,2]$ domain used in other studies.}
	\label{stokesfig:spectral_stokes_solution}
\end{figure}

For this problem we take the computational domain to be $C=[-6\pi, 6\pi]\times[-\pi,\pi]$, with channel walls located at $y=\pm2\zeta$, where $\zeta=12\pi/40$; and the cylinder of radius $R=\zeta$ is centered at $(0,0)$. This domain, together with the streamwise velocity $u$ produced using the IBSE-$2$ method with $n_y=512$, is shown in \Cref{stokesfig:spectral_stokes_solution}. The domain is discretized with a uniform Cartesian mesh, with $n_y$ points in the spanwise direction and $n_x=6 n_y$ points in the streamwise direction, with $\Delta x=2\pi/n_y$. The discrete meshpoints ${x_{ij}}$ are located at $x_{ij}=(-6\pi+i\Delta x,-\pi+j\Delta x)$ for $0\leq i<n_x$ and $0\leq j<n_y$. All differential operators are discretized in the usual way, and the spread and interpolation operators are discretized as described in \Cref{stokessection:numerics:operators}. Because the extension region $E$ (denoted by crosshatches in \Cref{stokesfig:spectral_stokes_solution}) is not localized to a small region of $C$, we choose the domain on which the extension equations are solved to coincide with $C$, i.e. $C_E=C$. In order to simulate a (rescaled) Poiseuille inflow condition, we solve the problem:
\begin{subequations}
	\label{stokeseq:stokes_equations_spectral}
	\begin{align}
		-\Delta\u + \grad p	&=	\beta\hat x	&	&\text{in }\Omega,	\\
		\grad\cdot\u		&=	0	&	&\text{in }\Omega,	\\
		\u					&=	0	&	&\text{on }\Gamma,	\\
		\int_{-2\zeta}^{2\zeta} u(-6\pi, y)\, dy	&= \zeta, \label{stokeseq:average_velocity_constraint}
	\end{align}
\end{subequations}
where $\Gamma$ denotes the internal boundaries (the surfaces of the cylinder and channel walls). The magnitude of the constant forcing $\beta$ is treated as a Lagrange multiplier, and is used to enforce the constraint that the spanwise average of the streamwise velocity ($u$) at the far left end of the domain is $\zeta$. The integral across the channel in \Cref{stokeseq:average_velocity_constraint} is discretized using Simpson's rule. In \cite{Dou2006,Alves2001,Claus2013}, the boundary condition specified at the inflow is that $u(-20, y) = 3(4-y)^2/8$. Our computational setup approximates this boundary condition but does not impose it exactly. Given the length of the channel compared to $R$, we expect our solutions to match closely with other studies; results are shown in \Cref{stokessection:stokes_test:flow_around_cylinder:results}. 

\subsection{Second-order finite difference discretization}
\label{stokessection:stokes_test:flow_around_cylinder:finite_difference}


In this section we describe how to couple the IBSE method to a standard finite difference scheme. The Stokes Equations \eqref{stokeseq:general_brinkmann_problem} are discretized using a second-order finite difference staggered grid discretization where the velocity unknowns ($u$, $v$) are stored on the vertical and horizontal cell edges, respectively, and the pressure unknowns $p$ are stored at the cell centers \cite{welch1965mac,mckee2008mac}. For this discretization, the computational domain is taken to be $C=[-20,20]\times[-4,4]$. There is no need to provide an extension region along the channel walls, as in the domain for the Fourier spectral discretization shown in \Cref{stokesfig:spectral_stokes_solution}. No-slip conditions ($\u=0$) are applied on the top and bottom walls at $y=-2$ and $y=2$; an inflow condition of $\u=(3(4-y^2)/8,0)$ is applied at $x=-20$; and an approximate outflow condition of $\partial\u/\partial x$ is applied at $x=20$. Note that for this discretization, the flow is driven by the boundary condition, and no additional forcing has to be added to the momentum equation as in \Cref{stokeseq:stokes_equations_spectral}. Inversion of the discrete Stokes system (\Cref{stokeseq:discrete_stokes_system}) is accomplished using a preconditioned GMRES method with a projection preconditioner \cite{Cai2013}. With the operator $\mathcal{L}=\theta\mathbb{I}-\nu\Delta_\u$, this preconditioner is given by
\begin{equation}
	P^{-1}	=	
	\begin{pmatrix}
		\mathbb{I}	&	\grad\Delta_p^{-1}	\\
		0			&	-\theta\Delta_p^{-1} + 2\nu\mathbb{I}
	\end{pmatrix}
	\begin{pmatrix}
		\mathbb{I}	&	0	\\
		-\grad\cdot	&	\mathbb{I}
	\end{pmatrix}
	\begin{pmatrix}
		\mathcal{L}^{-1}	&	0	\\
		0					&	\mathbb{I}
	\end{pmatrix}.
\end{equation}
The discrete operators $\mathcal{L}$ and $\Delta_p$ are inverted exactly using a combination of fast sine/cosine transforms in one direction and tri-diagonal matrix inverses in the other direction. In our numerical tests this preconditioner provides convergence to a relative residual of $10^{-10}$, independent of the grid spacing $\Delta x$, in approximately 30-40 iterations for the Stokes problem studied in this section, and approximately 3-4 iterations for the Navier-Stokes problems studied in \Cref{stokessection:navier_stokes:unsteady} and \Cref{stokessection:navier_stokes:steady}.

In principle, coupling the IBSE method to this discretization is simple. Examining the outline given in \Cref{stokessection:numerics:summary}, the only details that are not fully specified are the choice of $C_E$ (the computational domain for the extension equations), and the rules for transferring $\bxi_\u$ and $\xi_p$ from the grid or grids on which the extensions $\bxi_\u$ and $\xi_p$ are solved on to the staggered grids for $\U$.

\begin{figure}
	\centering
	\includegraphics[width=0.9\textwidth]{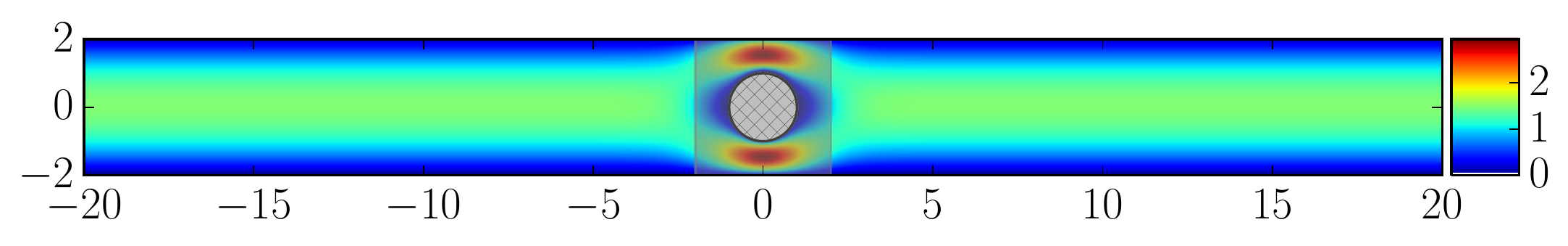}
	\caption[The streamwise velocity and domain (highlighting the localized extension domain) for the confined flow around a cylinder test problem used to validate the accuracy of the Immersed Boundary Smooth Extension, coupled to a finite difference solver, for a Stokes flow.]{The solution $u$ to the confined channel flow around a cylinder, produced using a second-order finite difference discretization to the Stokes equations as the underlying PDE solver for the IBSE-$2$ method with $n_y=512$, as described in \Cref{stokessection:stokes_test:flow_around_cylinder:finite_difference}. The extension region $E$ is denoted by crosshatches. The extension functions $\bxi_\u$ and $\xi_p$ are not solved for on the entirety of the computational domain, but rather in the much smaller region $C_E=[-2,2]\times[-2,2]$, localized around $E$, that is shaded in light gray.}
	\label{figure:fd_stokes_solution}
\end{figure}

In \Cref{figure:fd_stokes_solution}, we show the domain $C$ for this problem, together with the solution computed by the IBSE-$2$ method with $n_y=512$. Notice that the extension region $E$, denoted by crosshatches, is isolated to a small region of $C$. Solving the extension equations for $\bxi_\u$ and $\xi_p$ on all of $C$ would be computationally wasteful. Instead, $C_E$ is chosen for this problem to be $[-2,2]\times[-2,2]$, which is highlighted in gray in \Cref{figure:fd_stokes_solution}. In principle, the only requirement for $C_E$ is that it contain $E$; in practice, it should be chosen large enough so that the periodic boundary conditions used when solving the extension equation on $C_E$ do not induce any undesirably large gradients. The operators $\mathcal{H}^k$ and $\mathcal{H}^{k-1}$ are discretized and inverted using Fourier spectral methods.

Because the staggered discretization we use does not collocate the unknowns for the velocity fields $u$ and $v$ and the pressure field $p$, some choice must be made as to where the unknowns for the extension functions $\xi_u$, $\xi_v$, and $\xi_p$ should be located. We have explored two options:
\begin{enumerate}
	\item A single extension `grid' $C_E$, where $\xi_u$, $\xi_v$, and $\xi_p$ are collocated. Bilinear interpolation is used to transfer information between this grid and the staggered grids for $u$, $v$, and $p$.
	\item Three extension grids ${C_E}^u$, ${C_E}^v$, and ${C_E}^p$, on which the nodes for $\xi_u$, $\xi_v$, and $\xi_p$ are each collocated with the nodes for $u$, $v$, and $p$, respectively. No interpolation is necessary when transferring information between grids. The extension grid ${C_E}^u=[0,2]\times[\Delta x/2,2+\Delta x/2]$ is shifted from $[0,2]\times[0,2]$ by $\Delta x/2$ in the $y$ direction so that the locations of the unknowns coincide with the locations of the unknowns for the $u$ grid; the extension grid ${C_E}^v$ is shifted by $\Delta x/2$ in the $x$ direction, and the extension grid ${C_E}^p$ is shifted by $\Delta x/2$ in both the $x$ and $y$ directions.
\end{enumerate}

In our numerical experiments, the second option produced Schur-complement operators with slightly lower numerical condition numbers and improved stability when solving time-dependent equations. Although this choice requires slightly greater computational complexity, we use this option for this problem as well as the Navier-Stokes problems defined in \Cref{stokessection:navier_stokes}.

\subsection{Results: Stokes flow around a cylinder}
\label{stokessection:stokes_test:flow_around_cylinder:results}

\begin{figure}
	\centering
	\hspace*{\fill}
	\begin{subfigure}[b]{0.48\textwidth}
		\centering
		\includegraphics[width=\textwidth]{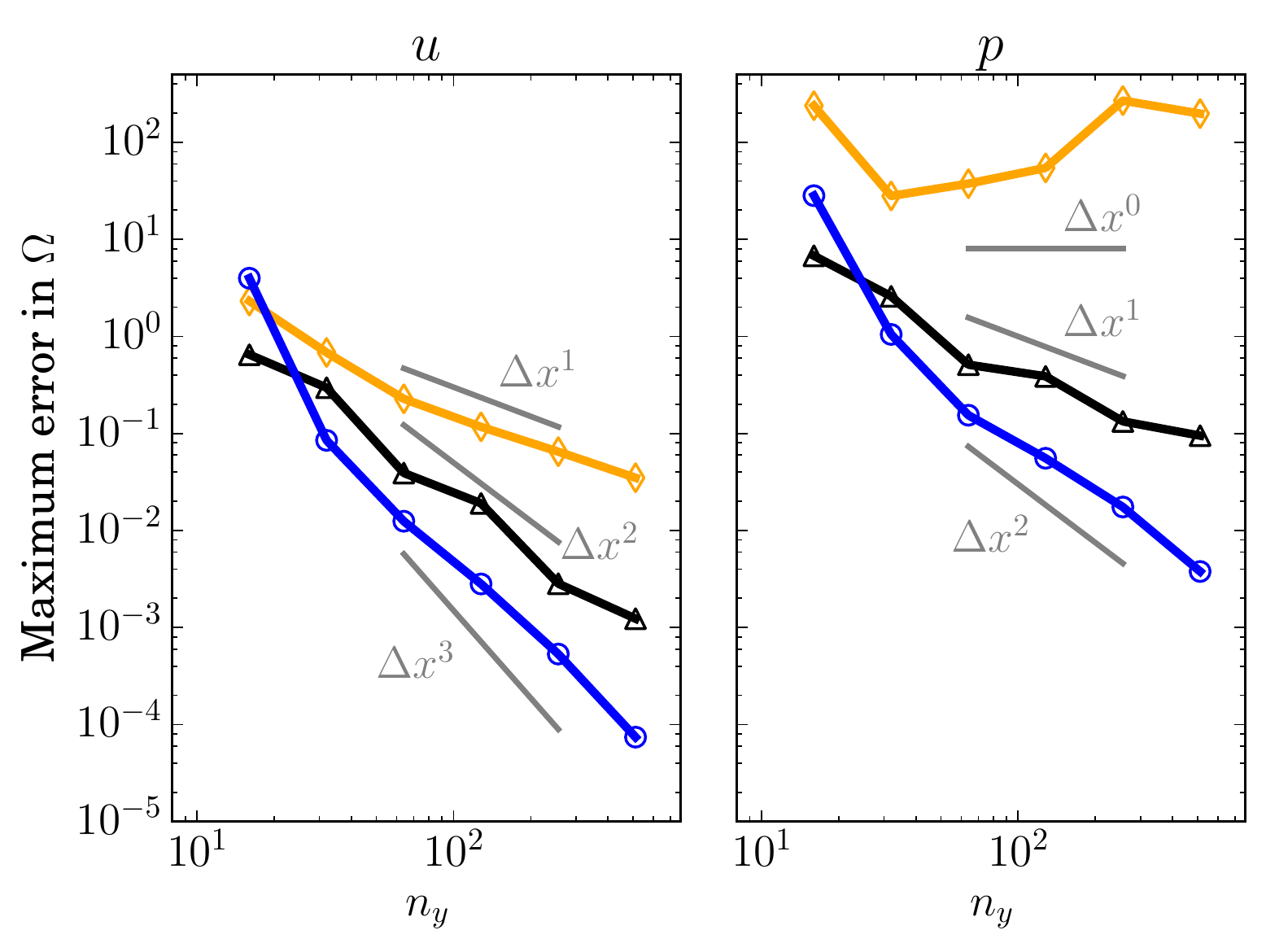}
		\caption{Spectral}
		\label{stokesfig:stokes_cylinder_error:spectral}
	\end{subfigure}
	\hfill
	\begin{subfigure}[b]{0.48\textwidth}
		\centering
		\includegraphics[width=\textwidth]{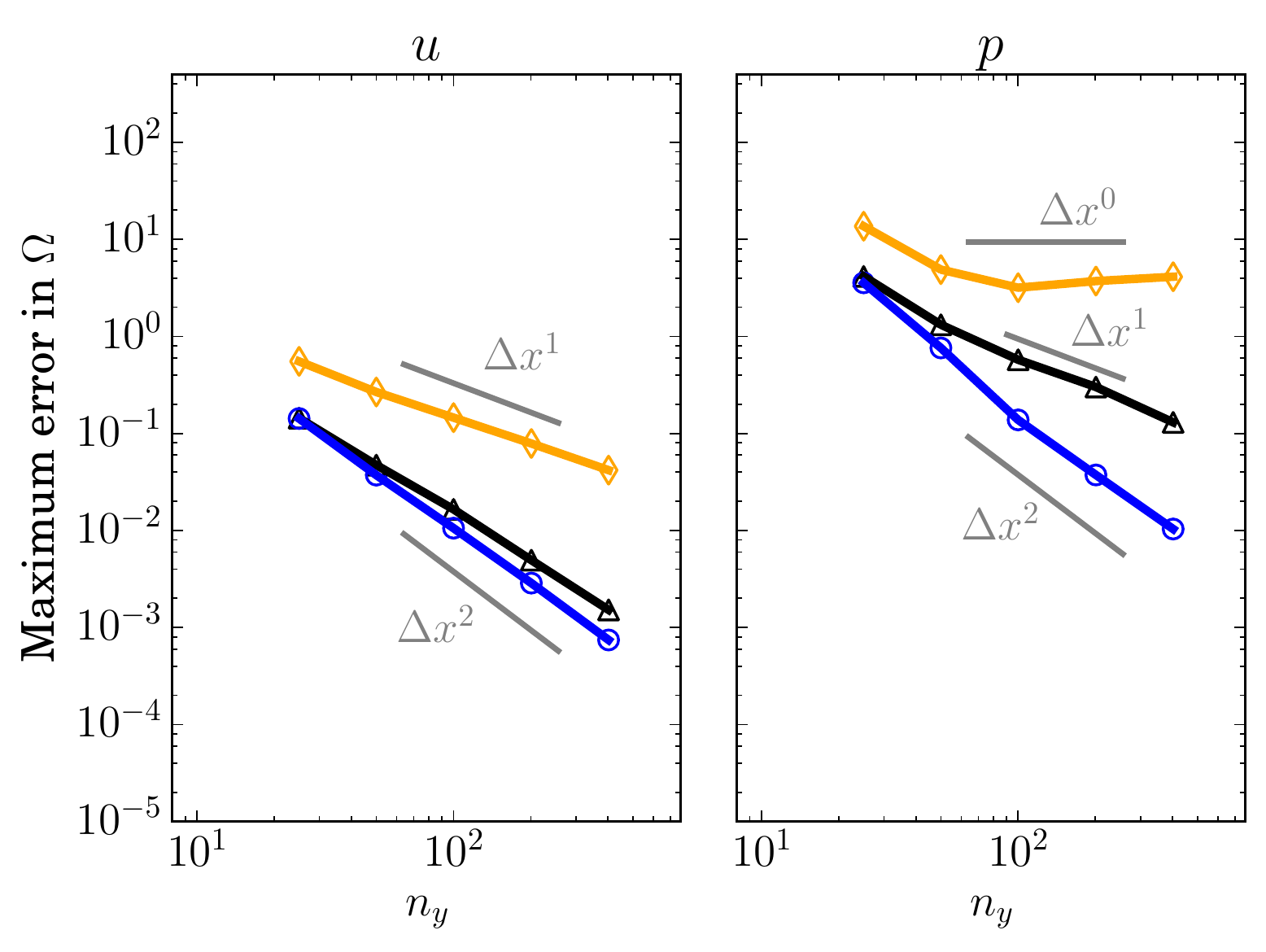}
		\subcaption{Finite difference}
		\label{stokesfig:stokes_cylinder_error:finite_difference}
	\end{subfigure}
	\hspace*{\fill}
	\caption[A refinement study validating the accuracy of the Immersed Boundary Smooth Extension method for Stokes flow on the confined channel flow around a cylinder problem.]{Refinement study showing the convergence of differences between successively refined solutions, in $L^\infty(\Omega)$, for the solutions $u$ and $p$ to the Stokes confined channel flow around a cylinder problem, computed by the \IBa ({\color{py_orange} $\lozenge$}), IBSE-$1$ ({{\color{black} $\triangle$}), and IBSE-$2$ ({\color{py_blue} \Circle}}) methods.. \Cref{stokesfig:stokes_cylinder_error:spectral} shows solutions produced using a Fourier spectral discretization, \Cref{stokesfig:stokes_cylinder_error:finite_difference} shows solutions produced using the finite difference discretization. The value of $n_y$ for the finite difference method has been scaled by a factor of $2\pi/4$, so that the $x$-axes represent the same effective resolution.}
	\label{stokesfig:stokes_cylinder_error}
\end{figure}

We now compare the results from the Fourier discretization defined in \Cref{stokessection:stokes_test:flow_around_cylinder:spectral} and the finite difference discretization defined in \Cref{stokessection:stokes_test:flow_around_cylinder:finite_difference}, across the \IBa, IBSE-$1$, and IBSE-$2$ methods. A convergence study showing the difference between successively refined solutions for $u$ and $p$ in $L^\infty(\Omega)$ is shown in \Cref{stokesfig:stokes_cylinder_error}. The same trends in convergence as observed in \Cref{stokesfig:manufactured_error_solutions:max} are observed for the Fourier spectral discretization. For the finite difference discretization, the $L^\infty(\Omega)$ convergence is limited by the accuracy of the underlying discretization to second-order. Note, however, that the IBSE-$1$ and IBSE-$2$ methods still differ in the convergence of the pressure function $p$: pointwise convergence of the pressure function is first-order for the IBSE-$1$ method and second-order for the IBSE-$2$ method. The convergence rates for all elements of the stress tensor are the same as for the pressure.

The drag coefficient is a simple, commonly used benchmark for comparing solutions across methods, and is computed as
\begin{equation}
	C_D	= \frac{2F_D}{\overline{U}D},
\end{equation}
where $\overline{U}$ is the average velocity of the fluid across the inflow boundary, $D$ is the diameter of the cylinder, and $F_D$ is the drag force:
\begin{equation}
	F_D = \int_\Gamma \sigma\cdot \n\cdot\hat x\,ds.
\end{equation}
For IBSE-$1$ and IBSE-$2$, the elements of the stress tensor are captured accurately up to the boundary. To compute $\sigma$, we simply compute the appropriate derivatives of $\u$ (spectrally or with second-order finite differences \emph{ignoring the boundary $\Gamma$}), and interpolate to the boundary by applying the operator $S^*$. For the spectral method, $\beta x$ must be added to the pressure, as the constant body forcing $\beta\hat x$ that is applied to drive the flow acts as an effective pressure gradient. For the IB method, $\u$ is only continuous, and differentiation produces large errors at the interface. The drag force can, however, be found by integrating the singular forces supported along the boundary \cite{Lai2000}:
\begin{equation}
	F_D = \int_\Gamma \G(s)\cdot\hat{x}\, ds,
\end{equation}
where $\G$ are the singular forces from \Cref{stokeseq:immersed_boundary}.

\begin{figure}[htb!]
	\centering
	\hspace*{\fill}
	\begin{subfigure}[b]{0.35\textwidth}
		\centering
		\includegraphics[width=\textwidth]{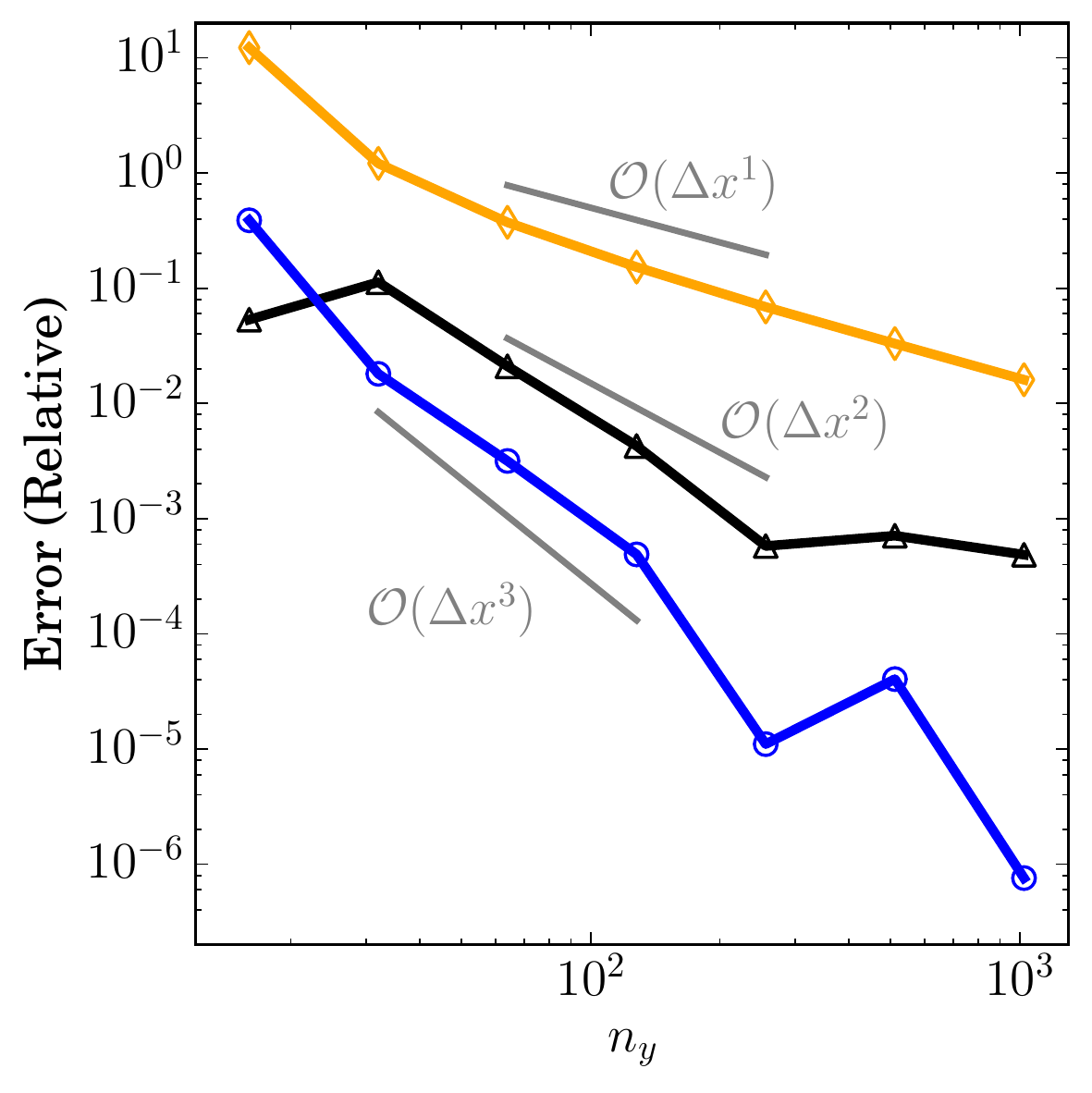}
		\caption{Spectral}
		\label{stokesfig:stokes_drag_error:spectral}
	\end{subfigure}
	\hfill
	\begin{subfigure}[b]{0.35\textwidth}
		\centering
		\includegraphics[width=\textwidth]{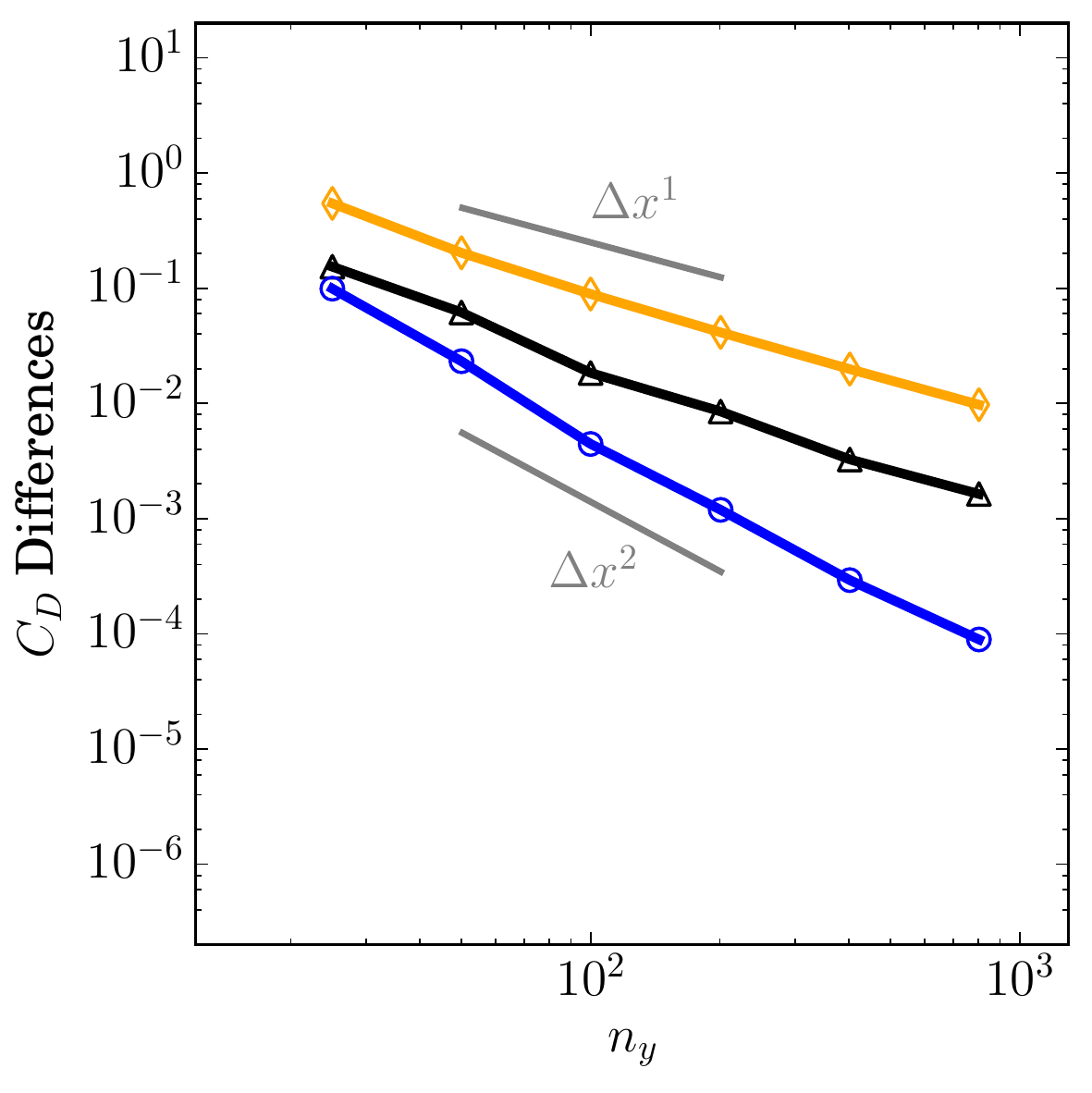}
		\subcaption{Finite difference}
		\label{stokesfig:stokes_drag_error:fd}
	\end{subfigure}
	\hspace*{\fill}
	\caption[Convergence of the drag coefficient for the confined channel flow around a cylinder problem, at Reynolds number $0$ (Stokes flow).]{Refinement study showing the convergence of the drag coefficient to the value $132.36$ from \cite{fan1999galerkin}, produced by the \IBa ({\color{py_orange} $\lozenge$}), IBSE-$1$ ({{\color{black} $\triangle$}), and IBSE-$2$ ({\color{py_blue} \Circle}}) methods. The error is normalized by the reference value. The value of $n_y$ for the finite difference method has been scaled by a factor of $2\pi/4$, so that the $x$-axes represent the same effective resolution.}
	\label{stokesfig:stokes_drag_error}
\end{figure}

In \Cref{stokesfig:stokes_drag_error}, we show a refinement study for the drag coefficient, as compared to the reference value of $132.36$, from \cite{fan1999galerkin}\footnote{This geometry is typically used to study polymeric fluids. The most accurate results currently available are from a recent high-order hp spectral element method \cite{Claus2013}, but a result for the Newtonian case is not provided in that study. The value that we use for the Newtonian case comes from an older hp finite element method \cite{fan1999galerkin}. At low elasticity where both \cite{fan1999galerkin} and \cite{Claus2013} report results, the drag coefficients reported agree to at least two decimal places.}. The IB method achieves first-order convergence of the drag coefficient in both methods;  the IBSE-$1$ and IBSE-$2$ methods achieve first- and second-order convergence with the finite difference discretization and convergence that is at least first- and second-order with the Fourier spectral method. It is surprising that the drag coefficient for the IB method converges at all: none of the elements of the stress tensor converge pointwise up to the boundary. Nevertheless, the singular forces accurately represent the stress the fluid is exerting on the solid. Individual elements of the stress tensor, however, do not converge near to the boundary. For the IBSE-$k$ method, the velocity field is $C^k(C)$, while the elements of the stress tensor are $C^{k-1}(C)$. For functions with $C^{k-1}(C)$ regularity, the operator $S^*$ is accurate to $\mathcal{O}(\Delta x^{k-1})$.

In \Cref{stokestable:stokes:cylinder:error}, we show the drag coefficient and the percent error computed by the Immersed Boundary, IBSE-$1$, and IBSE-$2$ methods for both the spectral and finite difference discretizations over a range of values of $n_y$. The percent error is computed relative to the value of $132.36$ from \cite{fan1999galerkin}. For the IBSE-$1$ and IBSE-$2$ methods, we also show the improvement in $\%$ error as compared to the \IBa$\ $method. The IBSE method provides the greatest advantage when tight error tolerances are required. Using the finite difference discretization, to achieve a $1\%$ error, the \IBa$\ $method would need to use $n_y\approx512$, while the IBSE-$2$ method requires only $n_y\approx 64$; however to achieve an error of $0.01\%$, the \IBa$\ $method requires $n_y\approx 65536$, while IBSE requires only $n_y\approx 512$.


\begin{table}
	\centering
	\begin{tabular}{c|cccc|cccc}
		\hline
		\multirow{2}{*}{Method}		&	\multicolumn{4}{c}{Spectral}	&	\multicolumn{4}{c}{Finite Difference}						\\
							&	$n_y$											&	$C_D$		&	\% Error		&	Improvement		&	$n_y$											&	$C_D$		&	\% Error	&	Improvement	\\
		\hline
		\IBa			&	\multirow{3}{*}{$64$}		&	181.956	&	37.47				&	-						&	\multirow{3}{*}{$32$}		&	159.357	&	20.40		&	-		\\
		IBSE-$1$	&													&	129.564	&	2.11				&	18					&													&	124.161	&	6.19		&	3		\\
		IBSE-$2$	&													&	132.781	&	0.32				&	117					&													&	129.270	&	2.33		&	8		\\
		\hline
		\IBa			&	\multirow{3}{*}{$256$}	&	141.505	&	6.91				&	-						&	\multirow{3}{*}{$128$}	&	137.860	&	4.16		&	-		\\
		IBSE-$1$	&													&	132.437	&	0.058				&	119					&													&	131.228	&	0.86		&	5		\\
		IBSE-$2$	&													&	132.361	&	$\leq0.004$	&	$\geq1700$	&													&	132.202	&	0.12		&	35	\\
		\hline
		\IBa			&	\multirow{3}{*}{$1024$}	&	134.495	&	1.61				&	-						&	\multirow{3}{*}{$512$}	&	133.657	&	0.98		&	-		\\
		IBSE-$1$	&													&	132.425	&	0.049				&	33					&													&	132.142	&	0.16		&	6		\\
		IBSE-$2$	&													&	132.360	&	$\leq0.004$	&	$\geq400$		&													&	132.348	&	0.0091	&	108	\\
		\hline
	\end{tabular}
	\caption{The drag coefficient $C_D$ and $\%$ error, computed with respect to the reference value of $132.60$ from \cite{fan1999galerkin}. Results are provided for the \IBa, IBSE-$1$, and IBSE-$2$ methods, for both the spectral and finite difference discretizations, and for a coarse, medium, and fine discretization. For the IBSE-$1$ and IBSE-$2$ methods, the improvement in the percent error, as compared to the percent error produced by the \IBa$\ $method, is shown in the column labeled `Improvement'. When $n_y=64$ for the spectral discretization, $\Delta x=2\pi/64\approx 0.098$; when $n_y=32$ for the finite difference discretization $\Delta x=4/32=0.125$.}
	\label{stokestable:stokes:cylinder:error}
\end{table}

\subsection{A comparison of the numerical requirements}
\label{stokessection:stokes:inefficiences}

In \Cref{stokessection:stokes_test:flow_around_cylinder:spectral}, we discretized and solved a flow-around a cylinder problem in a confined channel using a Fourier spectral method to discretize the underlying differential operators. This provides a convenient and simple way to discretize the Stokes equations to high-order, but it comes with a significant drawback: the Dirichlet $\u=0$ boundary condition along the geometrically simple channel walls must be enforced by adding a force distribution to the non-physical region where $y<-2\zeta$ or $y > 2\zeta$. This causes two significant inefficiencies:
\begin{enumerate}
	\item Spatial degrees of freedom (for $\u$, $p$, $\bxi_\u$, $\xi_p$) must be solved for in the extension regions where $y<-2\zeta$ or $y > 2\zeta$.
	\item Additional boundary nodes must be placed along the interior boundaries at $y=2\zeta$ and $y=-2\zeta$ (corresponding to the channel walls).
\end{enumerate}
The additional boundary nodes along the channel walls can become a significant problem for fine discretizations. Boundary nodes are placed a distance $\Delta s$ apart, where $\Delta s$ is set to be approximately twice the grid spacing: $\Delta s\approx 2\Delta x$. For this geometry, when $n_y$ boundary points are used to discretize the spatial grid in the $y$ direction, $3n_y$ boundary nodes are required for each of the channel walls.

In \Cref{table:sc_size}, we summarize the number of boundary nodes required for the computation, the size of the Schur-Complement that must be inverted for IBSE-$1$ and IBSE-$2$, and the time required to form the Schur complement as well as to compute the solution to one Stokes problem once the Schur complement has been formed. For the steady Stokes problem, the time to the form the Schur complement dominates the computation, the cost of a single Stokes solve is much smaller. For time dependent problems, the Schur-complement need only be formed once, facilitating efficient time-stepping.

For the finest discretization that we run for the spectral discretization ($n_y=1024$), the IBSE-$2$ method requires a Schur-complement that is $53011\times 53011$; the finite-difference scheme at this discretization requires a Schur-complement that is only  $3856\times 3856$. This can be a substantial issue. The matrix for the Fourier discretization requires approximately 2.9 billion elements to be stored; using double precision, this requires over 22 gigabytes of RAM. Although this is feasible in two dimensions, the discretization of long channel walls would be computationally prohibitive in three dimensions. For the finite difference discretization of the same size, only about 120 megabytes of memory are required.

Despite the disparity in memory requirements, solving the Stokes equations using a Fourier spectral discretization is very efficient, and the cost of each full solve is between 10 and 20 times faster using the spectral discretization. The time required to form the Schur complement is approximately equivalent: the much larger Schur-complement size is offset by the speed of the solve itself. This disparity is largest for the Stokes problem. For Navier-Stokes problems with moderate Reynolds numbers, the finite-difference solves are approximately 10 times faster due to the fact that fewer iterations of the preconditioned GMRES algorithm are required to reach the same tolerance. In this case the spectral solve is only moderately faster per timestep, and the finite-difference discretization enjoys a significant advantage in the time required to form the Schur-complement.
\begin{table}
	\centering
	\begin{tabular}{ccccccccccc}
		\hline
		\multirow{2}{*}{$n_y$}	& \multirow{2}{*}{$n_{wall}$}	&	\multirow{2}{*}{$n_\text{cylinder}$}	&	\multirow{2}{*}{$n_\text{bdy}$}	&\multicolumn{2}{c}{SC size}	&	\multicolumn{2}{c}{Solve Time} &	\multicolumn{2}{c}{Formation Time}	\\
		&&&& Spectral & FD	& Spectral & FD	& Spectral & FD	\\
		\hline
		16		&	48		&	7		&	103		&	827		&	56		&	53	&	272	& 40000		&	15000		\\
		32		&	96		&	15	&	207		&	1659	&	120		&	58	&	568	&	100000	&	70000		\\
		64		&	192		&	30	&	414		&	3315	&	240		&	66	&	465	&	220000	&	110000	\\
		128		&	384		&	60	&	828		&	6627	&	480		&	56	&	630	&	370000	&	300000	\\
		256		&	768		&	120	&	1656	&	13251	&	960		&	43	&	640	&	570000	&	610000	\\
		512		&	1536	&	241	&	3313	&	26507	&	1928	&	38	&	722	&	1000000	&	1400000	\\
		1024	&	3072	&	482	&	6626	&	53011	&	3856	&	40	& 812	&	2100000	&	3100000	\\
					&				&			&				&				&				&	\multicolumn{4}{c}{\text{(In units of time to compute an FFT)}}	\\
		\hline
	\end{tabular}
	\caption[A comparison of the numerical efficiency of the Immersed Boundary Smooth Extension method, for a spectral discretization and a finite difference discretization.]{A summary of the number of boundary nodes required to discretize the flow around the cylinder problem studied in \Cref{stokessection:stokes_test:flow_around_cylinder}, along with the associated size of the Schur-complement matrix defined in \Cref{stokeseq:define_schur_complement} for both the Fourier spectral and finite-difference discretizations. The time required to solve one problem, as well as to form the entire Schur complement, normalized by the time required to execute one (real) FFT are also shown. The columns, from left to right, give $n_y$, the number of points used to discretize the computational domain $C$ in the spanwise direction; $n_\text{wall}$, the number of boundary nodes used to discretize each wall; $n_\text{cylinder}$, the number of boundary nodes used to discretize the cylinder; $n_\text{bdy}$, the number of boundary nodes used to discretize all of the boundaries (for the spectral discretization), the size of the Schur complement (for both the spectral and finite-difference discretizations), and the (estimated) time to form the Schur-complement operator, again normalized by the time to execute one real FFT.}
	\label{table:sc_size}
\end{table}


\section{Navier-Stokes equations}
\label{stokessection:navier_stokes}

Because the IBSE method enables the efficient solution of the Stokes problem in the domain $\Omega$, unsplit time-stepping for the Navier-Stokes problem is simple when the non-linear terms are treated explicitly in time. We integrate the Navier-Stokes equation in time using a fourth-order implicit-explicit (IMEX) Backward Differentiation formula \cite{Hundsdorfer2007}:
\begin{multline}
	\frac{25}{12}\u^{n+1} - 4\u^n + 3\u^{n-1} - \frac{4}{3}\u^{n-2} + \frac{1}{4}\u^{n-3} = \\
		\Delta t \left[\mathcal{I}(\u^{n+1}) + 4\mathcal{E}(\u^n) - 6\mathcal{E}(\u^{n-1}) + 4\mathcal{E}(\u^{n-2}) - \mathcal{E}(\u^{n-3})\right],
\end{multline}
where $\mathcal{I}$ evaluates the \emph{stiff} terms in $\u_t$ while $\mathcal{E}$ evaluates the \emph{non-stiff} terms in $\u_t$.  For this problem, $\mathcal{I}(\u)$ is the Stokes operator and $\mathcal{E}(\u)=-\u\cdot\nabla \u$, discretized using second-order centered differences. The use of this type of timestepping eliminates errors due to a projection step and removes the need for artificial pressure boundary condtions, and is appropriate for low-Reynolds number problems and for problems where accuracy near to the boundary is critical.

\begin{figure}
	\centering
	\includegraphics[width=0.9\textwidth]{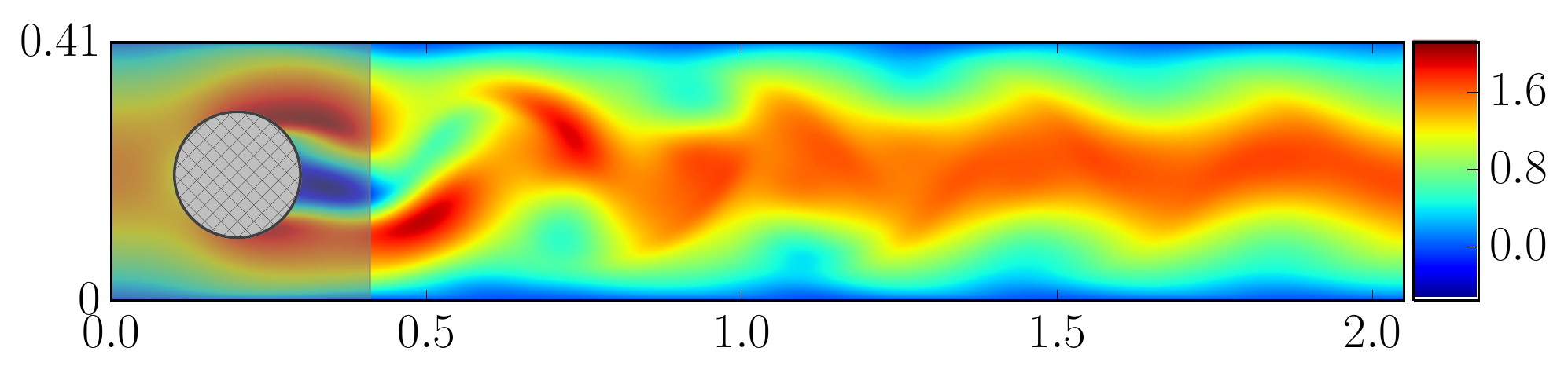}
	\caption[The streamwise velocity and domain (highlighting the localized extension domain) for the confined flow around a cylinder test problem used to validate the accuracy of the Immersed Boundary Smooth Extension, coupled to a finite difference solver, for a Navier-Stokes flow (for the un-steady case with Reynolds number $100$).]{The solution $u$ to the problem solved in \Cref{stokessection:navier_stokes:unsteady}, at $t=12.0$, produced using a second-order finite difference discretization to the Navier-Stokes equations as the underlying PDE solver for the IBSE-$2$ method with $n_y=512$. The extension region $E$ is shown in white with hatching. The extension functions $\bxi_\u$ and $\xi_p$ are not solved for on the entirety of the computational domain, but rather in the much smaller region $C_E=[0.41,0.41]\times[0.41,0.41]$, localized around $E$, that is shaded in light gray.}
	\label{fd_navier_stokes_unsteady_solution}
\end{figure}

To test the performance of the IBSE method on a time dependent problem, we set up and solve a flow around a cylinder problem in a standard benchmark flow geometry \cite{Schafer1996}, shown in \Cref{fd_navier_stokes_unsteady_solution}. A variety of simple scalar benchmarks for this problem are available and are reported as ranges which ensure a percent error of at most $0.15\%$ to $3\%$, depending on the benchmark \cite{Schafer1996}. To demonstrate the validity of the timestepping scheme, we first demonstrate the solver for a problem with an unsteady solution, at Reynolds number $100$, and show quantitative agreement with time-dependent benchmarks. Unfortunately, very high-accuracy benchmarks are unavailable for this problem. To provide a better measure of the accuracy of the IBSE method, we solve a steady problem, with Reynolds number $20$. For some subset of the benchmarks that we examine for this problem, results accurate to nearly machine precision are available from high-order spectral computations \cite{nabh1998high}.

For both the unsteady and steady cases, the test problem is the confined flow around a cylindrical obstacle. The computational domain $C$ is $[0,2.05]\times[0,0.41]$, the extension region $E=B_{0.05}(0.2,0.2)$ is the cylinder of radius $0.05$ centered at $(0.2,0.2)$, and the physical domain $\Omega=C\setminus E$. We will denote the internal boundary $\partial B_{0.05}(0.2,0.2)$ by $\Gamma$, while the boundary of the computational domain $C$ will be denoted by $\partial C$. $\partial C$ is divided into four parts: $\partial C=\cup\{\partial C_l, \partial C_r, \partial C_b, \partial C_t\}$, denoting the boundaries at $x=0$, $x=2.05$, $y=0$, and $y=0.41$, respectively. No-slip boundary conditions are applied on $\partial C_b$ and $\partial C_t$, homogeneous Neumann boundary conditions are applied on $\partial C_r$, and the inflow condition $\u=(4sU_my(0.41-y)/0.41^2,0)$ is applied on $\partial C_l$, where $U_m$ is the maximum value of the $u$ velocity at inflow and will be chosen to set the Reynolds number. To be precise, we solve
\begin{subequations}
	\label{stokeseq:navier_stokes}
	\begin{align}
		\u_t + \u\cdot\grad\u - \nu\Delta\u + \grad p	&=	0	&&\text{in }\Omega,	\\
		\grad\cdot\u	&=	0	&&\text{in }\Omega,	\\
		\u				&=	0	&&\text{on }\Gamma,\ \partial C_b,\text{ and }\partial C_t,	\\
		u(y)			&=	4U_my(0.41-y)/0.41^2	&&\text{on }\partial C_l,	\\
		v(y)			&=	0	&&\text{on }\partial C_l,	\\
		\partial\u/\partial n	&=	0	&&\text{on }\partial C_r.
	\end{align}
\end{subequations}

For the IBSE method, $C_E$, the computational domain for computing the extension functions $\bxi_\u$ and $\xi_p$ is chosen to be $[0,0.41]\times[0,0.41]$. This geometry, along with the $u$ value of the solution for the unsteady Reynolds number $100$ case, produced by the IBSE-$2$ method with $n_y=512$, is shown in \Cref{fd_navier_stokes_unsteady_solution}. We note that our geometry is slightly different than the benchmark geometry used in \cite{Schafer1996}. Their domain is taken to be a rectangle with width $2.2$ and height $0.41$, our rectangle is slightly less wide so that the rectangle has an aspect ratio of $5$, allowing the grid spacing in the $x$ direction and $y$ direction to be equal.

In practice, timestepping involves the explicit computation of a forcing function:
\begin{equation}
	\frac{25}{12}\f_\u = 4\u^n - 3\u^{n-1} + \frac{4}{3}\u^{n-2} - \frac{1}{4}\u^{n-3} + \Delta t\left[4\mathcal{E}(\u^n)-6\mathcal{E}(\u^{n-1})+4\mathcal{E}(\u^{n-2})-\mathcal{E}(\u^{n-3})\right],
\end{equation}
followed by solution of the equation
\begin{subequations}
	\label{stokeseq:navier_stokes_update}
	\begin{align}
		\left(\mathbb{I}-\frac{12}{25}\nu\Delta t\Delta\right) \u^{n+1}  + \grad p^{n+1}	&= \f_\u	& &\text{in }\Omega,	\\
		\grad\cdot\u	&=	0	&	&\text{in }\Omega,\\
		\u^{n+1}	&=	0	&	&\text{on }\Gamma,\ \partial C_b,\text{ and }\partial C_t,	\\
		u^{n+1}(y)			&=	4U_my(0.41-y)/0.41^2	&&\text{on }\partial C_l,	\label{stokeseq:navier_stokes_update:inflow}	\\
		v^{n+1}(y)			&=	0	&&\text{on }\partial C_l,	\\
		\partial\u^{n+1}/\partial n	&=	0	&&\text{on }\partial C_r.
	\end{align}
\end{subequations}

\subsection{Results: Unsteady flow (Reynolds number $100$)}
\label{stokessection:navier_stokes:unsteady}

To achieve a Reynolds number of $100$ for the unsteady case, $U_m$ is set to $1.5$ in \Cref{stokeseq:navier_stokes_update:inflow}. For this problem the Reynolds number is defined as $Re=\overline{U}D/\nu$, where $\overline U$ gives the average velocity across the inflow boundary. For the prescribed boundary condition, $\overline U = \frac{2}{3}U_m$. We integrate \Cref{stokeseq:navier_stokes} to $t=12.0$ with a timestep of $\Delta t=\Delta x/5$, for $n=2^7$ and $n=2^{9}$, using the IBSE-$2$ method when solving \Cref{stokeseq:navier_stokes_update}. The timestep $\Delta t$ need not be taken so small at this Reynolds number, but is chosen conservatively to ensure that error is spatially dominated. The $u$ velocity at $t=12.0$ is shown in \Cref{fd_navier_stokes_unsteady_solution}.

The benchmarks measured in \cite{Schafer1996} for this problem include the maximum drag coefficient ($C_D$), the maximum lift coefficient ($C_L$), the Strouhal number $St$, and the pressure difference ($\Delta P$), measured at $t=t_0+1/(2f)$, where $f$ is the frequency of separation and $t_0$ is the time at which the lift coefficient is maximized. The benchmark ranges, along with the values computed using the IBSE-$2$ method with $n_y=128$ and $n_y=512$ are shown in \Cref{stokestable:unsteady_results}. For the IBSE-$2$ method with $n=512$, the computed values of the benchmarks are consistent with the benchmark range, although the maximum value for the lift coefficient is slightly low. We note that method 9a from \cite{Schafer1996}, whose values agree very well with our results and the high-order results from \cite{nabh1998high} for the steady $\text{Re}=20$ case (see \Cref{stokessection:navier_stokes:steady}) also reports a lower maximum value of $C_L=0.9862$ that is in strong agreement with the IBSE-$2$, $n_y=512$ result for $C_L=0.9866$.

\begin{table}
	\centering
	\begin{tabular}{ccccc}
		\hline
										&	$C_D$						&	$C_L$	&	Strouhal Number	&	$\Delta P$	\\
		\hline
		Reference Range	&	$[3.22, 3.24]$	&	$[0.99, 1.01]$	&	$[0.295, 0.305]$	&	$[2.46, 2.50]$	\\
		IBSE-$2$, $n_y=128$	&	3.3123	&	0.9802	&	0.3005	&	2.4770	\\
		IBSE-$2$, $n_y=512$	&	3.2338	&	0.9866	&	0.3018	&	2.4979	\\
		\hline
	\end{tabular}
	\caption{The drag coefficient ($C_D$), lift coefficient ($C_L$), Strouhal number, and pressure difference computed by the IBSE-$2$ method for the unsteady Navier Stokes problem with Reynolds number $100$. The results are shown for $n_y=128$ and $n_y=512$, and compared to the reference ranges provided in \cite{Schafer1996}.}
	\label{stokestable:unsteady_results}
\end{table}

\subsection{Results: Steady flow (Reynolds number $20$)}
\label{stokessection:navier_stokes:steady}

\begin{figure}[b!]
	\centering
	\includegraphics[width=0.9\textwidth]{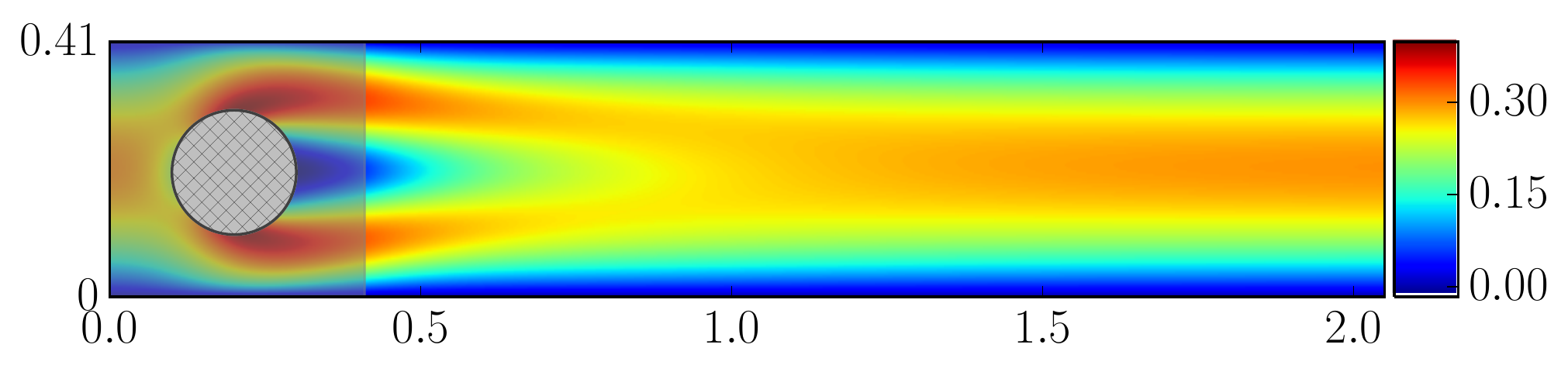}
	\caption[The streamwise velocity and domain (highlighting the localized extension domain) for the confined flow around a cylinder test problem used to validate the accuracy of the Immersed Boundary Smooth Extension, coupled to a finite difference solver, for a Navier-Stokes flow (for the steady case with Reynolds number $20$).]{The solution $u$ to the problem solved in \Cref{stokessection:navier_stokes:steady}, at steady state, produced using a second-order MAC finite difference discretization to the Navier-Stokes equations as the underlying PDE solver for the IBSE-$2$ method with $n_y=512$. The extension region $E$ is denoted with crosshatches. The extension functions $\bxi_\u$ and $\xi_p$ are not solved for on the entirety of the computational domain, but rather in the much smaller region $C_E=[0.41,0.41]\times[0.41,0.41]$, localized around $E$, that is shaded in light gray.}
	\label{fd_navier_stokes_steady_solution}
\end{figure}

\begin{figure}
	\centering
	\hspace*{\fill}
	\begin{subfigure}[b]{0.45\textwidth}
		\centering
	\includegraphics[width=\textwidth]{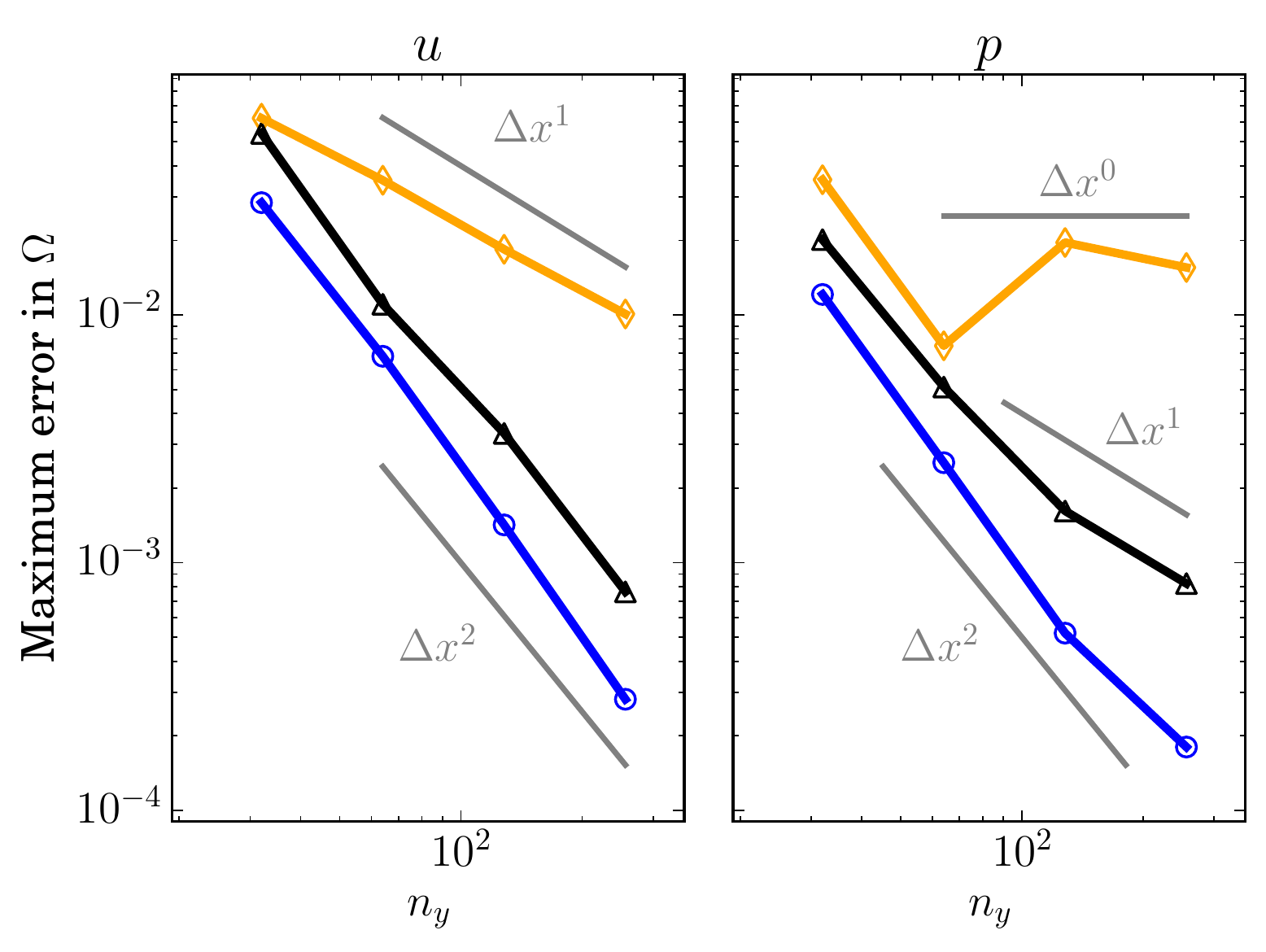}
		\subcaption{Refinement study for $u$, $p$ in $L^\infty(\Omega)$}
		\label{stokesfig:navier_stokes_steady:refine}
	\end{subfigure}
	\hfill
	\begin{subfigure}[b]{0.252\textwidth}
		\centering
		\includegraphics[width=\textwidth]{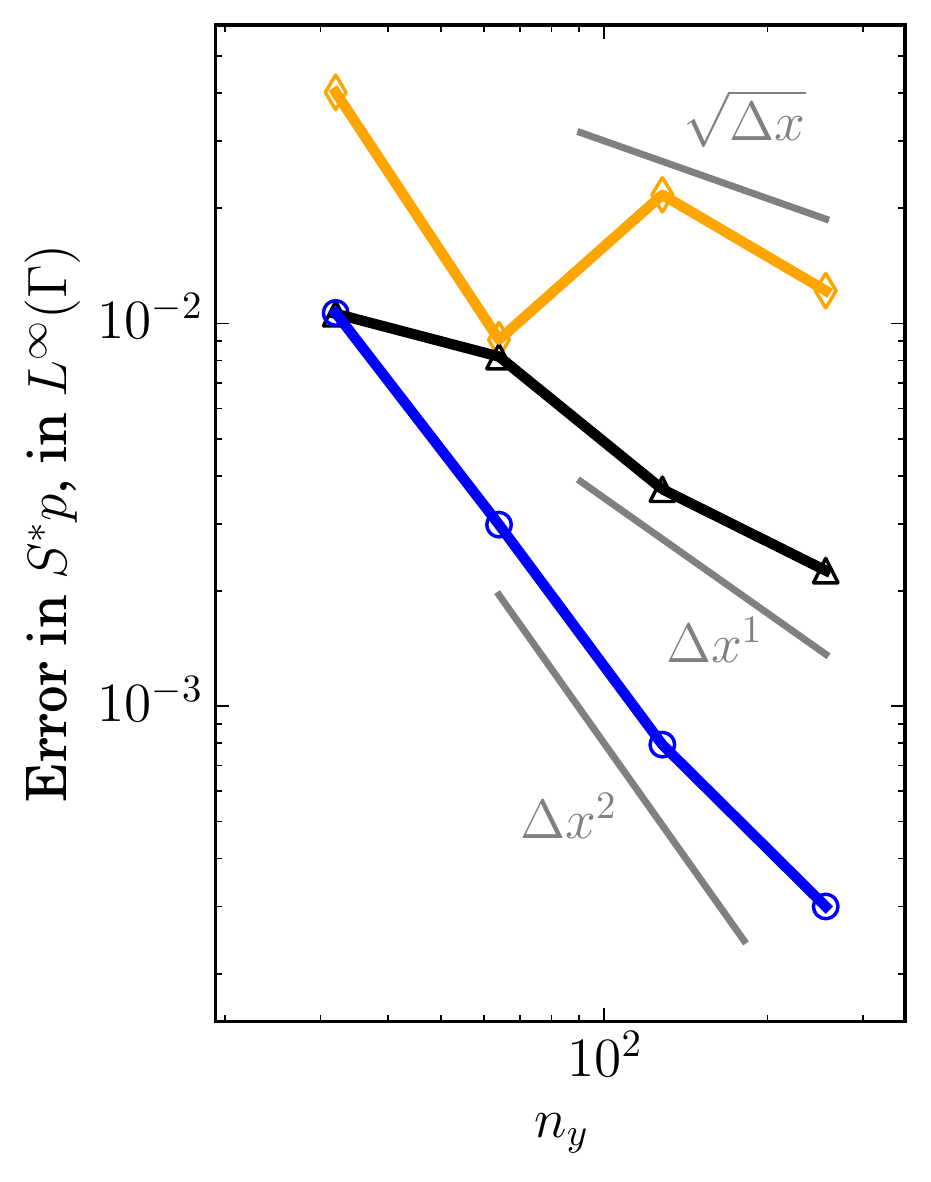}
		\subcaption{$L^\infty(\Gamma)$ difference in $S^*p$}
		\label{stokesfig:navier_stokes_steady:max_bp_error}
	\end{subfigure}
	\hfill
	\begin{subfigure}[b]{0.26\textwidth}
		\centering
		\includegraphics[width=\textwidth]{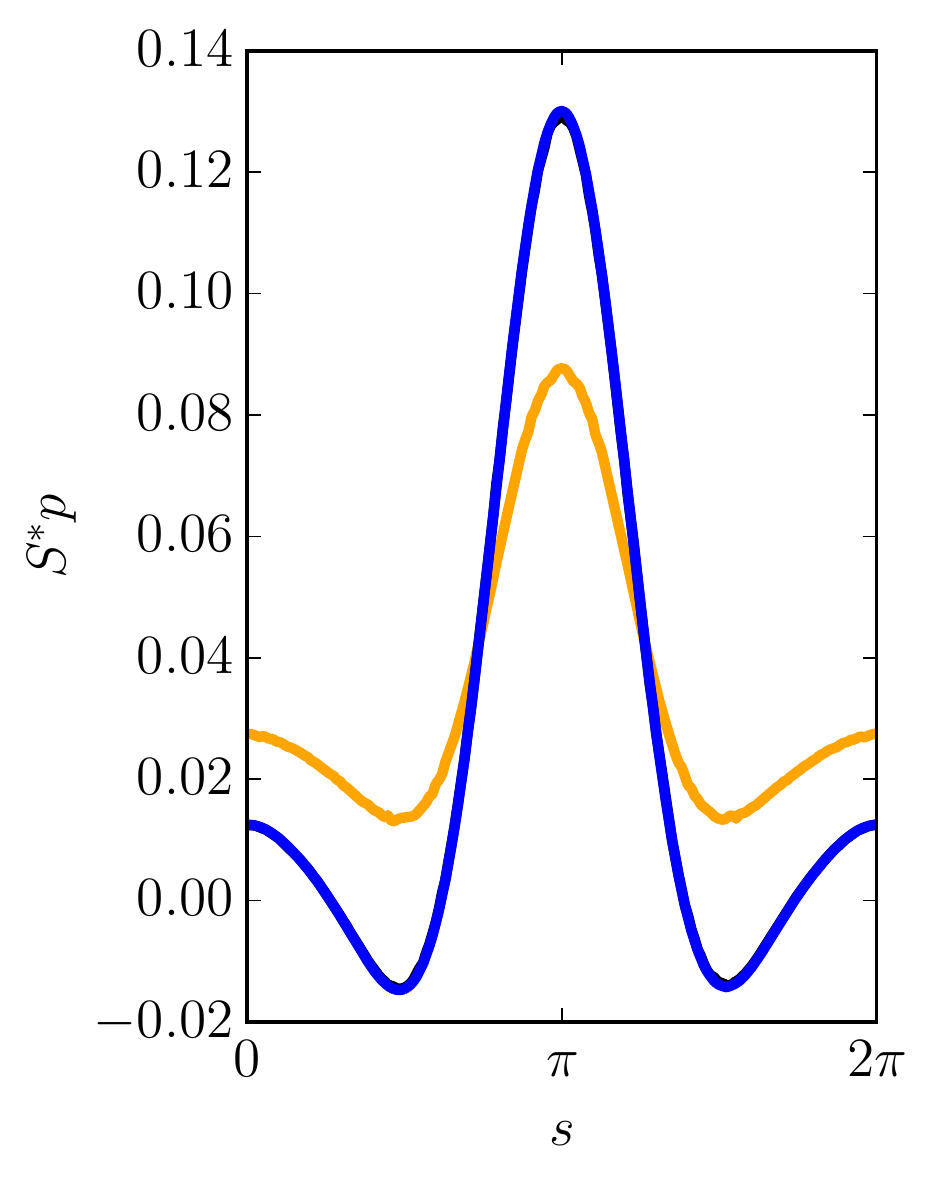}
		\subcaption{$S^*p(s)$}
		\label{stokesfig:navier_stokes_steady:bp}
	\end{subfigure}
	\hspace*{\fill}
	\caption[Refinement study showing errors for the Reynolds number $20$ Navier-Stokes test problem at steady state, for solutions produced by the Immersed Boundary Smooth Extension method.]{Refinement study showing the convergence of differences between successively refined solutions for $u$ and $p$ in $L^\infty(\Omega)$ (\Cref{stokesfig:navier_stokes_steady:refine}), and for $S^*p$ in $L^\infty(\Gamma)$ (\Cref{stokesfig:navier_stokes_steady:max_bp_error}), to the steady Navier-Stokes problem defined in \Cref{stokeseq:navier_stokes} at steady state ($t=12.0)$, produced using the \IBa ({\color{py_orange} $\lozenge$}), IBSE-$1$ ({{\color{black} $\triangle$}), and IBSE-$2$ ({\color{py_blue} \Circle}}) methods. \Cref{stokesfig:navier_stokes_steady:bp} shows the pressure function interpolated to the boundary, although this converges to something, as shown in \Cref{stokesfig:navier_stokes_steady:max_bp_error}, it is not converging to the correct value. Note that the lines for IBSE-$1$ and IBSE-$2$ visually coincide.}
	\label{stokesfig:navier_stokes_steady:max_error}
\end{figure}

For a steady case, $U_m$ is set to $0.3$ \cite{Schafer1996}, and $\nu$ is set to $0.001$, yielding a Reynolds number $Re=\overline{U}D/\nu=20$. We integrate \Cref{stokeseq:navier_stokes} to $t=12.0$ with a timestep of $\Delta t=\Delta x/5$, for $n=2^5$ to $n=2^{9}$, using the IB method as well as IBSE-$1$ and IBSE-$2$ when solving \Cref{stokeseq:navier_stokes_update}. As with the unsteady case, the timestep need not be taken so small but is chosen conservatively to ensure that error is spatially dominated. Convergence is assessed by comparing the ratio of the $L^\infty(\Omega)$ difference between solutions at successive levels of refinement (bilinear interpolation is used to transfer data between grids to do the refinement study).  Results are shown in \Cref{stokesfig:navier_stokes_steady:max_error}. As expected, the IB method achieves first order convergence in $L^\infty(\Omega)$ for the velocity, but fails to converge pointwise for the pressure. Both the IBSE-$1$ and IBSE-$2$ methods converge at second-order for the velocity field. Although the IBSE-$2$ method produces $C^2(C)$ solutions, third order convergence is not possible, as the underlying finite-difference discretization is only accurate to second order. Despite asymptotically similar convergence, the actual error is lower for all discretizations.  For the pressure, IBSE-$1$ converges at first order and IBSE-$2$ converges at second order.

To emphasize the improved convergence of the IBSE method for elements of the stress tensor, we examine the convergence of the pressure field $p$, interpolated to the boundary, in \Cref{stokesfig:navier_stokes_steady:max_bp_error}. Convergence is assessed by comparing the ratio of the $L^\infty(\Gamma)$ difference between $S^* p$ at different resolutions. From examining \Cref{stokesfig:navier_stokes_steady:max_bp_error}, it appears that the IB method converges, albeit slowly. However, looking at \Cref{stokesfig:navier_stokes_steady:bp}, we see that although the IB method is converging to \emph{something} on the boundary, that value is not consistent with the value determined by the IBSE method. For the IBSE-$1$ and IBSE-$2$ methods, $S^*p$ converges at first- and second-order, respectively. We note that the benchmark pressure difference, $\Delta P$, for the IBSE-$1$ and IBSE-$2$ methods is consistent with the ranges reported in \cite{Schafer1996} (see \Cref{table:benchmark_comparison_steady}), giving us good reason to believe that the pressure values (up to a constant) on the cylinder shown in \Cref{stokesfig:navier_stokes_steady:bp} are accurate.

For this problem the availability of high-accuracy values \cite{nabh1998high} for the drag coefficient ($C_D$), the lift coefficient ($C_L$) and the pressure difference ($\Delta P$) allows for precise benchmarking. We will also compute the length of the recirculation zone ($L_a$), for which only a benchmark range \cite{Schafer1996} is available. For the IBSE method, the drag and lift coefficients are computed by directly computing the stress, interpolating to the boundary, and computing the appropriate integral over $\Gamma$. For the IB method, this provides an inconsistent estimate, the appropriate integral over the singular forces is computed instead \cite{Lai2000}. No pressure difference ($\Delta P$) is reported for the IB method, as a consistent estimate of the pressure on the boundary cannot be obtained. This comparison is shown in \Cref{table:benchmark_comparison_steady}. The IB method does surprisingly well for all benchmarks: despite the fact that components of the stress tensor cannot be directly evaluated on the boundary, the singular forces $\G$ provide an alternative measure that is convergent. For $C_D$, $C_L$, and $\Delta P$, we expect the IBSE-$1$ method to provide first-order estimates due to the fact that the stress $\sigma$ is only continuous across the boundary $\Gamma$. The estimates for the benchmarks have approximately the same accuracy as the IB estimates. For the IBSE-$2$ method, the stress is differentiable across the boundary, and we expect these benchmarks to be second-order accurate. Unsurprisingly this method provides much more accurate values. For $n=512$, the IBSE-$2$ method yields a $\%$ error in the drag coefficient that is nearly two orders of magnitude smaller than the $\%$ error provided by the IB method. We again emphasize that with the IB method, consistent estimates of $p$ cannot be obtained at the boundary, and so we do not attempt to calculate a pressure difference $\Delta P$ for these methods.

\begin{table}[htb!]
	\centering
	\begin{tabular}{cc|cc|cc|c|cc}
		\hline
		Benchmark				&			&	\multicolumn{2}{c}{$C_D$} 	&	\multicolumn{2}{c}{$C_L$}	&	$L_a$	&	\multicolumn{2}{c}{$\Delta P$}		\\
		Range			&			&	\multicolumn{2}{c}{$[5.57, 5.59]$}	&	\multicolumn{2}{c}{$[0.0104, 0.0110]$}	&	$[0.0842,0.0852]$	&	\multicolumn{2}{c}{$[0.1172, 0.1176]$}	\\
		Value			&			&	\multicolumn{2}{c}{$5.57953523384$}	&	\multicolumn{2}{c}{$0.010618948146$}	&	-	&	\multicolumn{2}{c}{$0.11752016697$}	\\
								&			&	value	&	$\%_\text{err}$	&	value	&	$\%_\text{err}$	&	value	&	value	&	$\%_\text{err}$	\\
		\hline
		\multirow{3}{*}{$n=128$}&	\IBa		&	5.76813	&	3.38003			&	0.01306	&	22.9463			&	0.08803	&	-		&	-		\\
								&	IBSE-$1$&	5.52533	&	0.97142			&	0.00405	&	61.8221			&	0.08514	&	0.11625	&	1.08426	\\
								&	IBSE-$2$&	5.58697	&	0.13323			&	0.00875	&	17.5530			&	0.08481	&	0.11714	&	0.32420	\\
		\hline		
		\multirow{3}{*}{$n=512$}	&	\IBa	&	5.62468	&	0.80919			&	0.01081	&	1.81973			&	0.08567	&	-		&	-		\\
								&	IBSE-$1$&	5.55464	&	0.44618			&	0.01228	&	15.6430			&	0.08465	&	0.11677	&	0.63829	\\
								&	IBSE-$2$&	5.57883	&	0.01271			&	0.01069	&	0.71540			&	0.08464	&	0.11759	&	0.06075	\\
		\hline
	\end{tabular}
	\caption[A comparison of benchmark values with values produced by the Immersed Boundary Smooth Extension method for the Reynolds number $20$ confined channel flow around a cylinder Navier-Stokes problem.]{The drag coefficient $(C_D)$, lift coefficient $(C_L)$, length of the recirculation zone $(L_a)$, and the pressure difference $(\Delta P)$, computed by the \IBa, IBSE-$1$, and IBSE-$2$ methods, compared to the benchmark ranges computed in \cite{Schafer1996} (labeled as ``Benchmark Range''), as well as the results from the high-order spectral computations in \cite{nabh1998high} (labeled as ``Benchmark Value''). The percent errors reported are computed in reference to the benchmark value.}
	\label{table:benchmark_comparison_steady}
\end{table}


\section{Discussion}

In \cite{Stein2015}, we introduced the IBSE method for solving elliptic and parabolic PDE on general smooth domains to an arbitrarily high order of accuracy using simple Fourier spectral methods. In this paper, we extend the methodology to fluid problems with divergence constraints, including the incompressible Stokes and Navier Stokes equations. We decouple the extension procedure from the underlying PDE solver, and demonstrate accurate solutions with an underlying Fourier spectral solver as well as a second-order finite difference method. For the velocity field, we demonstrate up to third-order pointwise convergence, in contrast to the first-order pointwise convergence of the IB method. In addition, we show up to second-order pointwise convergence of the pressure function and all elements of the fluid stress tensor. In the IB method, these quantities do not converge pointwise.

In the IBSE method, accuracy is obtained by forcing the solution to be globally smooth on the entire computational domain.  This is accomplished by solving for the smooth extension of the \emph{unknown} solution, allowing for application of the IBSE method to the steady Stokes equations.  Remarkably, the fully coupled problem for the solution to the incompressible Stokes equations (combined with the equations that define the smooth extensions to the velocity and pressure fields) can be effectively reduced to a small, dense system of equations, and the cost to invert this system can be minimized by precomputing and storing its LU-factorization. The IBSE method thus provides an efficient algorithm for implicit or implicit-explicit (IMEX) timestepping on stationary domains.  The method requires minimal geometric information regarding the boundary: only its position, normals, and an indicator variable denoting whether points in the computational grid are inside of the physical domain $\Omega$ or not, enabling simple and robust code.  The method is flexible enough to allow high-order discretization in space and time for a wide range of nonlinear PDE using straightforward implicit-explicit timestepping schemes.

The IBSE method decouples the discretization of the extension problem from the discretization of the underlying PDE. By itself, this allows some gains in efficiency when the extension region is small compared to the physical domain $\Omega$, as is often the case when studying the flow around objects. More importantly, however, the methodology can be easily coupled to existing codes. This enables the simple and accurate inclusion of embedded boundaries into ongoing research without significant redevelopment of code.  Additionally, it allows for the ability to build on top of validated and optimized code for the solution of PDE, including those that that use adaptive mesh refinement, without sacrificing accuracy near to the embedded boundary.

The IBSE method shares some similarities with Fourier Continuation (FC) \cite{Lyon2010a,Lyon2010,Albin2011} and Active Penalty (AP) \cite{Shirokoff2013} methods.  The FC method uses Fourier methods to obtain a high-order discretization, and the idea that smooth, non-periodic functions can be turned into smooth and periodic functions by extending them into a larger domain in some appropriate way.  There are two key differences between the FC method and the IBSE method.  The FC method ({\small\emph{i}}) uses dimensional splitting to reduce the problem to a set of one-dimensional problems and ({\small\emph{ii}}) relies on data from the previous timestep in order to generate the Fourier continuations that allow for their high-order spatial accuracy.  This forces the FC method to use an Alternating-Direction Implicit scheme to take large stable timesteps when solving parabolic equations, complicating the implementation of high-order timestepping.  In addition, the FC method requires the use of an iterative solver for computing solutions to the Poisson equation, complicating an efficient discretization of the \emph{incompressible} Navier-Stokes equations, although the FC method has been used to solve the compressible Navier-Stokes equations to high order \cite{Albin2011}.  In contrast to the conditioning issues faced by the IBSE method, the one-dimensional nature of the Fourier-continuation problem allows the FC method to use high-precision arithmetic in certain precomputations to effectively eliminate the conditioning problem inherent in constructing smooth extensions \cite{Platte2011}.  This enables the FC method to obtain stable methods that converge to higher order than can be achieved by the IBSE method in double-precision arithmetic.

The AP method, like the FC method, relies on data from previous timesteps in the way that it imposes smoothness on the solution of the PDE.  A large drag force is applied that penalizes deviations of the solution in the extension region from the smooth extension of the solution at the previous timestep.  This dependence on data from previous timesteps forces the AP method to use \emph{explicit} timestepping when solving parabolic equations. The AP algorithm for the Navier-Stokes equations uses a projection method, requires explicit time-stepping, and is only able to achieve second- and first-order convergence for the velocity and pressure fields, respectively.

In contrast, the IBSE method directly solves the steady Stokes problem, smoothly extending both the unknown velocity and pressure fields. The extension equations are coupled directly to the unknown solutions through the unsplit Stokes problem, eliminating the errors inherent in projection methods at low Reynolds number. The price that we pay for extending the unknown solution is that the high-efficiency that we obtain is only currently achievable on stationary domains.  This is because the inversion method used for the IBSE method relies on an expensive precomputation that depends on the physical domain and the discretization.  However, there is no fundamental obstacle to the application of the IBSE method to problems with moving domains.  Instead, the challenge is to find a robust and efficient method to invert the IBSE system in \Cref{stokeseq:the_stokes_system} that does not require substantial precomputation.  Recent progress has been made for preconditioning similar systems of equations for the simulation of rigid-body motion in an Immersed Boundary framework \cite{kallemov2016immersed}.  The integration of these ideas into the IBSE method to allow for simulation of moving boundary problems will be an area of active future research.

We have only implemented two-dimensional examples in this paper.  The IBSE method extends to three dimensions without any changes, although there is one minor difficulty that must be resolved: the production of an accurate enough quadrature for the discretization of the boundary integrals in the spread operators $S$ and $T_k$.  For two-dimensional problems, this comes nearly for free since the simple quadrature rule given in \Cref{stokessection:numerics:operators} is spectrally accurate for closed one-dimensional curves.  High-order quadrature rules for two-dimensional surfaces are more complicated but well-developed, and high-order surface representation has been incorporated into the Immersed Boundary method \cite{griffith2012hybrid,Shankar2012}.



\section*{Acknowledgements}

This work was supported in part by the National Science Foundation under Grant DMS-1160438.



\clearpage

\bibliography{bib,manual_bib}{}

\begin{thebibliography}{49}
\providecommand{\natexlab}[1]{#1}
\providecommand{\url}[1]{\texttt{#1}}
\expandafter\ifx\csname urlstyle\endcsname\relax
  \providecommand{\doi}[1]{doi: #1}\else
  \providecommand{\doi}{doi: \begingroup \urlstyle{rm}\Url}\fi

\bibitem[Peskin(2002)]{Peskin2002}
Charles~S. Peskin.
\newblock {The immersed boundary method}.
\newblock \emph{Acta Numerica}, 11:\penalty0 479--517, 2002.
\newblock ISSN 0962-4929.
\newblock \doi{10.1017/S0962492902000077}.

\bibitem[Mittal and Iaccarino(2005)]{Mittal2005}
Rajat Mittal and Gianluca Iaccarino.
\newblock {Immersed Boundary Methods}.
\newblock \emph{Annual Review of Fluid Mechanics}, 37\penalty0 (1):\penalty0
  239--261, 2005.
\newblock ISSN 0066-4189.
\newblock \doi{10.1146/annurev.fluid.37.061903.175743}.

\bibitem[Olson and Layton(2014)]{olson2014simulating}
Sarah~D Olson and Anita~T Layton.
\newblock Simulating biofluid-structure interactions with an immersed boundary
  framework--a review.
\newblock \emph{Biological Fluid Dynamics: Modeling, Computations, and
  Applications}, 628:\penalty0 1, 2014.

\bibitem[Kallemov et~al.(2016)Kallemov, Bhalla, Griffith, and
  Donev]{kallemov2016immersed}
Bakytzhan Kallemov, Amneet Bhalla, Boyce Griffith, and Aleksandar Donev.
\newblock An immersed boundary method for rigid bodies.
\newblock \emph{Communications in Applied Mathematics and Computational
  Science}, 11\penalty0 (1):\penalty0 79--141, 2016.

\bibitem[Taira and Colonius(2007)]{Taira2007}
Kunihiko Taira and Tim Colonius.
\newblock {The immersed boundary method: A projection approach}.
\newblock \emph{Journal of Computational Physics}, 225\penalty0 (2):\penalty0
  2118--2137, 2007.
\newblock ISSN 00219991.
\newblock \doi{10.1016/j.jcp.2007.03.005}.

\bibitem[Teran and Peskin(2009)]{Teran2009}
Joseph~M Teran and Charles~S Peskin.
\newblock {Tether force constraints in Stokes flow by the immersed boundary
  method on a periodic domain}.
\newblock \emph{SIAM Journal on Scientific Computing}, 31\penalty0
  (5):\penalty0 3404--3416, 2009.

\bibitem[Li and Ito(2006)]{li2006immersed}
Zhilin Li and Kazufumi Ito.
\newblock \emph{The immersed interface method: numerical solutions of PDEs
  involving interfaces and irregular domains}, volume~33.
\newblock Siam, 2006.

\bibitem[Fedkiw et~al.(1999)Fedkiw, Aslam, Merriman, and Osher]{Fedkiw1999}
Ronald~P Fedkiw, Tariq Aslam, Barry Merriman, and Stanley Osher.
\newblock {A Non-oscillatory Eulerian Approach to Interfaces in Multimaterial
  Flows (the Ghost Fluid Method)}.
\newblock \emph{Journal of Computational Physics}, 152\penalty0 (2):\penalty0
  457--492, 1999.
\newblock ISSN 00219991.
\newblock \doi{10.1006/jcph.1999.6236}.
\newblock URL
  \url{http://www.sciencedirect.com/science/article/pii/S0021999199962368}.

\bibitem[Angot et~al.(1999)Angot, Bruneau, and Fabrie]{Angot1999}
Philippe Angot, Charles-Henri Bruneau, and Pierre Fabrie.
\newblock {A penalization method to take into account obstacles in
  incompressible viscous flows}.
\newblock \emph{Numerische Mathematik}, 81\penalty0 (4):\penalty0 497--520,
  1999.
\newblock ISSN 0029-599X.
\newblock \doi{10.1007/s002110050401}.

\bibitem[Lai and Peskin(2000)]{Lai2000}
Ming-Chih Lai and Charles~S. Peskin.
\newblock {An Immersed Boundary Method with Formal Second-Order Accuracy and
  Reduced Numerical Viscosity}.
\newblock \emph{Journal of Computational Physics}, 160\penalty0 (2):\penalty0
  705--719, 2000.
\newblock ISSN 00219991.
\newblock \doi{10.1006/jcph.2000.6483}.
\newblock URL
  \url{http://linkinghub.elsevier.com/retrieve/pii/S0021999100964830}.

\bibitem[Mark and van Wachem(2008)]{Mark2008}
Andreas Mark and Berend G~M van Wachem.
\newblock {Derivation and validation of a novel implicit second-order accurate
  immersed boundary method}.
\newblock \emph{Journal of Computational Physics}, 227\penalty0 (13):\penalty0
  6660--6680, 2008.
\newblock ISSN 00219991.
\newblock \doi{10.1016/j.jcp.2008.03.031}.

\bibitem[Linnick and Fasel(2005)]{Linnick2005}
Mark~N. Linnick and Hermann~F. Fasel.
\newblock {A high-order immersed interface method for simulating unsteady
  incompressible flows on irregular domains}.
\newblock \emph{Journal of Computational Physics}, 204\penalty0 (1):\penalty0
  157--192, 2005.
\newblock ISSN 00219991.
\newblock \doi{10.1016/j.jcp.2004.09.017}.

\bibitem[Liu and Zheng(2014)]{Liu2014}
Jian~Kang Liu and Zhou~Shun Zheng.
\newblock {Efficient high-order immersed interface methods for heat equations
  with interfaces}.
\newblock \emph{Applied Mathematics and Mechanics}, 35\penalty0
  (51174236):\penalty0 1189--1202, 2014.
\newblock ISSN 02534827.
\newblock \doi{10.1007/s10483-014-1851-6}.

\bibitem[Xu and Wang(2006)]{Xu2006}
Sheng Xu and Z.~Jane Wang.
\newblock {An immersed interface method for simulating the interaction of a
  fluid with moving boundaries}.
\newblock \emph{Journal of Computational Physics}, 216\penalty0 (2):\penalty0
  454--493, 2006.
\newblock ISSN 00219991.
\newblock \doi{10.1016/j.jcp.2005.12.016}.

\bibitem[Zhong(2007)]{Zhong2007}
Xiaolin Zhong.
\newblock {A new high-order immersed interface method for solving elliptic
  equations with imbedded interface of discontinuity}.
\newblock \emph{Journal of Computational Physics}, 225\penalty0 (1):\penalty0
  1066--1099, 2007.
\newblock ISSN 00219991.
\newblock \doi{10.1016/j.jcp.2007.01.017}.

\bibitem[Yu et~al.(2007)Yu, Zhou, and Wei]{Yu2007}
S~Yu, Y~Zhou, and G~Wei.
\newblock {Matched interface and boundary (MIB) method for elliptic problems
  with sharp-edged interfaces}.
\newblock \emph{Journal of Computational Physics}, 224\penalty0 (2):\penalty0
  729--756, 2007.
\newblock ISSN 00219991.
\newblock \doi{10.1016/j.jcp.2006.10.030}.
\newblock URL \url{http://dx.doi.org/10.1016/j.jcp.2006.10.030}.

\bibitem[Zhou et~al.(2006)Zhou, Zhao, Feig, and Wei]{Zhou2006}
Y.~C. Zhou, Shan Zhao, Michael Feig, and G.~W. Wei.
\newblock {High order matched interface and boundary method for elliptic
  equations with discontinuous coefficients and singular sources}.
\newblock \emph{Journal of Computational Physics}, 213\penalty0 (1):\penalty0
  1--30, 2006.
\newblock ISSN 00219991.
\newblock \doi{10.1016/j.jcp.2005.07.022}.

\bibitem[Gibou and Fedkiw(2005)]{Gibou2005}
Fr{\'{e}}d{\'{e}}ric Gibou and Ronald Fedkiw.
\newblock {A fourth order accurate discretization for the Laplace and heat
  equations on arbitrary domains, with applications to the Stefan problem}.
\newblock \emph{Journal of Computational Physics}, 202\penalty0 (2):\penalty0
  577--601, 2005.
\newblock ISSN 00219991.
\newblock \doi{10.1016/j.jcp.2004.07.018}.

\bibitem[Boyd(2005)]{Boyd2005}
John~P. Boyd.
\newblock Fourier embedded domain methods: extending a function defined on an
  irregular region to a rectangle so that the extension is spatially periodic
  and $\textnormal{C}^\infty$.
\newblock \emph{Applied Mathematics and Computation}, 161\penalty0
  (2):\penalty0 591--597, 2005.
\newblock ISSN 00963003.

\bibitem[Bueno-Orovio(2006)]{Bueno-Orovio2006}
Alfonso Bueno-Orovio.
\newblock Fourier embedded domain methods: Periodic and $c^\infty$ extension of
  a function defined on an irregular region to a rectangle via convolution with
  gaussian kernels.
\newblock \emph{Applied Mathematics and Computation}, 183\penalty0
  (2):\penalty0 813--818, 2006.
\newblock ISSN 00963003.
\newblock \doi{10.1016/j.amc.2006.06.029}.

\bibitem[Lui(2009)]{Lui2009}
S.~H. Lui.
\newblock {Spectral domain embedding for elliptic PDEs in complex domains}.
\newblock \emph{Journal of Computational and Applied Mathematics}, 225\penalty0
  (2):\penalty0 541--557, 2009.
\newblock ISSN 03770427.
\newblock \doi{10.1016/j.cam.2008.08.034}.
\newblock URL \url{http://dx.doi.org/10.1016/j.cam.2008.08.034}.

\bibitem[Sabetghadam et~al.(2009)Sabetghadam, Sharafatmandjoor, and
  Norouzi]{Sabetghadam2009}
Feriedoun Sabetghadam, Shervin Sharafatmandjoor, and Farhang Norouzi.
\newblock {Fourier spectral embedded boundary solution of the Poisson's and
  Laplace equations with Dirichlet boundary conditions}.
\newblock \emph{Journal of Computational Physics}, 228\penalty0 (1):\penalty0
  55--74, 2009.
\newblock ISSN 00219991.
\newblock \doi{10.1016/j.jcp.2008.08.018}.
\newblock URL \url{http://dx.doi.org/10.1016/j.jcp.2008.08.018}.

\bibitem[Albin and Bruno(2011)]{Albin2011}
Nathan Albin and Oscar~P. Bruno.
\newblock {A spectral FC Solver for the compressible Navier-Stokes equations in
  general domains I: Explicit time-stepping}.
\newblock \emph{Journal of Computational Physics}, 230\penalty0 (16):\penalty0
  6248--6270, jul 2011.
\newblock ISSN 00219991.
\newblock \doi{10.1016/j.jcp.2011.04.023}.
\newblock URL
  \url{http://linkinghub.elsevier.com/retrieve/pii/S0021999111002695}.

\bibitem[Lyon and Bruno(2010{\natexlab{a}})]{Lyon2010a}
Mark Lyon and Oscar~P. Bruno.
\newblock {High-order unconditionally stable FC-AD solvers for general smooth
  domains I. Basic Elements}.
\newblock \emph{Journal of Computational Physics}, 229\penalty0 (9):\penalty0
  3358--3381, 2010{\natexlab{a}}.
\newblock ISSN 00219991.
\newblock \doi{10.1016/j.jcp.2010.01.006}.
\newblock URL \url{http://dx.doi.org/10.1016/j.jcp.2009.11.020}.

\bibitem[Lyon and Bruno(2010{\natexlab{b}})]{Lyon2010}
Mark Lyon and Oscar~P. Bruno.
\newblock {High-order unconditionally stable FC-AD solvers for general smooth
  domains II. Elliptic, parabolic and hyperbolic PDEs; theoretical
  considerations}.
\newblock \emph{Journal of Computational Physics}, 229\penalty0 (9):\penalty0
  3358--3381, 2010{\natexlab{b}}.
\newblock ISSN 00219991.
\newblock \doi{10.1016/j.jcp.2010.01.006}.
\newblock URL \url{http://dx.doi.org/10.1016/j.jcp.2010.01.006}.

\bibitem[Shirokoff and Nave(2015)]{Shirokoff2013}
David Shirokoff and J-C Nave.
\newblock {A Sharp-Interface Active Penalty Method for the Incompressible
  Navier-Stokes Equations}.
\newblock \emph{Journal of Scientific Computing}, 62\penalty0 (1):\penalty0
  53----77, 2015.
\newblock URL \url{http://arxiv.org/abs/1303.5681}.

\bibitem[Stein et~al.(2015)Stein, Guy, and Thomases]{Stein2015}
David~B. Stein, Robert~D. Guy, and Becca Thomases.
\newblock {Immersed Boundary Smooth Extension: A high-order method for solving
  PDE on arbitrary smooth domains using Fourier spectral methods}.
\newblock \emph{Journal of Computational Physics}, 304:\penalty0 252--274,
  2015.
\newblock ISSN 10902716.
\newblock \doi{10.1016/j.jcp.2015.10.023}.
\newblock URL \url{http://arxiv.org/abs/1506.07561}.

\bibitem[Chorin(1968)]{chorin1968numerical}
Alexandre~Joel Chorin.
\newblock Numerical solution of the navier-stokes equations.
\newblock \emph{Mathematics of computation}, 22\penalty0 (104):\penalty0
  745--762, 1968.

\bibitem[Cai et~al.(2014)Cai, Nonaka, Bell, Griffith, and
  Donev]{cai2014efficient}
Mingchao Cai, Andy Nonaka, John~B Bell, Boyce~E Griffith, and Aleksandar Donev.
\newblock Efficient variable-coefficient finite-volume stokes solvers.
\newblock \emph{Communications in Computational Physics}, 16\penalty0
  (05):\penalty0 1263--1297, 2014.

\bibitem[Brinkman(1949)]{brinkman1949calculation}
HC~Brinkman.
\newblock A calculation of the viscous force exerted by a flowing fluid on a
  dense swarm of particles.
\newblock \emph{Applied Scientific Research}, 1\penalty0 (1):\penalty0 27--34,
  1949.

\bibitem[Adams and Fournier(2003)]{adams2003sobolev}
Robert~A Adams and John~JF Fournier.
\newblock \emph{Sobolev spaces}, volume 140.
\newblock Academic press, 2003.

\bibitem[Akiki and Balachandar(2016)]{Akiki2016}
G.~Akiki and S.~Balachandar.
\newblock {Immersed boundary method with non-uniform distribution of Lagrangian
  markers for a non-uniform Eulerian mesh}.
\newblock \emph{Journal of Computational Physics}, 307\penalty0
  (November):\penalty0 34--59, 2016.
\newblock ISSN 00219991.
\newblock \doi{10.1016/j.jcp.2015.11.019}.
\newblock URL
  \url{http://linkinghub.elsevier.com/retrieve/pii/S0021999115007597}.

\bibitem[Stein(2016)]{Stein2016}
David~B Stein.
\newblock \emph{{The Immersed Boundary Smooth Extension (IBSE) Method: A
  Flexible and Accurate Fictitious Domain Method, and Applications to the Study
  of Polymeric Flow in Complex Geometries}}.
\newblock PhD thesis, 2016.

\bibitem[John(1982)]{john1982partial}
Fritz John.
\newblock Partial differential equations, volume 1 of applied mathematical
  sciences, 1982.

\bibitem[Anderson et~al.(1999)Anderson, Bai, Bischof, Blackford, Demmel,
  Dongarra, Du~Croz, Greenbaum, Hammarling, McKenney, and Sorensen]{laug}
E.~Anderson, Z.~Bai, C.~Bischof, S.~Blackford, J.~Demmel, J.~Dongarra,
  J.~Du~Croz, A.~Greenbaum, S.~Hammarling, A.~McKenney, and D.~Sorensen.
\newblock \emph{{LAPACK} Users' Guide}.
\newblock Society for Industrial and Applied Mathematics, Philadelphia, PA,
  third edition, 1999.
\newblock ISBN 0-89871-447-8 (paperback).

\bibitem[Mori(2008)]{Mori2008}
Yoichiro Mori.
\newblock {Convergence proof of the velocity field for a stokes flow immersed
  boundary method}.
\newblock \emph{Communications on Pure and Applied Mathematics}, 61\penalty0
  (9):\penalty0 1213--1263, 2008.
\newblock ISSN 00103640.
\newblock \doi{10.1002/cpa.20233}.

\bibitem[Dou and Phan-thien(2006)]{Dou2006}
Hua-shu Dou and Nhan Phan-thien.
\newblock {The flow of an Oldroyd-B fluid past a cylinder in a channel :
  adaptive viscosity vorticity ( DAVSS- 3 ) formulation}.
\newblock 87\penalty0 (1999), 2006.

\bibitem[Alves et~al.(2001)Alves, Pinho, and Oliveira]{Alves2001}
M.a. Alves, F.T. Pinho, and P.J. Oliveira.
\newblock {The flow of viscoelastic fluids past a cylinder: finite-volume
  high-resolution methods}.
\newblock \emph{Journal of Non-Newtonian Fluid Mechanics}, 97\penalty0
  (2-3):\penalty0 207--232, feb 2001.
\newblock ISSN 03770257.
\newblock \doi{10.1016/S0377-0257(00)00198-1}.
\newblock URL
  \url{http://linkinghub.elsevier.com/retrieve/pii/S0377025700001981}.

\bibitem[Claus and Phillips(2013)]{Claus2013}
S.~Claus and T.N. Phillips.
\newblock {Viscoelastic flow around a confined cylinder using spectral/hp
  element methods}.
\newblock \emph{Journal of Non-Newtonian Fluid Mechanics}, 200:\penalty0
  131--146, oct 2013.
\newblock ISSN 03770257.
\newblock \doi{10.1016/j.jnnfm.2013.03.004}.
\newblock URL
  \url{http://linkinghub.elsevier.com/retrieve/pii/S0377025713000785}.

\bibitem[Welch et~al.(1965)Welch, Harlow, Shannon, and Daly]{welch1965mac}
J~Eddie Welch, Francis~Harvey Harlow, John~P Shannon, and Bart~J Daly.
\newblock The mac method-a computing technique for solving viscous,
  incompressible, transient fluid-flow problems involving free surfaces.
\newblock Technical report, Los Alamos Scientific Lab., Univ. of California, N.
  Mex., 1965.

\bibitem[McKee et~al.(2008)McKee, Tom{\'e}, Ferreira, Cuminato, Castelo, Sousa,
  and Mangiavacchi]{mckee2008mac}
S~McKee, MF~Tom{\'e}, VG~Ferreira, JA~Cuminato, A~Castelo, FS~Sousa, and
  N~Mangiavacchi.
\newblock {The MAC method}.
\newblock \emph{Computers \& Fluids}, 37\penalty0 (8):\penalty0 907--930, 2008.

\bibitem[Cai et~al.(2013)Cai, Nonaka, Bell, Griffith, and Donev]{Cai2013}
M.~Cai, A.~J. Nonaka, J.~B. Bell, B.~E. Griffith, and A.~Donev.
\newblock {Efficient Variable-Coefficient Finite-Volume Stokes Solvers}.
\newblock pages 1--25, 2013.
\newblock URL \url{http://arxiv.org/abs/1308.4605}.

\bibitem[Fan et~al.(1999)Fan, Tanner, and Phan-Thien]{fan1999galerkin}
Yurun Fan, RI~Tanner, and N~Phan-Thien.
\newblock Galerkin/least-square finite-element methods for steady viscoelastic
  flows.
\newblock \emph{Journal of Non-Newtonian Fluid Mechanics}, 84\penalty0
  (2):\penalty0 233--256, 1999.

\bibitem[Hundsdorfer and Ruuth(2007)]{Hundsdorfer2007}
Willem Hundsdorfer and Steven~J. Ruuth.
\newblock {IMEX extensions of linear multistep methods with general
  monotonicity and boundedness properties}.
\newblock \emph{Journal of Computational Physics}, 225\penalty0 (2):\penalty0
  2016--2042, 2007.
\newblock ISSN 00219991.
\newblock \doi{10.1016/j.jcp.2007.03.003}.

\bibitem[Sch{\"{a}}fer and Turek(1996)]{Schafer1996}
M.~Sch{\"{a}}fer and S.~Turek.
\newblock {Benchmark Computations of Laminar Flow Around a Cylinder}.
\newblock \emph{Flow Simulation with High-Performance Computers II, Volume 52
  of Notes on Numerical Fluid Mechanics, Vieweg}, 52:\penalty0 547 -- 566,
  1996.
\newblock ISSN 0001-1452.
\newblock \doi{10.1007/978-3-322-89849-4_39}.

\bibitem[Nabh(1998)]{nabh1998high}
Guido Nabh.
\newblock \emph{On high order methods for the stationary incompressible
  Navier-Stokes equations}.
\newblock Interdisziplin{\"a}res Zentrum f{\"u}r Wiss. Rechnen der Univ.
  Heidelberg, 1998.

\bibitem[Platte et~al.(2011)Platte, Trefethen, and Kuijlaars]{Platte2011}
Rodrigo~B. Platte, Lloyd~N. Trefethen, and Arno B.~J. Kuijlaars.
\newblock {Impossibility of Fast Stable Approximation of Analytic Functions
  from Equispaced Samples}.
\newblock \emph{SIAM Review}, 53\penalty0 (2):\penalty0 308--318, 2011.
\newblock ISSN 0036-1445.
\newblock \doi{10.1137/090774707}.

\bibitem[Griffith and Luo(2012)]{griffith2012hybrid}
Boyce~E Griffith and Xiaoyu Luo.
\newblock Hybrid finite difference/finite element version of the immersed
  boundary method.
\newblock \emph{Submitted in revised form}, 2012.

\bibitem[Shankar et~al.(2012)Shankar, Wright, Fogelson, and Kirby]{Shankar2012}
Varun Shankar, Grady~B. Wright, Aaron~L. Fogelson, and R.~M. Kirby.
\newblock {A Study of Different Modeling Choices For Simulating Platelets
  Within the Immersed Boundary Method}.
\newblock pages 1--33, oct 2012.
\newblock URL \url{http://arxiv.org/abs/1210.1885v1}.

\end{thebibliography}
\bibliographystyle{unsrtnat}





\end{document}